\documentclass[11 pt,a4 paper, oneside]{amsart}
\usepackage{latexsym}
\usepackage[top=3cm, bottom=3.5cm]{geometry}
\usepackage[utf8]{inputenc}

\usepackage[T1]{fontenc}

\usepackage[english]{babel}
\date{}

\usepackage{bm}
\usepackage{amsthm}
\usepackage{amssymb}
\usepackage{amsmath}
\usepackage{amsrefs}
\usepackage{indentfirst}
\usepackage{graphicx}
\usepackage{mathrsfs}

\usepackage{multicol}
\usepackage{multirow}
\usepackage{lscape}

\usepackage{marginnote}
\usepackage{enumitem}

\usepackage[hidelinks]{hyperref}
\hypersetup{
	colorlinks,
	linkcolor={black},
	citecolor={blue},
	urlcolor={black}
}

\usepackage{caption}

\usepackage{changepage}

\usepackage{cleveref}

\usepackage{rotating}
\usepackage{float}
\usepackage{flafter} 
\usepackage{tabularx}

\usepackage{tikz} 
\usepackage{tikz-cd}
\usepackage{amsfonts}
\newtheorem{theorem}{Theorem}[section]
\newtheorem{lemma}[theorem]{Lemma}
\newtheorem{proposition}[theorem]{Proposition}
\newtheorem{corollary}[theorem]{Corollary}
\theoremstyle{definition}
\newtheorem{definition}[theorem]{Definition}
\newtheorem{remark}[theorem]{Remark}
\newtheorem{result}{Main result}
 
\newtheorem{thmintro}[result]{Theorem}

\title{Some computations on trivial canonical-bundle solvmanifolds}
\author{Lapo Rubini}

\newcommand{\olsi}[1]{\,\overline{\!{#1}}} 
\newcommand{\barT}{\olsi{T}}
\newcommand{\uno}{\bar{1}}
\newcommand{\due}{\bar{2}}
\newcommand{\tre}{\bar{3}}
\DeclareMathOperator{\Ima}{Im}

\newcommand{\g}{\mathfrak{g}}
\newcommand{\C}{\mathbb{C}}
\newcommand{\R}{\mathbb{R}}
\newcommand{\Z}{\mathbb{Z}}

\begin{document}

\begin{abstract}
	We compute the Dolbeault and the Bott-Chern cohomology of six dimensional solvmanifolds endowed with a complex structure of splitting type, introduced by Kasuya, and with trivial canonical bundle. We build, following results by Angella and Kasuya, finite dimensional double subcomplexes $(C_\Gamma^{\bullet,\bullet},\partial,\bar{\partial})\subseteq(\wedge^{\bullet,\bullet}G/\Gamma,\partial,\bar{\partial})$ for which the inclusion is an isomorphism in cohomology. We decompose such double complexes into indecomposable ones. Lastly, we study some notions of formality for this class of manifolds, giving a characterization of the $\partial\bar{\partial}$-Lemma property in general complex dimension, and we compute triple ABC-Massey products on them.
\end{abstract}

\pagestyle{plain}

\maketitle
\tableofcontents

\section*{Introduction}

The class of nilmanifolds and solvmanifolds, i.e. compact quotients of connected and simply-connected nilpotent, respectively solvable, Lie groups by a lattice, provides plenty of examples in non-K\"ahler geometry. Thanks to their structure, some of their properties and invariants can be deduced from the associated Lie algebras (\cite{Nomizu, Rollenske2010dolbeaultcohomologynilmanifoldsleftinvariant, Hattori_60, Sakane1976OnCC, console2001dolbeaultcohomologycompactnilmanifolds}). Most of the results regarding the nilpotent case involve the presence of a left-invariant complex structure on the manifold $G/\Gamma$, that is a complex structure induced by a complex structure on the Lie algebra $\g$ associated to $G$. A nilmanifold carrying such a complex structure has, by \cite{SALAMON2001311}, holomorphically trivial canonical bundle (i.e. admits an invariant non-zero closed $(n,0)$-form, where $n$ is the complex dimension of the manifold). This property is interesting in itself, as it can be considered as a generalization of the Calabi-Yau property in the non-K\"ahler setting (\cite{Tosatti_2015}).
	
In the solvable case, the presence of an invariant complex structure does not guarantee the triviality of the canonical bundle. This has led to various results regarding the classification of solvmanifolds with or without this property (\cite{AngellaOtalUgarte_2017,Bock_2016,fino2015dimensionalsolvmanifoldsholomorphicallytrivial,tolcachier2024sixdimensionalcomplexsolvmanifoldsnoninvariant}). Among these manifolds Kasuya introduced in \cite{kasuya2012techniquescomputationsdolbeaultcohomology} the subclass of solvmanifolds endowed with a structure of splitting type, i.e. the Lie group defining the manifold can be expressed as a semidirect product of $\C^n$ and a nilpotent Lie group with certain compatibility properties with the complex structure. On these manifolds, Kasuya and Angella in \cite{kasuya2012techniquescomputationsdolbeaultcohomology,AngellaKasuya_2017} have constructed a finite-dimensional double complex that allows the explicit computation of the $\partial$ and $\bar{\partial}$ Dolbeault and Bott-Chern cohomologies of such manifolds (and hence, thanks to  \cite[Theorem C]{stelzig2023pluripotentialhomotopytheory}, also Aeppli cohomology). Here, by a double complex (or bicomplex) $(A,\partial,\bar{\partial})$ over the field $K$ we mean a bigraded $K$-vector space $A=\bigoplus_{p,q\in\Z} A^{p,q}$ with two endomorphisms $\partial, \bar{\partial}$ of bidegree respectively $ (1,0) $ and $(0,1)$ such that $\partial\circ\partial=0$, $\bar{\partial}\circ\bar{\partial}=0$, $\partial\circ\bar{\partial}+\bar{\partial}\circ\partial=0$ (see \cite{Stelzig_2021} for more details). If $(A,\partial,\bar{\partial})$ is a bicomplex, Dolbeault, Bott-Chern and Aeppli cohomologies, introduced in \cite{dolbeault1956}, \cite{BottChern_65} and \cite{Aeppli_65} respectively, are defined as
\begin{equation*}
	H^{\bullet,\bullet}_{\bar{\partial}}(A):=\dfrac{\ker\bar{\partial}}{\Ima\bar{\partial}},\quad
	H^{\bullet,\bullet}_{BC}(A):=\dfrac{\ker\partial\cap\ker\bar{\partial}}{\Ima\partial\bar{\partial}}\quad\text{ and }\quad H^{\bullet,\bullet}_{A}(A):=\dfrac{\ker\partial\bar{\partial}}{\Ima\partial+\Ima\bar{\partial}}.
\end{equation*}
If $X$ is a complex manifold, its cohomologies are defined as the cohomologies of the double complex of its differential forms $(\Lambda^{\bullet,\bullet}X,\partial,\bar{\partial})$. The interest in Bott-Chern cohomology is justified by many results (\cite{Angella_2015,AngellaKasuya_2017,Angella_2012,DeligneGriffiths1975,Stelzig_2021,schweitzer2007autourlacohomologiebottchern,stelzig2023pluripotentialhomotopytheory} and many others), among which we find its relation with the $\partial\bar{\partial}$-Lemma property, introduced in \cite{DeligneGriffiths1975}: namely a compact complex manifold satisfies the $\partial\bar{\partial}$-Lemma if a complex form that is $\partial$-closed, $\bar{\partial}$-closed and $d$-exact is $\partial\bar{\partial}$-exact. This property holds if and only if Bott-Chern and Dolbeault cohomologies are naturally isomorphic. 

Since a compact K\"ahler manifold satisfies the $\partial\bar{\partial}$-Lemma, Bott-Chern cohomology is particularly studied in non-K\"ahler geometry. Moreover, some notions of formality regarding Bott-Chern cohomology has recently been introduced and studied in \cite{Angella_2015,milivojevic2024bigradednotionsformalityaepplibottchernmassey,placini2024nontrivialmasseyproductscompact,stelzig2023pluripotentialhomotopytheory,tardini2015geometricbottchernformalitydeformations}, giving new invariants which can be computed starting from the bigraded bidifferential algebra $\wedge^{\bullet,\bullet}X$ of differential forms, such as Aeppli-Bott-Chern-Massey products, introduced in \cite{Angella_2015}. Among the various notions of formality we can find strong formality, defined in \cite[Definition 4.1]{milivojevic2024bigradednotionsformalityaepplibottchernmassey}: a manifold is said to be strongly formal if its algebra of differential forms is equivalent (in a certain sense) to a trivial bigraded algebra $(H,\partial=0,\bar{\partial}=0)$ (see Section \ref{section4} for a precise definition). In particular, it has been proved in \cite{placini2024nontrivialmasseyproductscompact} and in \cite{sferruzza2024bottchernformalitymasseyproducts} that strong formality is a property independent of the K\"ahlerianity of the manifold and of the $\partial\bar{\partial}$-Lemma, meaning that it is not implied by these conditions.\\ 

This paper has the aim to give explicit cohomology computations on six dimensional solvmanifolds carrying a complex structure of splitting type and with trivial canonical bundle and to study their formality. As recalled before, thanks to \cite{kasuya2012techniquescomputationsdolbeaultcohomology} and \cite{AngellaKasuya_2017}, Dolbeault and Bott-Chern cohomologies can be explicitly computed, passing through the definition of a finite dimensional sub-bicomplex of $\wedge^{\bullet,\bullet}X$ built starting from the structure of the solvmanifold. Indeed there exists a finite dimensional subcomplex $C^{\bullet,\bullet}_\Gamma\subseteq\wedge^{\bullet,\bullet}X$ such that the inclusion induces isomorphisms in Dolbeault (and its conjugate) and Bott-Chern cohomologies, hence thanks to \cite[Theorem C]{stelzig2023pluripotentialhomotopytheory} it can be used to compute any cohomological invariant (e.g. de Rham cohomology, Aeppli cohomology...). Thanks to the classifications in \cite{AngellaOtalUgarte_2017} of splitting type complex structures on six-dimensional solvable Lie algebras and in \cite{fino2015dimensionalsolvmanifoldsholomorphicallytrivial,tolcachier2024sixdimensionalcomplexsolvmanifoldsnoninvariant} of six dimensional solvmanifolds with holomorphically trivial canonical bundle, we have selected the algebras $\g_1,\g_2^{\alpha}$ and $\g_8$ and computed their cohomology by applying the results in \cite{kasuya2012techniquescomputationsdolbeaultcohomology, AngellaKasuya_2017}, completing the work started by Angella, Kasuya and Otal in \cite{kasuya2012techniquescomputationsdolbeaultcohomology,AngellaKasuya_2017,otalthesis}. We can summarize our computations as follow. 

\begin{result}
	The splitting type solvmanifolds with underlying Lie algebras $\g_1$, $\g_2^{\alpha}$ and $\g_8$ have respectively $3$, $5$ and $7$ possible cohomologies. Moreover, there are $3$ recurrent cohomology structures (see Section \ref{section2}), meaning that they are realized by all these three algebras.
\end{result}

The Dolbeault and the Bott-Chern cohomologies, together with the decompositions of the bicomplexes $C_\Gamma^{\bullet,\bullet}$, are summarized in \crefrange{table1}{tabellaBCg8Otal} in Appendix \ref{tables} and in \crefrange{figure1}{figure15} in Appendix \ref{decompositions}.

With the bicomplex $C^{\bullet,\bullet}_{\Gamma}$ defined by Angella and Kasuya we have also studied various notions of bigraded formality for splitting type solvmanifolds with holomorphically trivial canonical bundle. In addition to strong formality and to the $\partial\bar{\partial}$-Lemma property, we have analysed weak formality and geometrically-Bott-Chern formality, introduced in \cite{milivojevic2024bigradednotionsformalityaepplibottchernmassey} and in \cite{Angella_2015} respectively: namely, weak formality refers to bigraded algebras cohomologically equivalent to an algebra such that $\partial\bar{\partial}=0$, geometrically-Bott-Chern formality to manifolds whose Bott-Chern cohomology ring is indeed a bigraded differential algebra (precise definitions are given in Section \ref{section4}). We have proved the following result, concerning both complex dimension $3$ and general dimension. 
\begin{thmintro}\label{thmintro}
	For splitting type solvmanifolds of the form $\C^n\ltimes_\varphi\C^m/\Gamma$, the $\partial\bar{\partial}$-Lemma property and strong formality are equivalent notions. For complex dimension $3$ splitting type solvmanifolds with holomorphically trivial canonical bundle, these notions are also equivalent to geometrically-Bott-Chern-formality and to weak formality.
\end{thmintro}

In Section \ref{section1} we recall the definitions and the classifications of 6 dimensional splitting type solvmanifolds with holomorphically trivial canonical bundle. We recall the notion of splitting type complex structure defined by Kasuya in \cite{kasuya2012techniquescomputationsdolbeaultcohomology} and, thanks to \cite{fino2015dimensionalsolvmanifoldsholomorphicallytrivial,tolcachier2024sixdimensionalcomplexsolvmanifoldsnoninvariant} and \cite{AngellaOtalUgarte_2017}, we focus on the Lie algebras $\g_1, \g_2^{\alpha}$ and $\g_8$. In addition to the lattices for the corresponding connected and simply-connected Lie groups considered in \cite[Proposition 2.10]{fino2015dimensionalsolvmanifoldsholomorphicallytrivial}, we give a countable family of lattices for the Lie group $G_1$ underlying the Lie algebra $ \g_1 $.

Section \ref{section2} presents computations made in \cite{kasuya2012techniquescomputationsdolbeaultcohomology,AngellaKasuya_2017,otalthesis} together with ours, passing through the definition of the double complex $C^{\bullet,\bullet}_\Gamma$ and its explicit expression in all the cases involved. We list in Appendix \ref{tables} such bicomplexes and the cohomologies obtained.

In Section \ref{section3} we analyse the structure of the double complexes built in the previous section. We recall the definition of the indecomposable double complexes and \cite[Theorem 3]{Stelzig_2021}, that states the uniqueness (up to isomorphisms) of the decomposition of a bicomplex into indecomposable ones. We also recall some useful results in cohomology, proved by Stelzig in \cite{Stelzig_2021}, that involve such decomposition. We compute such decomposition for the bicomplexes $C_\Gamma^{\bullet,\bullet}$, giving explicit expressions in Appendix \ref{decompositions}. 

In Section \ref{section4} we study bigraded notions of formality on the class of manifolds introduced in Section \ref{section1}. We recall the definitions of geometrically-Bott-Chern-formality and ABC-Massey products introduced in \cite{Angella_2015}, as well as weak and strong formality defined in \cite{milivojevic2024bigradednotionsformalityaepplibottchernmassey}. We then proceed to prove Theorem \ref{thmintro}: we show in Corollary \ref{corollarymio} the equivalence between $\partial\bar{\partial}$-Lemma property and strong formality for splitting type solvmanifolds of the form $\C^n\ltimes_\varphi\C^m/\Gamma$ and, lastly, we compute triple ABC-Massey products for the manifolds considered in the previous sections, verifying that in complex dimension 3 a splitting type solvmanifold with holomorphically trivial canonical bundle is geometrically-Bott-Chern-formal if and only if it satisfies the $\partial\bar{\partial}$-Lemma. 

\subsubsection*{Acknowledgements:} This work has been supported by the project PRIN 2022 "Real and Complex Manifolds: Geometry and holomorphic dynamics" (code
2022AP8HZ9). The author would like to express his sincere gratitude to Daniele Angella and Jonas Stelzig for their continuous support and valuable suggestions.

\section{\for{toc}{Solvmanifolds with trivial canonical bundle}\except{toc}{Solvmanifolds with trivial canonical bundle and their invariant complex structures}}\label{section1}

In this section we recall some results regarding the class of manifolds analysed in this paper. For more details, see for example \cite{Bock_2016}, \cite{AngellaOtalUgarte_2017}, \cite{otalthesis}, \cite{fino2015dimensionalsolvmanifoldsholomorphicallytrivial}, \cite{tolcachier2024sixdimensionalcomplexsolvmanifoldsnoninvariant}.

We consider solvmanifolds, namely manifolds defined as compact quotients $G/\Gamma$, where $G$ is a solvable Lie group and $\Gamma\subseteq G$ is a lattice. A (left-)invariant complex structure on a Lie group $G$ is a complex structure $J$ that is invariant with respect to the left multiplication on $G$. This condition is equivalent to a linear complex structure on the Lie algebra $\g$ associated to $G$. Moreover, an invariant complex structure on a Lie group $G$ naturally induces an invariant complex structure on any quotient $G/\Gamma$.

Six-dimensional solvmanifolds with an invariant complex structure and trivial canonical bundle have been classified by Fino, Otal and Ugarte in \cite[Theorem 2.8]{fino2015dimensionalsolvmanifoldsholomorphicallytrivial} and by Tolcachier in \cite[Theorem 4.2]{tolcachier2024sixdimensionalcomplexsolvmanifoldsnoninvariant}. They include nilmanifolds (that have trivial canonical bundle by \cite{barberis2007}), classified in \cite{ugarte2007}, and the following.

\begin{theorem}\cite[Theorem 2.8]{fino2015dimensionalsolvmanifoldsholomorphicallytrivial}, \cite[Theorem 4.2]{tolcachier2024sixdimensionalcomplexsolvmanifoldsnoninvariant}\label{classificazione}
	Let $\g$ be a strongly unimodular (non-nilpotent) solvable Lie algebra of dimension $6$. Then $\g$ admits a complex structure with non-zero closed $(3,0)$-form (i.e. with holomorphically trivial canonical bundle) if and only if it is isomorphic to one of the following:
	
	\begin{list}{}{\setlength{\itemsep}{.5ex}}
		\item $\g_1=A^{-1,-1,1}_{5,7}\oplus\mathbb{R}=(e^{15}, -e^{25}, -e^{35}, e^{45},0,0)$;
		\item $\g_2^\alpha=A^{\alpha,-\alpha,1}_{5,17}\oplus\mathbb{R}=(\alpha e^{15}+e^{25},-e^{15}+\alpha e^{25},-\alpha e^{35}+e^{45},-e^{35}-\alpha e^{45},0,0), \,\,\alpha\ge0$;
		\item $\g_3=\mathfrak{e}(2)\oplus\mathfrak{e}(1,1)=(0,-e^{13}, e^{12}, 0, -e^{46}, -e^{45})$;
		\item $\g_4=A^{0,0,1}_{6,37}=(e^{23}, -e^{36}, e^{26}, -e^{56}, e^{46}, 0)$;
		\item $\g_5=A^{0,1,1}_{6,82}=(e^{24}+e^{35}, e^{26}, e^{36}, -e^{46}, -e^{56}, 0)$;
		\item $\g_6=A^{0,0,1}_{6,88}=(e^{24}+e^{35}, -e^{36}, e^{26}, -e^{56}, e^{46}, 0)$;
		\item $\g_7=B^1_{6,6}=(e^{24}+e^{35}, e^{46}, e^{56}, -e^{26}, -e^{36}, 0)$;
		\item $\g_8=N_{6,18}^{0,-1,-1}=(e^{16}-e^{25}, e^{15}+e^{26}, -e^{36}+e^{45}, -e^{35}-e^{46}, 0, 0)$;
		\item $\g_9=B^1_{6,4}=(e^{45}, e^{15}+e^{36}, e^{14}-e^{26}+e^{56},-e^{56}, e^{46}, 0)$;
		\item $\g_{10}=\mathfrak{s}_{6,147}^0=(e^{26}-e^{35}, -e^{16}+e^{36}-e^{45}, e^{46},-e^{36}, 0, 0)$.
	\end{list}
\end{theorem}

We follow the notations introduced in \cite{Bock_2016} and in \cite{snobl_winternitz2014} specifically for the listed algebras and in \cite{SALAMON2001311} for the general expressions. For example, we use the notation $(e^{15}, -e^{25}, -e^{35}, e^{45},0,0)$ to mean that there exists a basis of the dual of the Lie algebra  $\{e^i\}_{i=1}^6$ satisfying $de^1=e^1\wedge e^5,$ $ de^2=-e^2\wedge e^5$, $de^3=-e^3\wedge e^5$, $de^4=e^4\wedge e^5$ and $de^5=de^6=0$.

Thanks to \cite{AngellaKasuya_2017}, \cite{AngellaOtalUgarte_2017} and \cite{kasuya2012techniquescomputationsdolbeaultcohomology}, we can straightforwardly compute Dolbeault and Bott-Chern cohomology of splitting type solvmanifolds. Hence, we focus on this class of manifolds, defined as follow.

\begin{definition}\cite[Assumption 1.1]{kasuya2012techniquescomputationsdolbeaultcohomology}
	A solvmanifold $M=G/\Gamma$ endowed with an invariant complex structure $J$ is said to be of splitting type if $G$ is a semi-direct product $\mathbb{C}^n\ltimes_\varphi N$, where:
	\begin{enumerate}
		\item $N$ is a connected simply-connected $2m$-dimensional nilpotent Lie group endowed with a left-invariant complex structure $J_N$;
		
		\item for any $z\in\mathbb{C}^n$, the automorphism $\varphi(z)\in Aut(N)$ is a holomorphic automorphism of $N$ with respect to $J_N$;
		
		\item $\varphi$ induces a semi-simple action on the Lie algebra $\mathfrak{n}$ associated to $N$;
		
		\item the lattice $\Gamma$ can be written as $\Gamma=\Gamma_{\mathbb{C}^n}\ltimes_\varphi\Gamma_N$, where $\Gamma_{\mathbb{C}^n} $ and $\Gamma_N$ are lattices of $\mathbb{C}^n$ and $N$ respectively and it holds $\varphi(z)(\Gamma_N)\subseteq\Gamma_N$ for any $z\in\Gamma_{\mathbb{C}^n}$;
		
		\item the inclusion 		$({\wedge}^{\bullet,\bullet}\mathfrak{n}^*,\bar{\partial})\hookrightarrow({\wedge}^{\bullet,\bullet}N/\Gamma_N,\bar{\partial})$ induces an isomorphism in cohomology
		\begin{equation*}
		H^{\bullet,\bullet}_{\bar{\partial}}(\mathfrak{n})\simeq H^{\bullet,\bullet}_{\bar{\partial}}(N/\Gamma_N).
		\end{equation*} 
	\end{enumerate}
\end{definition}

If we compare Theorem \ref{classificazione} with the classification made in \cite[Theorem 1.7]{AngellaOtalUgarte_2017}, we see that $\g_1, \g_2^{\alpha} $ and $\g_8$ are the only Lie algebras that carry a structure of splitting type with holomorphically trivial canonical bundle.

The complex structures we are interested in are classified, up to equivalence, in \cite[Proposition 3.3, Proposition 3.4]{fino2015dimensionalsolvmanifoldsholomorphicallytrivial} (we recall that two complex structures $J$ and $J'$ on $\g$ are said to be equivalent if there exists an automorphism $F$ of $\g$ such that $F\circ J=J'\circ F$). Since we work with differential forms, it is useful to write the structure equations obtained from these complex structures, which define the images of the differentials $\partial$ and $\bar{\partial}$ on a $(1,0)$-basis $\{\omega^1,\omega^2,\omega^3\}$ of $(\g^{1,0})^*$:
\begin{align}
&(\mathfrak{g}_1,J):\begin{cases}\label{str1}
d\omega^1=\omega^1\wedge(\omega^3+\omega^{\tre})\\
d\omega^2=-\omega^2\wedge(\omega^3+\omega^{\tre})\\
d\omega^3=0
\end{cases};\\
&(\mathfrak{g}_2^0,J):\begin{cases}\label{str20}
d\omega^1=i\omega^1\wedge(\omega^3+\omega^{\tre})\\
d\omega^2=-i\omega^2\wedge(\omega^3+\omega^{\tre})\\
d\omega^3=0
\end{cases};\\
&(\mathfrak{g}_2^{\alpha=\frac{\cos\theta}{\sin\theta}},J^\pm):\begin{cases}\label{str2a}
d\omega^1=(\pm\cos\theta+i\sin\theta)\omega^1\wedge(\omega^3+\omega^{\tre})\\
d\omega^2=-(\pm\cos\theta+i\sin\theta)\omega^2\wedge(\omega^3+\omega^{\tre})\\
d\omega^3=0
\end{cases}\text{where }\theta\in(0,\pi/2);\\
&(\g_8, J^A):\begin{cases}\label{struttura8}
d\omega^1=-(A-i)\omega^{13}-(A+i)\omega^{1\tre}\\
d\omega^2=(A-i)\omega^{23}+(A+i)\omega^{2\tre}\\
d\omega^3=0
\end{cases}
\end{align}
where $A\in\mathbb{C}$ with $\mathfrak{Im} (A)\neq0$ (here we denote $\bar{\omega}^i$ by $\omega^{\bar{i}}$ and the wedge products $\omega^i\wedge\omega^j$ by $\omega^{ij}$). We observe that these complex structures allow us to reduce our study to splitting type solvmanifolds, as justified by the following result.

\begin{proposition}\cite[Proposition 4.2.15]{otalthesis}\label{PropCaratteri}
	Let $G$ be a Lie group with underlying algebra $\g_1$, $\g_2^{\alpha}$ with $\alpha\ge0$ or $\g_8$ endowed with a left-invariant complex structure $J$ defined by the equations \eqref{str1}, \eqref{str20}, \eqref{str2a} or \eqref{struttura8}. Then $G$ is a semi-direct product $G=\mathbb{C}\ltimes_{\varphi_A}\mathbb{C}^2$, where $\varphi_A:\mathbb{C}\to\emph{GL}(\mathbb{C}^2)$ is defined by the diagonal matrix:
	\begin{equation}\label{action}
	\varphi_A(z):=\left(\begin{array}{cc}
	\alpha_1(z) & 0\\
	0 & \alpha_2(z)\\
	\end{array}
	\right)
	\end{equation}
	with $\alpha_1,\alpha_2:\mathbb{C}\to\mathbb{C}^*$ characters such that $\alpha_2=\alpha_1^{-1}$ and defined by
	\begin{equation}
	\alpha_1(z):=\begin{cases}
	e^{A(z+\bar{z})}&\text{if $J$ satisfies \eqref{str1}, \eqref{str20} or \eqref{str2a}}\\
	e^{-(A-i)z-(A+i)\bar{z}}\quad&\text{if $J$ satisfies \eqref{struttura8}}
	\end{cases}
	\end{equation}	
	with $A=1$ in case \eqref{str1}, $A=i$ in case \eqref{str20}, $A=\pm\cos\theta+i\sin\theta$, $\theta\in(0,\frac{\pi}{2})$ in case \eqref{str2a}.
\end{proposition}

Moreover, we observe that thanks to this result we can focus on splitting type solvmanifolds defined from a Lie group $G$ that is the semi-direct product of two complex vector spaces, $G=\mathbb{C}\ltimes_{\varphi_A}\mathbb{C}^2$.

\begin{remark}
	The complex structure on $\g_1$ given by \eqref{str1} and the complex structure on $\g_8$ given by \eqref{struttura8} for $A=-i$ correspond to the complex structures studied by Kasuya in \cite[Example 1, Example 2 case (a)]{kasuya2012techniquescomputationsdolbeaultcohomology}, called respectively the completely solvable Nakamura manifold and the complex parallelizable Nakamura manifold.
\end{remark}

The next step consists in finding lattices for the connected and simply-connected Lie groups $G_1, G_2^{\alpha}$ and $G_8$ associated to the Lie algebras $\g_1$, $\g_2^{\alpha}$ and $\g_8$. The next result guarantees the existence of such lattices.

\begin{proposition}\cite[Proposition 2.10]{fino2015dimensionalsolvmanifoldsholomorphicallytrivial}
	The connected and simply-connected Lie groups $G_1$ and $G_8$ with underlying Lie algebras $\mathfrak{g}_1$ and $\g_8$ respectively admit lattices.
	
	Moreover, there exists a countable number of distinct $\alpha$'s (including $\alpha=0$) for which the connected, simply-connected Lie group $G_2^{\alpha}$ underlying the Lie algebra $\mathfrak{g}_2^{\alpha}$ admits a lattice. 
\end{proposition}

Using the same techniques used to prove this Proposition, we also have found a countable family of lattices $\Gamma_n$ for the connected and simply-connected Lie group with underlying algebra $\g_1 \simeq A^{-1,-1,1}_{5,7}\oplus\R$ following the technique used for $\g_2^{\alpha}$ in the proof of \cite[Proposition 2.10]{fino2015dimensionalsolvmanifoldsholomorphicallytrivial}. We can write the nilpotent ideal of the algebra as the semi-direct product $A^{-1,-1,1}_{5,7}=\mathbb{R}\ltimes_{\text{ad}_{e_5}}\mathbb{R}^4$, with coordinate matrix of $\text{ad}_{e_5}$ given by $ B=\text{diag}(-1,1,1,-1) $. Then, the coordinate matrix of $d_e(\mu(t))$ is $\text{diag}(e^{-t}, e^{t},e^{t},e^{-t})$. Its characteristic polynomial $p(\lambda)=(1-(e^t+e^{-t})\lambda+\lambda^2)^2$ is integer if $t=t_n=\log(\frac{n+\sqrt{n^2-4}}{2})$ with $n\in\mathbb{Z}$. Moreover, for $n\neq\pm2$, we have $e^{t_nB}=M^{-1}C_nM$, where
	\begin{equation*}
		M^{-1}=\left(\begin{array}{cccc}
		1 & 0 & a_n^{-1} & 0\\
		0 & 1 & 0 & a_n\\
		1 & 0 & a_n & 0 \\
		0 & 1 & 0 & a_n^{-1}\\
		\end{array}\right),\quad 
		C_n=\left(\begin{array}{cccc}
		n & 0 & -1 & 0\\
		0 & n & 0 & -1\\
		1 & 0 & 0 & 0 \\
		0 & 1 & 0 & 0\\
		\end{array}\right),	
	\end{equation*}
with $a_n=-\frac{n-\sqrt{n^2-4}}{2}$. Hence, if we take the basis 
\begin{equation*}
X_1=e_1+e_3,\, X_2=e_2+e_4,\, X_3=a_n^{-1}e_1+a_ne_3,\, X_4=a_ne_2+a_n^{-1}e_4,
\end{equation*}
applying \cite[Theorem 2.4]{Yamada_16} we have a lattice 
\begin{equation}\label{lattice}
\Gamma_n =\log\left(\frac{n+\sqrt{n^2-4}}{2}\right)\mathbb{Z}\ltimes_{\mu}\exp^{\mathbb{R}^4}(\mathbb{Z}\langle X_1,X_2,X_3,X_4\rangle)
\end{equation}
for each $n\in\mathbb{Z}$, $n\neq\pm2$. Hence, $\Gamma=\Gamma_n\times\mathbb{Z}$ is a lattice for the Lie group $G_1$ with underlying Lie algebra $\g_1$.

\section{\for{toc}{Dolbeault and Bott-Chern cohomology of solvmanifolds}\except{toc}{Dolbeault and Bott-Chern cohomology of splitting type solvmanifolds with holomorphically trivial canonical bundle}}\label{section2}

The class of splitting type solvmanifolds has been introduced by Kasuya in \cite{kasuya2012techniquescomputationsdolbeaultcohomology} to have a class of manifolds on which compute explicitly the cohomology. Indeed, for a solvmanifold $G/\Gamma$ of splitting type it is possible to build a finite-dimensional double complex that has the same Dolbeault cohomology of the manifold (see \cite{kasuya2012techniquescomputationsdolbeaultcohomology}) starting from the Lie algebra $\g$ underlying $G$ and the lattice $\Gamma$ defining the solvmanifold. Moreover, Angella and Kasuya in \cite{AngellaKasuya_2017} show an analogue result for the computation of Bott-Chern cohomology.

More precisely, let $G=\mathbb{C}^n\ltimes_\varphi N$ be a Lie group endowed with an invariant complex structure of splitting type. Let $\{X_1,...,X_n\}$ be the standard basis of $\C^n$ as a complex vector space and $\{Y_1,...,Y_m\}\subseteq\mathfrak{n}^{1,0}$ an invariant $(1,0)$-basis for the complex structure $J_N$ of $N$ such that the induced action $\varphi$ on $\mathfrak{n}^{1,0}$ is represented in this basis by the diagonal matrix
	\begin{equation*}
		\varphi(z)=\left(\begin{array}{ccc}
			\alpha_1&&\\
			&\ddots&\\
			&&\alpha_m
		\end{array}\right)
	\end{equation*} 
for $\alpha_1,\dots...,\alpha_m:\C^n\to\C^*$ characters. Let $\{x_1,..,x_n,\alpha_1^{-1}y_1,...,\alpha_m^{-1}y_m\}$ be a basis of $(\g^{1,0})^*$, dual to the basis $\{X_1,...,X_n,\alpha_1Y_1,...,\alpha_mY_m\}$ of $\g^{1,0}$. 

Using \cite[Lemma 2.2]{kasuya2012techniquescomputationsdolbeaultcohomology} we find unique unitary characters $\beta_i,\gamma_j:\C^n\to\C^*$ such that the products $\alpha_i\beta_i^{-1}$ and $\bar{\alpha}_j\gamma_j^{-1}$ are holomorphic for $i,j=1,..,m$. Starting from this construction, Kasuya proved the following result.

\begin{theorem}\cite[Corollary 4.2]{kasuya2012techniquescomputationsdolbeaultcohomology} \label{thmKasuya}
	Let $G/\Gamma$ be a solvmanifold endowed with an invariant complex structure of splitting type and $B^{\bullet,\bullet}_{\Gamma}$ the finite-dimensional double complex given by 
	\begin{align}
	B^{p,q}_\Gamma:=\mathbb{C}\langle x_I\wedge(\alpha_J^{-1}\beta_J)y_J\wedge\bar{x}_K\wedge(\bar{\alpha}_L^{-1}\gamma_L)\bar{y}_L\vert&\,|I|+|J|=p,|K|+|L|=q\nonumber\\
	&\text{ such that }(\beta_J\gamma_L)_{|_\Gamma}=1\rangle\label{bicompl}
	\end{align}
	for $(p,q)\in\mathbb{N}^2$. Then the inclusion $(B^{\bullet,\bullet}_\Gamma,\bar{\partial})\hookrightarrow(\wedge^{\bullet,\bullet}(G/\Gamma),\bar{\partial})$ induces an isomorphism in cohomology:
	\begin{equation*}
		H^{\bullet,\bullet}_{\bar{\partial}}(B^{\bullet,\bullet
		}_{\Gamma})\cong H^{\bullet,\bullet}_{\bar{\partial}}(G/\Gamma).
	\end{equation*}
\end{theorem} 

Here, for characters $\{\eta_i\}_{i=1}^k$, we use the following notation: we write $\eta_I:=\eta_{i_1}\cdot...\cdot\eta_{i_l}$ for $I=\{i_1,...,i_l\}$, $i_1<...<i_l$.

The Bott-Chern cohomology of a solvmanifold endowed with a complex structure of splitting type can be computed thanks to the following result proved by Angella and Kasuya.

\begin{theorem}\cite[Corollary 2.15, Theorem 2.16]{AngellaKasuya_2017}\label{thmAngKas}
	Let $G/\Gamma$ be a solvmanifold endowed with an invariant complex structure of splitting type and $C^{\bullet,\bullet}_\Gamma\subset\wedge^{\bullet,\bullet}(G/\Gamma)$ be the finite-dimensional subspace given by
	\begin{equation}
		C^{\bullet,\bullet}_\Gamma:=B^{\bullet,\bullet}_\Gamma\,+\,\bar{B}^{\bullet,\bullet}_\Gamma,
	\end{equation}
	where $B^{\bullet,\bullet}_\Gamma$ is defined by \eqref{bicompl} and $\bar{B}^{\bullet,\bullet}_\Gamma:=\{\bar{\omega}\in\wedge^{\bullet,\bullet}(G/\Gamma)\,|\,\omega\in B^{\bullet,\bullet}_\Gamma\}$. Then the inclusion $(C^{\bullet,\bullet}_\Gamma,\partial,\bar{\partial})\hookrightarrow(\wedge^{\bullet,\bullet}(G/\Gamma),\partial,\bar{\partial})$ induces isomorphisms in cohomology:
	\begin{align*}
		&H^{\bullet,\bullet}_{\bar{\partial}}(C^{\bullet,\bullet}_\Gamma)\cong H^{\bullet,\bullet}_{\bar{\partial}}(G/\Gamma);\\
		&H^{\bullet,\bullet}_{\partial}(C^{\bullet,\bullet}_\Gamma)\cong H^{\bullet,\bullet}_{\partial}(G/\Gamma);\\
		&H^{\bullet,\bullet}_{BC}(C^{\bullet,\bullet}_\Gamma)\cong H^{\bullet,\bullet}_{BC}(G/\Gamma);\\
		&H^{\bullet,\bullet}_{A}(C^{\bullet,\bullet}_\Gamma)\cong H^{\bullet,\bullet}_{A}(G/\Gamma).
	\end{align*}
	
\end{theorem}

As a consequence of the two previous theorems, it has been proved in \cite{otalthesis} the following result.

\begin{lemma}\label{lemmaAngKas}\cite[Lemma 4.2.13]{otalthesis}
	Let $G/\Gamma$ be a complex solvmanifold of splitting type and $ (B^{\bullet, \bullet}_\Gamma,\partial,\bar{\partial}) $ the complex defined by \eqref{bicompl}. If $\partial_{|_{B^{\bullet, \bullet}_\Gamma}}=\bar{\partial}_{|_{B^{\bullet, \bullet}_\Gamma}}=0$ and $B^{q,p}_\Gamma=\overline{B^{p,q}_\Gamma}$ for all $p,q\in \mathbb{N}$, then $(G/\Gamma, J)$ satisfies the $\partial\bar{\partial}$-Lemma.
\end{lemma}

To apply Theorem \ref{thmKasuya} and Theorem \ref{thmAngKas} we also need to choose lattices $\Gamma=\Gamma'\ltimes_\varphi\Gamma''$ in $G$, where $\Gamma''$ is a lattice of $\C^2$ and $\Gamma'$ is a lattice of $\C$ compatible with the splitting, that is $\varphi(z_3)(\Gamma'')\subseteq\Gamma''$ for all $z_3\in\Gamma'$. This implies that the matrix $\varphi(z_3)_{|_{\Gamma'}}$ needs to be conjugated to an integer matrix for all $z_3\in\Gamma'$. This condition can be directly verified applying the following Lemma.

\begin{lemma}\cite[Lemma 4.2.17]{otalthesis}\label{compatibile}
	A matrix $M_f=\emph{diag}(e^f, e^{-f})\in\emph{GL}(2,\C)$ with $f\in\C$ is in the class of conjugation of an integer matrix if and only if $f=\log(\frac{n+\sqrt{n^2-4}}{2})$ with $n\in\mathbb{Z}$.
\end{lemma}

Hence the compatibility of the lattices given by \cite[Proposition 2.10]{fino2015dimensionalsolvmanifoldsholomorphicallytrivial} and by \eqref{lattice} with the splitting can be verified by applying this Lemma.

We now proceed by considering each Lie algebra separately. From now on, to make the discussion more concise, we will denote solvmanifolds carrying a complex structure of splitting type and with holomorphically trivial canonical bundle by splitting Calabi-Yau manifolds (\cite{Tosatti_2015}).

\subsection{Dolbeault and Bott-Chern cohomology for $\g_1$}

This case correspond to the completely solvable Nakamura manifold and has already been studied by Angella and Kasuya in \cite{kasuya2012techniquescomputationsdolbeaultcohomology} and \cite{AngellaKasuya_2017}. The expression of the action $\varphi_A$ given by \eqref{action} for $A=1$ corresponds to the complex structure \eqref{str1}. The lattices of the connected and simply-connected Lie group $G_1$ with underlying Lie algebra $\g_1$ are $\Gamma_{n,b}=(\log(\frac{n+\sqrt{n^2-4}}{2})\mathbb{Z}\times ib\mathbb{Z})\ltimes_{\varphi_A}\Gamma'$ for $0\neq b\in\R$ and $\Gamma'$ a lattice of $\C^2$. Then Lemma \ref{compatibile} grants that for such lattices ${\varphi_A}_{|_{\Gamma_{n,b}}}$ is in the conjugation class of an integer matrix, hence the lattices $\Gamma_{n,b}$ are compatible with the splitting. 

The computation of the double complexes $B^{\bullet,\bullet}_{\Gamma_{n,b}}$ and $C^{\bullet,\bullet}_{\Gamma_{n,b}}$ as well as the computation of cohomology groups $H^{\bullet,\bullet}_{\bar{\partial}}(G_1/\Gamma_{n,b})$ and $H^{\bullet,\bullet}_{BC}(G_1/\Gamma_{n,b})$ are done in \cite[Example 1]{kasuya2012techniquescomputationsdolbeaultcohomology} and \cite[Example 3.1]{AngellaKasuya_2017}. The three cases considered in the references correspond, in our notation, to $b=k\pi$ for case $(i)$, to $b=(2k+1)\frac{\pi}{2}$ for case $(ii)$ and to $b\neq k\frac{\pi}{2}$ for case $(iii)$. The double complex $B^{\bullet,\bullet}_\Gamma$ depends only on the parameter $b$, hence the lattices that we have found in Section \ref{section1} give the same cohomologies obtained in \cite[Example 1]{kasuya2012techniquescomputationsdolbeaultcohomology} and in \cite[Example 3.1]{AngellaKasuya_2017}.

The dimensions of the cohomology groups are summarized in Table \ref{table1} (see \cite{kasuya2012techniquescomputationsdolbeaultcohomology} and \cite{AngellaKasuya_2017} for explicit computations). We remark that solvmanifolds in case $(iii)$ satisfy $\partial\bar{\partial}$-Lemma. Hence we can easily compute the de Rham cohomology: if a compact complex manifold satisfies the $\partial\bar{\partial}$-Lemma, then it holds that
\begin{equation*}
H^k_{dR}(M)\cong \bigoplus_{p+q=k}H^{p,q}_{\bar{\partial}}(M).
\end{equation*} 
Thanks to \cite{Hattori_60}, we have that the de Rham cohomology for a completely solvable solvmanifold $G/\Gamma$ does not depend on the lattice $\Gamma$, hence it is the same as for cases $(i)$ and $(ii)$.

\subsection{Dolbeault and Bott-Chern cohomology for $\g_2^\alpha$, $\alpha\ge0$}\label{dolbeault_g2}

For this case, Otal computed in \cite{otalthesis} the Dolbeault cohomology. We have completed his work by computing the bicomplex $C^{\bullet,\bullet}_\Gamma$ and the Bott-Chern cohomology of the considered manifolds.

The characters $\alpha_j$ in the description of the Lie group $G$ as a semi-direct product $G=\C\ltimes_{\varphi_A}\C^2$ are defined in Proposition \ref{PropCaratteri} for $A=e^{i\theta}, \theta\in(0,\pi)$. Considering $\{z_1,z_2\}$ coordinates on the factor $\C^2$ and $z_3$ a coordinate on the factor $\C$, we have a basis of invariant $(1,0)$-forms
\begin{equation*}
	\omega^1=\alpha_1^{-1}dz_1,\quad\omega^2=\alpha_2^{-1}dz_2,\quad\omega^3=dz_3.
\end{equation*}
The unitary characters $\beta_1,\beta_2,\gamma_1,\gamma_2$ such that $\alpha_1\beta_1^{-1}, \alpha_2\beta_2^{-1}, \bar{\alpha}_1\gamma_1^{-1}, \bar{\alpha}_2\gamma_2^{-1}$ are holomorphic are the following:
\begin{equation*}
	\begin{array}{cc}
	\beta_1(z_3)=e^{-\bar{A}z_3+A\bar{z}_3},\quad&\beta_2(z_3)=\beta_1(z_3)^{-1}=e^{\bar{A}z_3-A\bar{z}_3},\\
	\gamma_1(z_3)=e^{-Az_3+\bar{A}\bar{z}_3},\quad&\gamma_2(z_3)=\gamma_1(z_3)^{-1}=e^{Az_3-\bar{A}\bar{z}_3}.
	\end{array}
\end{equation*}
Depending on the value of these characters on the considered lattice $\Gamma$, the generators of the complex $B^{\bullet,\bullet}_{\Gamma}$ are taken among the following differential forms:
\begin{equation*}
\begin{cases}
\beta_1\omega^1=e^{-(A+\bar{A})z_3}dz_1\\
\beta_2\omega^2=e^{(A+\bar{A})z_3}dz_2\\
\omega^3=dz_3
\end{cases}
\begin{cases}
\gamma_1\omega^{\uno}=e^{-(A+\bar{A})z_3}d\bar{z}_1\\
\gamma_2\omega^{\due}=e^{(A+\bar{A})z_3}d\bar{z}_2\\
\omega^{\tre}=d\bar{z}_3
\end{cases}
\end{equation*} 

It is shown in \cite{otalthesis} that the lattices $\Gamma'\subseteq\C$ such that $\Gamma=\Gamma'\ltimes_{\varphi_A}\Gamma''$ is compatible with the splitting are $\Gamma'=x_3\mathbb{Z}\oplus ib\mathbb{Z}$, where
\begin{itemize}
	\item for $\g_2^0$, $x_3\in\{\frac{\pi}{2},\frac{\pi}{3},\frac{\pi}{4},\frac{\pi}{6}\}$ with $A=e^{i\frac{\pi}{2}}=i$;
	\item for $\g_2^{\alpha}$, $\alpha>0$, $x_3(n)=\frac{\pi}{2\mathfrak{Im} A_n}$ with $A_n=e^{i\theta_n}$ such that $\theta_n\in(0,\pi)$ satisfies
	\begin{equation*}
		\tan\theta_n=\dfrac{\pi}{\log(|\frac{n+\sqrt{n^2-4}}{2}|)}\quad\text{for }n\le-3.
	\end{equation*}
	Then we have $\alpha=\alpha_n=\frac{1}{\pi}|\log(|\frac{n+\sqrt{n^2-4}}{2}|)|$.
\end{itemize}

We start with $\g_2^0$. We have that ${\beta_1}_{|_\Gamma}\neq1$, ${\gamma_1}_{|_\Gamma}\neq1$ and $(\beta_1\gamma_1)_{|_\Gamma}=1$ for any admissible $x_3$; $(\beta_1\gamma_1^{-1})_{|_\Gamma}=1$ if and only if $x_3=\frac{\pi}{2}$. The computation of the double complex $B^{\bullet,\bullet}_\Gamma$, made by Otal in \cite[Table 4.2]{otalthesis}, is summarized in Table \ref{tabellag02}. For any value $x_3$ it shows that the solvmanifolds $G_2^0/\Gamma$ satisfy the $\partial\bar{\partial}$-Lemma thanks to Lemma \ref{lemmaAngKas}. Moreover, in \cite{AngellaOtalUgarte_2017} it is shown that such manifolds are K\"ahler. Hence it holds
\begin{equation*}
H_{BC}^{\bullet,\bullet}(G_2^0/\Gamma)\cong H_{BC}^{\bullet,\bullet}(C^{\bullet,\bullet}_\Gamma)\cong C^{\bullet,\bullet}_\Gamma=B^{\bullet,\bullet}_\Gamma\cong H^{\bullet,\bullet}_{\bar{\partial}}(B^{\bullet,\bullet}_\Gamma)
\cong H^{\bullet,\bullet}_{\bar{\partial}}(G_2^0/\Gamma).
\end{equation*}

We consider now $\g_2^{\alpha_n}$ for $\alpha_n=\frac{1}{\pi}|\log(|\frac{n+\sqrt{n^2-4}}{2}|)|$, $n\le-3$. In this case, the double complex $B^{\bullet,\bullet}_\Gamma$ depends on the real value $0\neq b\in\R$ in the lattice $\Gamma'=\frac{\pi}{2\mathfrak{Im} A_n}\mathbb{Z}\oplus ib\mathbb{Z}$. Following \cite{otalthesis}, we have:
\begin{equation}
\begin{cases}\label{casi}
(\beta_1\gamma_1^{-1})_{|_\Gamma}=1 &\text{for any }b\in\R;\\
(\beta_1\gamma_1)_{|_\Gamma}=1 &\text{if and only if }b=k\frac{\pi}{2\mathfrak{Re} A_n}\text{ with }k\in\mathbb{Z};\\
{\beta_1}_{|_\Gamma}={\gamma_1}_{|_\Gamma}=1 &\text{if and only if }b=(2k+1)\frac{\pi}{2\mathfrak{Re} A_n}\text{ with }k\in\mathbb{Z}.
\end{cases}
\end{equation}
We summarize the computation made in \cite{otalthesis} of $B^{\bullet,\bullet}_\Gamma$ in Table \ref{tabellaBg2}. 

We now compute the double complexes $C^{\bullet,\bullet}_\Gamma$ depending on the value of $b$ varying as in \eqref{casi}. Our results on Bott-Chern cohomology are summarized in Table \ref{tabellaCg2}. For $b\neq k\frac{\pi}{2\mathfrak{Re}A_n}$ the complex $B^{\bullet,\bullet}_\Gamma$ satisfies the hypothesis of Lemma \ref{lemmaAngKas}, then we have the isomorphism $H^{\bullet,\bullet}_{BC}(G_2^{\alpha_n}/\Gamma)\cong B^{\bullet,\bullet}_\Gamma$ and the validity of the $\partial\bar{\partial}$-Lemma. 

We observe that the cohomologies that we have obtained in the three cases distinguished for the algebra $\g_2^{\alpha}$ with $\alpha>0$ are isomorphic to the cohomologies obtained for $\g_1$ by Angella and Kasuya in \cite{kasuya2012techniquescomputationsdolbeaultcohomology} and \cite{AngellaKasuya_2017} (actually, there is an isomorphism between the corresponding bicomplexes $C_\Gamma^{\bullet,\bullet}$). This could suggest the existence of biholomorphisms between the splitting Calabi-Yau manifolds risen from these algebras.

\subsection{Dolbeault and Bott-Chern cohomology for $\g_8$}\label{dolbeault_g8}

Some of the manifolds included in this case have been studied yet in \cite{otalthesis}, \cite{kasuya2012techniquescomputationsdolbeaultcohomology} and \cite{AngellaKasuya_2017}. In particular, in \cite[Example 2]{kasuya2012techniquescomputationsdolbeaultcohomology} and \cite[Example 3.4]{AngellaKasuya_2017} Angella and Kasuya computed the Dolbeault and the Bott-Chern cohomologies of the complex parallelizable Nakamura manifold, in \cite[Table 4.5]{otalthesis} Otal summarise his computations of Dolbeault cohomology. We have concluded the study completing the computations for all the splitting Calabi-Yau manifolds we can construct from the Lie algebra $\g_8$.

The characters $\alpha_j$ in the description of the Lie group $G$ as a semi-direct product $G=\C\ltimes_{\varphi_A}\C^2$, as stated in Proposition \ref{PropCaratteri}, are 
\begin{equation}\label{caratteri8}
\alpha_1(z_3)=e^{-(A-i)z_3-(A+i)\bar{z}_3},\quad \alpha_2(z_3)=e^{(A-i)z_3+(A+i)\bar{z}_3}
\end{equation}
where $A\in\C$, $\mathfrak{Im}A\neq0$. We consider $\{z_1,z_2\}$ coordinates on the factor $\C^2$ and $z_3$ a coordinate on $\C$. We then have a basis of invariant $(1,0)$-forms
\begin{equation*}
\omega^1=\alpha_1^{-1}dz_1,\quad\omega^2=\alpha_2^{-1}dz_2,\quad\omega^3=dz_3.
\end{equation*}
The unitary characters $\beta_1,\beta_2,\gamma_1,\gamma_2$ such that the products $\alpha_1\beta_1^{-1}, \alpha_2\beta_2^{-1}, \bar{\alpha}_1\gamma_1^{-1}, \bar{\alpha}_2\gamma_2^{-1}$ are holomorphic functions are the following:
\begin{equation*}\label{characters8}
\begin{array}{cc}
\beta_1(z_3)=e^{(\bar{A}-i)z_3-(A+i)\bar{z}_3},\quad&\beta_2(z_3)=\beta_1(z_3)^{-1}=e^{-(\bar{A}-i)z_3+(A+i)\bar{z}_3},\\
\gamma_1(z_3)=e^{(A-i)z_3-(\bar{A}+i)\bar{z}_3},\quad&\gamma_2(z_3)=\gamma_1(z_3)^{-1}=e^{-(A-i)z_3+(\bar{A}+i)\bar{z}_3}.
\end{array}
\end{equation*}
Following \cite{otalthesis}, the generators of the complex $B^{\bullet,\bullet}_\Gamma$ are taken, depending on the lattice $\Gamma$, among the following forms:
\begin{equation*}
\begin{cases}
\beta_1\omega^1=e^{(A+\bar{A}-2i)z_3}dz_1\\
\beta_2\omega^2=e^{-(A+\bar{A}-2i)z_3}dz_2\\
\omega^3=dz_3
\end{cases}
\begin{cases}
\gamma_1\omega^{\uno}=e^{(A+\bar{A}-2i)z_3}d\bar{z}_1\\
\gamma_2\omega^{\due}=e^{-(A+\bar{A}-2i)z_3}d\bar{z}_2\\
\omega^{\tre}=d\bar{z}_3
\end{cases}
\end{equation*} 

Similarly to $\g_2^\alpha$, it is shown in \cite{otalthesis} how to construct some lattices $\Gamma$ which are compatible with the splitting for each fixed $A$ with $\mathfrak{Im}A\neq0$. Let $F_A:\Z\to\C$ be the function defined by
\begin{equation*}
	F_A(n):=\dfrac{1}{2\mathfrak{Im}A}\text{Arg}\left(\frac{n+\sqrt{n^2-4}}{2}\right)+\dfrac{i}{2}\left(\log\left|\frac{n+\sqrt{n^2-4}}{2}\right|-\dfrac{\mathfrak{Re}A}{\mathfrak{Im}A}\text{Arg}\left(\frac{n+\sqrt{n^2-4}}{2}\right)\right).
\end{equation*}
Then, if $\Gamma'=F_A(n)\Z\oplus F_A(n')\Z\subseteq\C$ with $n\neq n'$ and $\Gamma''$ is a lattice for the factor $\C^2$, the lattice $\Gamma=\Gamma'\ltimes_{\varphi_A}\Gamma''$ is compatible with the splitting.

In order to apply Theorem \ref{thmKasuya} and Theorem \ref{thmAngKas}, we tabulate the different values of the map $F_A:\Z\to\C$ in the following table.

\begin{table}[H]
\begin{center}
\begin{tabular}{c|c}
	$n$&$F_A(n)$\\
	\hline
	$\le-3$&$ \frac{\pi}{2\mathfrak{Im} A}+\frac{i}{2}(\log|\frac{n+\sqrt{n^2-4}}{2}|-\frac{\mathfrak{Re}A}{\mathfrak{Im}A}\pi) $\\
	$ -2 $ & $\frac{\pi}{2\mathfrak{Im} A}(1-i\mathfrak{Re}A)$\\
	$-1$ & $\frac{\pi}{3\mathfrak{Im} A}(1-i\mathfrak{Re}A)$\\
	$0$ & $\frac{\pi}{4\mathfrak{Im} A}(1-i\mathfrak{Re}A)$\\
	$1$ & $\frac{\pi}{6\mathfrak{Im} A}(1-i\mathfrak{Re}A)$\\
	$2$ & $0$\\
	$ \ge 3 $ & $ \frac{i}{2}\log(\frac{n+\sqrt{n^2-4}}{2}) $
\end{tabular}
\end{center}
\caption{Values of the function $F_A(n)$.}\label{tabellaF}
\end{table}

In \cite{otalthesis}, the author specializes the values $n$ and $n'$ considering lattices of the form $\Gamma'=F_A(-2)\Z\oplus F_A(n)\Z$ for $n\ge3$. With such lattices he obtains 3 classes of manifolds with different cohomology. We have generalized the computation, obtaining the following result.
	
\begin{theorem}\label{teorema}
	Let $G=\C\ltimes_{\varphi_A}\C^2$ be a Lie group endowed with an invariant complex structure of splitting type, where $\varphi_A$ is described by \eqref{action} and \eqref{caratteri8}. Let $\Gamma=\Gamma'\ltimes_{\varphi_A}\Gamma''$ be a lattice of $G$ compatible with the splitting, where $\Gamma'=F_A(n)\Z\oplus F_A(n')\Z$ for $n,n'\in\Z$, $n\neq n'$. Then, varying $A\in\C$ and $n,n'\in\Z$, we obtain one of the following cases:
	\begin{enumerate}[label=$(\roman*)$,itemsep=.8ex]
		\item\label{caso1} ${\beta_1}_{|_\Gamma},{\gamma_1}_{|_\Gamma},(\beta_1\gamma_1^{-1})_{|_\Gamma},(\beta_1\gamma_1)_{|_\Gamma}\neq1$;
		\item\label{caso2} ${\beta_1}_{|_\Gamma}=1$, ${\gamma_1}_{|_\Gamma},(\beta_1\gamma_1^{-1})_{|_\Gamma}, (\beta_1\gamma_1)_{|_\Gamma}\neq1$;
		\item\label{caso3} ${\gamma_1}_{|_\Gamma}=1$, ${\beta_1}_{|_\Gamma},(\beta_1\gamma_1^{-1})_{|_\Gamma}, (\beta_1\gamma_1)_{|_\Gamma}\neq1$;
		\item\label{caso4} $(\beta_1\gamma_1^{-1})_{|_\Gamma}=1$, ${\beta_1}_{|_\Gamma},{\gamma_1}_{|_\Gamma},(\beta_1\gamma_1)_{|_\Gamma}\neq1$;
		\item\label{caso5} $(\beta_1\gamma_1)_{|_\Gamma}=1$, ${\beta_1}_{|_\Gamma},{\gamma_1}_{|_\Gamma},(\beta_1\gamma_1^{-1})_{|_\Gamma}\neq1$;
		\item\label{caso6} $(\beta_1\gamma_1^{-1})_{|_\Gamma}=(\beta_1\gamma_1)_{|_\Gamma}=1$, ${\beta_1}_{|_\Gamma},{\gamma_1}_{|_\Gamma}\neq1$;
		\item\label{caso7} ${\beta_1}_{|_\Gamma}={\gamma_1}_{|_\Gamma}=(\beta_1\gamma_1^{-1})_{|_\Gamma}=(\beta_1\gamma_1)_{|_\Gamma}=1$.		
	\end{enumerate}

Moreover, different cases have different cohomologies, hence we obtain seven non biholomorphic splitting Calabi-Yau manifolds starting from the Lie algebra $\g_8$.
\begin{remark}
	Notice that the seven cases considered in this theorem are all the combinations allowed by Theorem \ref{thmKasuya} and Theorem \ref{thmAngKas}, hence this result covers all the possible cohomology types that we can obtain in this context.
\end{remark}	
\begin{proof}	
The proof consists, for each case stated in the Theorem, in choosing suitable  $A$, $n$ and $n'$ such that the products of the characters defined by \eqref{characters8} assume the wanted values. To do so we rewrite the characters $\beta_1(z_3)$, $\gamma_1(z_3)$, $\beta_1(z_3)\gamma_1(z_3)$, $\beta_1(z_3)\gamma_1^{-1}(z_3)$ in the following way:
\begin{equation}\label{char}
\begin{array}{cc}
\beta_1(z_3)=e^{-2i\mathfrak{Re}(z_3)+2i\mathfrak{Im}(\bar{A}z_3)},&\quad\gamma_1(z_3)=e^{-2i\mathfrak{Re}(z_3)+2i\mathfrak{Im}(Az_3)},\\
\beta_1(z_3)\gamma_1(z_3)=e^{-4i\mathfrak{Re}(z_3)+4i\mathfrak{Re}(A)\mathfrak{Im}(z_3)},&\quad\beta_1(z_3)\gamma_1^{-1}(z_3)=e^{-4i\mathfrak{Im}(A)\mathfrak{Re}(z_3)}.
\end{array}
\end{equation}

Since we have $\beta_2=\beta_1^{-1}$ and $\gamma_2=\gamma_1^{-1}$, we only need to check the values assumed by the listed characters. It is possible to achieve all the cases except case \ref{caso5} by taking $\mathfrak{Re}A=0$. As we can see from Table \ref{tabellaF}, this assumption simplifies the value of the function $F_A$ and of the characters in \eqref{char} computed on the lattice $\Gamma$.

Cases \ref{caso4}, \ref{caso6} and \ref{caso7} has been obtained by Otal in \cite[Proposition 4.2.25]{otalthesis}: taking $n=-2$ and $n'\ge 3$, we have
\begin{equation*}
\begin{cases}
(\beta_1\gamma_1^{-1})_{|_\Gamma}=1 &\text{for any }A\in\C;\\
(\beta_1\gamma_1)_{|_\Gamma}=1 &\text{if and only if }A=\frac{i}{k}\text{ with }0\neq k\in\mathbb{Z};\\
{\beta_1}_{|_\Gamma}={\gamma_1}_{|_\Gamma}=1 &\text{if and only if }A=\frac{i}{2k+1}\text{ with }k\in\mathbb{Z}.
\end{cases}
\end{equation*}

Case \ref{caso2} corresponds to Example 3.4, case $b)$ in \cite{AngellaKasuya_2017}; in our notation this example corresponds to $A=-i$ and, without loss of generality, $n'\in\{-1,0,1\}$ and $n^2\ge 9$. 

We now study the remaining cases. For case \ref{caso1} we observe that, for $A=it$ with $t\in\R, t\neq\pm1$, we have 
\begin{align*}
&{\beta_1}_{|_\Gamma}=1\iff t=\frac{1}{kl-1}\quad\text{for }k\in\Z,\\
&{\gamma_1}_{|_\Gamma}=1\iff t=\frac{1}{kl+1}\quad\text{for }k\in\Z.
\end{align*}
Here, $l$ is the integer such that $\mathfrak{Re}(F_A(n))=\frac{\pi}{l\mathfrak{Im}A}$ as we can see from Table \ref{tabellaF}. Moreover, the product $\beta_1\gamma_1^{-1}$ is $1$ on the lattice $\Gamma$ if and only if $\mathfrak{Im}(A)\mathfrak{Re}(z_3)=k\frac{\pi}{2}$ for a $k\in\Z$. This condition in verified if $n,n'\notin\{-1,0,1\}$ as we can see from Table \ref{tabellaF}. So case \ref{caso1} can be achieved, for example, choosing $n'\in\{-1,0,1\},$ $ n\ge3$ and $A=it$, $ t\in\R$, $t\neq\frac{1}{kl\pm1}$ for $k\in\Z$, where $l\in\Z$ is such that $\mathfrak{Re}(F_A(n'))=\frac{\pi}{l\mathfrak{Im}A}$ for $ n'\in\{-1,0,1\}$.

Case \ref{caso3} is analogous to case \ref{caso2}: it is given, for example, by $A=i$, $n'\in\{-1,0,1\}$ and $n^2\ge 9$. Indeed, for $A=i$, we have $\mathfrak{Im}(Az_3)=\mathfrak{Re}(z_3)$ and $\mathfrak{Im}(\bar{A}z_3)=-\mathfrak{Re}(z_3)$, that implies
\begin{equation*}
{\beta_1}_{|_\Gamma}=e^{-4i\mathfrak{Re}(z_3)}\text{ for }z_3\in\Gamma',\quad {\gamma_1}_{|_\Gamma}\equiv1.
\end{equation*}
So by looking again at Table \ref{tabellaF}, we need to choose (without loss of generality) $n'\in\{-1,0,1\}$ to have ${\beta_1}_{|_\Gamma}\neq1$. Then we need $n^2\ge 9$ to have that $\Gamma'$ is a lattice.

To achieve case \ref{caso5} we need to take $\mathfrak{Re}A\neq0$; since $\beta_1\gamma_1^{-1}$ does not depend on $\mathfrak{Re}(A)$, we can proceed as in case $(i)$ to choose $n$ and $n'$ such that we have ${\beta_1\gamma_1^{-1}}_{|_\Gamma}\neq1$. We then fix $n'\in\{-1,0,1\}$ and $n'\neq n\in\{-2,-1,0,1\}$. Choosing $\mathfrak{Re}(A)=1$, we have ${\beta_1\gamma_1}_{|_\Gamma}\equiv1$. The other two characters satisfy
\begin{align*}
&{\beta_1}_{|_\Gamma}=1\iff \forall z_3\in\Gamma',\quad\mathfrak{Re}(z_3)-\mathfrak{Im}(\bar{A}z_3)=k\pi\,\, \text{ for some }k\in\Z;\\
&{\gamma_1}_{|_\Gamma}=1\iff \forall z_3\in\Gamma',\quad\mathfrak{Re}(z_3)-\mathfrak{Im}(Az_3)=k\pi\,\, \text{ for some }k\in\Z.
\end{align*}
It is easy to check that these conditions are, in general, not verified. Then, for a suitable choice of $\mathfrak{Im}(A)$, we get case \ref{caso5}.
\end{proof}
\end{theorem}

We can now proceed with the computation of the complexes $B^{\bullet,\bullet}_\Gamma$ and $C^{\bullet,\bullet}_\Gamma$ and the corresponding cohomologies. Both Dolbeault and Bott-Chern cohomologies of case \ref{caso2} and case \ref{caso7} have been computed by Angella and Kasuya in \cite{AngellaKasuya_2017}. Dolbeault cohomology of cases \ref{caso4} and \ref{caso6} is presented in \cite{otalthesis}. We then have computed the remaining cases, both for Dolbeault and Bott-Chern. 

Our results are summarized in Appendix \ref{tables}: \crefrange{tabellaBg8123}{tabellaBCg8Otal} show the complexes $B^{\bullet,\bullet}_\Gamma$ and $C^{\bullet,\bullet}_\Gamma$ and the corresponding cohomology for each case listed in Theorem \ref{teorema}. Computations made by Otal, Angella and Kasuya together with ours show that only cases \ref{caso1} and \ref{caso4} satisfy the $\partial\bar{\partial}$-Lemma. Cases \ref{caso4}, \ref{caso6} and \ref{caso7} have the same cohomologies as cases $(iii), (ii)$ and $(i)$ respectively of both $\g_1$ and $\g_2^{\alpha}$, $\alpha>0$ (as observed for $\g_2^{\alpha}$, there is an isomorphism between the corresponding bicomplexes $C_\Gamma^{\bullet,\bullet}$). Again, this fact could suggest the presence of biholomorphisms between the solvmanifolds obtained in these cases. We also have that cases \ref{caso4} and \ref{caso5} of Theorem \ref{teorema} lead to the same Dolbeault cohomology, but only the former case satisfies the $\partial\bar{\partial}$-Lemma.

\section{Decomposition of the double complex $(C_\Gamma^{\bullet,\bullet},\partial, \bar{\partial})$}\label{section3}

In this section, we recall the decomposition of bicomplexes into indecomposable ones, then we proceed with studying the decompositions of the bicomplexes $(C^{\bullet,\bullet}_\Gamma,\partial,\bar{\partial})$ corresponding to the splitting Calabi-Yau manifolds considered in the previous sections (more details on this subject can be found, for example, in \cite{Stelzig_2021} and \cite{khovanov2020}).

A double complex $A$ is called indecomposable if there is no non-trivial decomposition $A=A_1\oplus A_2$. The following double complexes are indecomposable (the drawn components are one dimensional and the drawn maps are isomorphisms, while all components and arrows which are non-drawn are zero):
\begin{itemize}
	\item squares
	\begin{equation*}
	\begin{tikzcd}
		A^{p-1,q}\arrow{r}{\partial}& A^{p,q}\\
		A^{p-1,q-1}\arrow{u}{\bar{\partial}}\arrow{r}{\partial}&A^{p,q-1}\arrow{u}
		{\bar{\partial}}
	\end{tikzcd}
	\end{equation*}
	\item zigzags
	\begin{equation*}
	\begin{tikzcd} 
	A^{p,q} 
	\end{tikzcd}\quad,\quad
	\begin{tikzcd}
	A^{p,q+1} \\
	A^{p,q}\arrow{u}{\bar{\partial}}
	\end{tikzcd}\quad,\quad
	\begin{tikzcd}
	A^{p,q}\arrow{r}{\partial} &A^{p+1,q}
	\end{tikzcd}\quad,\quad
	\begin{tikzcd}
	A^{p,q+1}\\
	A^{p,q}\arrow{u}{\bar{\partial}} \arrow{r}{\partial} &A^{p+1,q}
	\end{tikzcd}\quad,
	\end{equation*}
	\begin{equation*}
	\begin{tikzcd}
	A^{p-1,q}\arrow{r}{\partial}& A^{p,q}\\
	&A^{p,q-1}\arrow{u}{\bar{\partial}}
	\end{tikzcd}\quad,\quad
	\begin{tikzcd}
	A^{p,q}\arrow{r}{\partial}& A^{p+1,q}&\\
	&A^{p+1,q-1}\arrow{u}{\bar{\partial}}\arrow{r}{\partial}&A^{p+2,q-1}
	\end{tikzcd}\quad,\quad\dots
	\end{equation*}
\end{itemize}
The shape of a square or a zigzag $A$ is defined as 
\begin{equation*}
	S(A):=\{(p,q)\in\Z^2\,|\,A^{p,q}\neq 0\}.
\end{equation*}
The isomorphism class of a square or a zigzag is determined uniquely by its shape. We can then associate to each shape $S$ a square or a zigzag of that shape $C(S)$ in which all non-zero components are the field $K$ and all non-zero arrows are $\pm\text{Id}$.

The aim of this section is to apply the following Theorem to the double complex $(C_\Gamma^{\bullet,\bullet},\partial,\bar{\partial})$ associated to a six-dimensional splitting Calabi-Yau manifold.

\begin{theorem}\cite[Theorem 3]{Stelzig_2021}\label{decomposition}
	For any bounded double complex $A$ over a field $K$ (i.e.
	$A^{p,q}=0$ for almost all $(p,q)\in\Z^2$), there exist unique cardinal numbers $\emph{mult}_S(A)$ and an isomorphism $A\simeq\bigoplus_S C(S)^{\oplus\emph{mult}_S(A)}$, where $S$ runs over the shapes of squares and zigzags.
\end{theorem}

This decomposition can be used to obtain interesting results on the cohomology of a double complex $(A,\partial,\bar{\partial})$. Counting corners, startings and endings of zigzags in a certain bidegree gives the dimensions of Dolbeault, Bott-Chern and Aeppli cohomologies. Formally, the following result holds.

\begin{lemma}\cite[Lemma 8]{Stelzig_2021}
	Let $(A,\partial,\bar{\partial})$ be a bounded double complex with a decomposition in elementary complexes with pairwise distinct support $\varphi:\bigoplus A_i\to A$. 
	\begin{itemize}
		\item For every $(p,q)\in \Z^2$, the maps induced by $\varphi$
		
		\begin{equation*}
			\bigoplus_{\substack{A_i \text{ zigzag}\\ (p,q)\in S(A_i)\\ (p-1,q),(p+1,q)\notin S(A_i)}}  H^{p,q}_{\partial}(A_i)\rightarrow H^{p,q}_{\partial}(A)
		\end{equation*}
		
		\begin{equation*}
			\bigoplus_{\substack{A_i \text{ zigzag}\\ (p,q)\in S(A_i)\\ (p,q-1),(p,q+1)\notin S(A_i)}} H^{p,q}_{\bar{\partial}}(A_i)\rightarrow H^{p,q}_{\bar{\partial}}(A)
		\end{equation*}
		are isomorphisms.
		
		\item For every $(p,q)\in \Z^2$, the maps induced by $\varphi$
		
		\begin{equation*}
		\bigoplus_{\substack{A_i \text{ zigzag}\\ (p,q)\in S(A_i)\\ (p+1,q),(p,q+1)\notin S(A_i)}}  H^{p,q}_{BC}(A_i)\rightarrow H^{p,q}_{BC}(A)
		\end{equation*}
		
		\begin{equation*}
		\bigoplus_{\substack{A_i \text{ zigzag}\\ (p,q)\in S(A_i)\\ (p-1,q),(p,q-1)\notin S(A_i)}} H^{p,q}_{A}(A_i)\rightarrow H^{p,q}_{A}(A)
		\end{equation*}
		are isomorphisms.
	\end{itemize}
\end{lemma}

We now proceed with the decomposition in irreducible complexes of the double complexes $(C_\Gamma^{\bullet,\bullet}, \partial,\bar{\partial})$ for the manifolds built from the Lie algebras $\g_1$, $\g_2^{\alpha}$ and $\g_8$ considered in the previous sections. We give in Appendix \ref{decompositions} an explicit and a graphic representation of these decompositions. We can achieve such a decomposition by considering the explicit expression of the bicomplexes $(C^{\bullet,\bullet},\partial,\bar{\partial})$ wrote in Appendix \ref{tables} (see \cite{AngellaKasuya_2017} for $\g_1$) and by computing the images of the differentials $\partial$ and $\bar{\partial}$ for each generator.

As we can see by the structure equations \eqref{str1}, \eqref{str20}, \eqref{str2a} and \eqref{struttura8}, for six-dimensional splitting Calabi-Yau manifolds the bicomplex $B^{\bullet,\bullet}_{\Gamma}$ is generated by holomorphic forms (and hence $\bar{\partial}\equiv0$ on $B^{\bullet,\bullet}_{\Gamma}$) and $\bar{B}^{\bullet,\bullet}_{\Gamma}$ is generated by antiholomorphic forms (${\partial}\equiv0$ on $\bar{B}^{\bullet,\bullet}_{\Gamma}$). This is in fact true in general dimension: if we consider a splitting type solvmanifold $X=G/\Gamma$ with defining Lie group $G=\C^n\ltimes_\varphi\C^m$, we can see by definition \eqref{bicompl} that $\bar{\partial}\equiv0$ on $B^{\bullet,\bullet}_{\Gamma}$ and ${\partial}\equiv0$ on $\bar{B}^{\bullet,\bullet}_{\Gamma}$.

Hence, in the decomposition of the bicomplex $C^{\bullet,\bullet}_{\Gamma}=B^{\bullet,\bullet}_{\Gamma}+\bar{B}^{\bullet,\bullet}_{\Gamma}$ in indecomposable ones we can only have length $2$ zigzags and dots. This condition implies for the bicomplex $C^{\bullet,\bullet}_\Gamma$ to have the page-$1$-$\partial\bar{\partial}$-property, namely its Fr\"olicher spectral sequence degenerates at page $ 2 $ and the Hodge filtration induces a pure Hodge structure on the de Rham cohomology (see \cite[Definition 1]{kasuya2020frolicherspectralsequencehodge} and \cite[Theorem and Definition 1.2]{popovici2022higherpagehodgetheorycompact}). Indeed, the page-$1$-$\partial\bar{\partial}$-property is equivalent to the fact that the bicomplex decomposes in a direct sum of squares, dots and length $2$ zigzags.

This was actually proved by Kasuya and Stelzig to be true in \cite[Theorem 9]{kasuya2020frolicherspectralsequencehodge} for solvmanifolds endowed with a bi-invariant complex structure of splitting type, but if we consider splitting type solvmanifolds with defining Lie group $G=\C^n\ltimes_\varphi\C^m$ we can refine the decomposition: if $X=G/\Gamma$ is a splitting type solvmanifold with $G=\C^n\ltimes_\varphi\C^m$, the bicomplex $C^{\bullet,\bullet}_\Gamma$ can't have squares in its decomposition, then it is a direct sum of dots and length $2$ zigzags. Thanks to \cite[Theorem and Definition 1.2]{popovici2022higherpagehodgetheorycompact}, \cite[(5.21)]{DeligneGriffiths1975} we can also conclude that a bicomplex $C^{\bullet,\bullet}_\Gamma$ defined by a splitting type solvmanifold $X=\C^n\ltimes_\varphi\C^m/\Gamma$ whose decompositions has no lines (e.g. it is direct sums of dots) satisfies the $\partial\bar{\partial}$-Lemma.

\section{On formality of splitting type solvmanifolds}\label{section4}

The aim of this section is to study some notions of formality for six-dimensional splitting Calabi-Yau manifolds and, more in general, for splitting type solvamnifolds $\C^n\ltimes_\varphi N/\Gamma$ with $N=\C^m$. We recall that a cbba is a bicomplex $A$ carrying the additional structure of a graded-commutative product, where the graded commutativity considers the total degree. A weak equivalence is a map $f:A\to B$ between two cbbas whose underlying map of double complexes is a bigraded quasi-isomorphism, i.e. satisfies one of the following equivalent conditions (for more details see \cite[Theorem 1.21, Definition 1.22]{stelzig2023pluripotentialhomotopytheory}):
\begin{enumerate}
	\item it is a bigraded chain homotopy equivalence, i.e. it is an isomorphism in the homotopy category $\text{Ho(BiCo)}$ of bicomplexes;
	\item $\text{cone}(f):=\text{coker}(A\to B\oplus(\square[-1,-1]\otimes A))$, where $\square[-1,-1]$ is the square with bottom-left corner in bidegree $(-1,-1)$, is a direct sum of squares;
	\item for any decomposition $A=A^{sq}\oplus A^{zig}$ and $A=B^{sq}\oplus B^{zig}$ in direct sum of squares and zigzags, the map $f:A^{zig}\to B^{zig}$ is an isomorphism;
	\item the induced maps $H_{BC}(f)$ and $H_A(f)$ are isomorphisms.
\end{enumerate}
If the cbbas $A$ and $B$ are locally bounded (i.e. they do not contain infinite-length zigzags in their decompositions), these conditions are also equivalent to $H_{\partial}(f)$ and $H_{\bar{\partial}}(f)$ being isomorphisms. Then, according to \cite[Theorem C]{stelzig2023pluripotentialhomotopytheory}, for any cohomological functor $H$ (i.e. an additive functor $H$ from the category of bicomplexes to an additive category such that $H(P)=0$ for any direct sum of squares $ P $, see \cite[Definition 8.2]{STELZIG2022108560}) the map $H(f)$ is an isomorphism. If we then have that two cbbas $A$ and $B$ are connected by a chain of bigraded quasi-isomorphisms, then they have isomorphic cohomologies. By a chain of bigraded quasi-isomorphism we mean a diagram
\begin{equation*}
	A\overset{f_1}{\longrightarrow}C_1\overset{f_2}{\longleftarrow}C_2\overset{f_3}{\longrightarrow}...\overset{f_{n-1}}{\longleftarrow}C_{n-1}\overset{f_n}{\longrightarrow}B,
\end{equation*}
where $C_i$ are cbbas and $f_i$ are bigraded quasi-isomorphisms.

We start with the notion of geometrically-Bott-Chern formality. Let $X$ be a compact complex manifold. Let $g$ be a Hermitian metric on $X$ and $*^g$ be the complex Hodge operator associated with $g$. Let $\tilde{\Delta}_{BC}^g$ be the $4$-th order elliptic self-adjoint differential operator (\cite[Proposition 5]{5a8ecba0-bf3c-3a04-94d0-f3fe0a96f498}, \cite[\S 2.b]{schweitzer2007autourlacohomologiebottchern}) defined as
\begin{equation*}
\tilde{\Delta}_{BC}^g:=(\partial\bar{\partial})(\partial\bar{\partial})^*+(\partial\bar{\partial})^*(\partial\bar{\partial})+\left(\bar{\partial}^*\partial\right)\left(\bar{\partial}^*\partial\right)^*+\left(\bar{\partial}^*\partial\right)^*\left(\bar{\partial}^*\partial\right)+\bar{\partial}^*\bar{\partial}+\partial^*\partial.
\end{equation*}
According to \cite{schweitzer2007autourlacohomologiebottchern}, the isomorphism $H_{BC}^{\bullet,\bullet}(X)\simeq(\ker\tilde{\Delta}^g_{BC})^{\bullet,\bullet}$ holds. In general, the bicomplex $(\ker\tilde{\Delta}^g_{BC})^{\bullet,\bullet}$ does not carry an algebra structure, hence it is not a cbba.

Motivated by \cite{Kotschick_2001}, Angella and Tomassini introduced in \cite{Angella_2015} the following notion of formality related to Bott-Chern cohomology.

\begin{definition}\cite[Definitions 1.1, 1.2]{Angella_2015} 
	A Hermitian metric $g$ on $X$ is said to be geometrically-Bott-Chern-formal if the inclusion $\iota_{BC}^g:(\ker\tilde{\Delta}_{BC}^g,0,0)\hookrightarrow(\wedge^{\bullet,\bullet}X,\partial,\bar{\partial})$ is a morphism of cbbas, and it is such that $H_{BC}(\iota_{BC}^g)$ is an isomorphism (equivalently if $\ker\tilde{\Delta}_{BC}^g$ is an algebra according to \cite[Definition 1.1]{tardini2015geometricbottchernformalitydeformations}). A compact complex manifold $X$ is said to be geometrically-Bott-Chern-formal if there exists a geometrically-Bott-Chern-formal metric on $X$.
\end{definition}

As an obstruction to this property there are triple Aeppli-Bott-Chern-Massey products, defined by Angella and Tomassini in \cite{Angella_2015} as follows.

\begin{definition}\cite[Definitions 2.1, 2.2]{Angella_2015}
	Let $(A^{\bullet,\bullet},\partial,\bar{\partial})$ be a cbba. Take cohomology classes 
	\begin{equation*}
	\mathfrak{a}_{12}=[\alpha_{12}]\in H_{BC}^{p,q}(A^{\bullet,\bullet}),\quad \mathfrak{a}_{23}=[\alpha_{23}]\in H_{BC}^{r,s}(A^{\bullet,\bullet}), \quad \mathfrak{a}_{34}=[\alpha_{34}]\in H_{BC}^{u,v}(A^{\bullet,\bullet})
	\end{equation*}
	such that $\mathfrak{a}_{12}\cup\mathfrak{a}_{23}=0$ in $H^{p+r,q+s}_{BC}(A^{\bullet,\bullet})$ and $\mathfrak{a}_{23}\cup\mathfrak{a}_{34}=0$ in $H^{r+u,s+v}_{BC}(A^{\bullet,\bullet})$ and let $\alpha_{13}$ and $\alpha_{24}$ such that
	\begin{equation*}
	(-1)^{p+q}\alpha_{12}\wedge\alpha_{23}=\partial\bar{\partial}\alpha_{13}, \quad (-1)^{r+s}\alpha_{23}\wedge\alpha_{34}=\partial\bar{\partial}\alpha_{24}.
	\end{equation*}
	The triple Aeppli-Bott-Chern-Massey product $\langle\mathfrak{a}_{12},\mathfrak{a}_{23},\mathfrak{a}_{34}\rangle_{ABC}$ is defined as
	\begin{align}
		\langle\mathfrak{a}_{12},\mathfrak{a}_{23}&,\mathfrak{a}_{34}\rangle_{ABC}:=[(-1)^{p+q}\alpha_{12}\wedge\alpha_{24}-(-1)^{r+s}\alpha_{13}\wedge\alpha_{34}]\\\nonumber
		\in&\dfrac{H_A^{p+r+u-1,q+s+v-1}(A^{\bullet,\bullet})}{\mathfrak{a}_{12
		}\cup H_A^{r+u-1,s+v-1}(A^{\bullet,\bullet})+H_A^{p+r-1,q+s-1}(A^{\bullet,\bullet})\cup\mathfrak{a}_{34}}.
	\end{align}
	
	If $X$ is a complex manifold, the triple Aeppli-Bott-Chern-Massey products are defined as the triple Aeppli-Bott-Chern-Massey products of the cbba $(\wedge^{\bullet,\bullet}X,\partial,\bar{\partial})$.
\end{definition}
The relation between geometrically-Bott-Chern-formality and ABC-Massey products is given by the following theorem.

\begin{theorem}\cite[Theorem 2.4]{Angella_2015}
	 On compact complex geometrically-Bott-Chern-formal manifolds, triple Aeppli-Bott-Chern-Massey products vanish. 
\end{theorem}

Related to complex manifolds, we are also interested in the following notions of formality concerning the bigrading of the algebra of differential forms.

\begin{definition}
	A cbba $(A^{\bullet,\bullet},\partial,\bar{\partial})$ is called:
	\begin{enumerate}
		\item Dolbeault formal (\cite{7c0de90b-83a5-334f-b43c-69140714ab91}) if it is connected by a chain of morphisms inducing isomorphisms in $H_{\bar{\partial}}$-cohomology to a cbba $(H^{\bullet,\bullet}, \partial_H, \bar{\partial}_H\equiv0)$;
		
		\item weakly formal (\cite{milivojevic2024bigradednotionsformalityaepplibottchernmassey}) if it is connected by a chain of weak equivalences to a cbba $(H^{\bullet,\bullet},\partial_H, \bar{\partial}_H)$ such that $\partial_H\bar{\partial}_H=0$;
		
		\item strongly formal (\cite{milivojevic2024bigradednotionsformalityaepplibottchernmassey}) if it is connected by a chain of weak equivalences to a cbba $(H^{\bullet,\bullet}, \partial_H\equiv0, \bar{\partial}_H\equiv0)$.
	\end{enumerate} 
A complex manifold $X$ is said to be formal in the sense of $(1), (2)$ or $(3)$ if the cbba of forms $(\bigwedge^{\bullet,\bullet}X,\partial,\bar{\partial})$ is formal in the sense of $(1), (2)$ or $(3)$ respectively.
\end{definition}

\begin{remark}
	We have that a splitting type solvmanifold $X=(\C^n\ltimes_\varphi N)/\Gamma$ with $N=\C^m$ is always Dolbeault formal: by definition \eqref{bicompl} we have $\bar{\partial}\equiv0$ on $B^{\bullet,\bullet}_{\Gamma}$ and the inclusion $B^{\bullet,\bullet}_{\Gamma}\subseteq \wedge^{\bullet, \bullet}X$ induces the wanted isomorphism in cohomology. 
\end{remark}

It is observed in \cite{milivojevic2024bigradednotionsformalityaepplibottchernmassey} that notion $(3)$ of formality implies the $\partial\bar{\partial}$-Lemma property for any compact complex manifold and, if for a complex manifold the $\partial\bar{\partial}$-Lemma holds, notions $(2)$ and $(3)$ are equivalent. If the manifold $X$ is one of the splitting type solvmanifolds we have considered in the previous sections, we can say more. Indeed the following property holds.

\begin{lemma}\label{lemmamio}
	Let $X=(\C^n\ltimes_\varphi N)/\Gamma$ be a splitting type solvmanifold with $N=\C^m$. Then the bicomplex
	$C^{\bullet,\bullet}_{\Gamma}$ satisfies the $\partial\bar{\partial}$-Lemma (and hence $X$ do so) if and only if $B^{\bullet,\bullet}_{\Gamma}=\bar{B}^{\bullet,\bullet}_{\Gamma}$ (and $\partial_{|_{B^{\bullet, \bullet}_\Gamma}}=0$, $\bar{\partial}_{|_{B^{\bullet, \bullet}_\Gamma}}=0$).
	\begin{proof} 
		The $(\impliedby)$ implication is \cite[Lemma 4.2.13]{otalthesis}. 
		
		Since $N=\C^m$ and by the definition of the bicomplexes $B^{\bullet,\bullet}_{\Gamma}$ and $\bar{B}^{\bullet,\bullet}_{\Gamma}$, we have $\bar{\partial}_{|_{B^{\bullet,\bullet}_{\Gamma}}}=0$, $\partial_{|_{\bar{B}^{\bullet, \bullet}_\Gamma}}=0$ (and so $H^{\bullet,\bullet}_{\bar{\partial}}(X)\simeq B^{\bullet,\bullet}_{\Gamma}$). By \cite[Proposition (5.17)]{DeligneGriffiths1975}, a bicomplex satisfying the $\partial\bar{\partial}$-Lemma is a direct sum of dots and squares. Since $\partial\bar{\partial}=0$ on $C^{\bullet,\bullet}_{\Gamma}=B^{\bullet,\bullet}_{\Gamma}+\bar{B}^{\bullet,\bullet}_{\Gamma}$ as observed, we conclude that $C^{\bullet,\bullet}_{\Gamma}$ is a direct sum of dots. So differentials are trivial on $C^{\bullet,\bullet}_{\Gamma}$ and, thanks to \cite[Corollary 2.15]{AngellaKasuya_2017} and to the $\partial\bar{\partial}$-Lemma, we get
		\begin{equation*}
		C^{\bullet,\bullet}_{\Gamma}\simeq H_{BC}^{\bullet,\bullet}(X)\simeq H_{\bar{\partial}}^{\bullet,\bullet}(X)\simeq B^{\bullet,\bullet}_{\Gamma}.
		\end{equation*}
		Hence, since $B^{\bullet,\bullet}_{\Gamma}\subseteq C^{\bullet,\bullet}_{\Gamma}$ and such bicomplexes are of finite dimension, we have the thesis.
	\end{proof} 
\end{lemma}

According to \cite{milivojevic2024bigradednotionsformalityaepplibottchernmassey}, strong formality implies both weakly formality and the $\partial\bar{\partial}$-Lemma. As a consequence of Lemma \ref{lemmamio}, since the bicomplex $B^{\bullet,\bullet}_\Gamma$ is an algebra, we have the following statement.

\begin{corollary}\label{corollarymio}
	Let $X=\C^n\ltimes_\varphi N/\Gamma$ be a splitting type solvmanifold with $N=\C^m$. Then $X$ satisfies the $\partial\bar{\partial}$-Lemma if and only if $X$ is strongly formal.
\end{corollary}

As a consequence of Lemma \ref{lemmamio}, we have that if a splitting type solvmanifold of the form $X=\C^n\ltimes_\varphi\C^m/\Gamma$ satisfies the $\partial\bar{\partial}$-Lemma, then it is geometrically-Bott-Chern formal since in this case we have cbbas isomorphisms $C_\Gamma^{\bullet,\bullet}\simeq H_{BC}^{\bullet,\bullet}(X)\simeq(\ker\tilde{\Delta}^g_{BC})^{\bullet,\bullet}$, i.e. $(\ker\tilde{\Delta}^g_{BC})^{\bullet,\bullet}$ is an algebra.

The next step is the computation of triple ABC-Massey products for six-dimensional splitting Calabi-Yau manifolds in order to have more informations regarding the notions of formality we have just introduced. To do so, we firstly have to find the smallest cbba in $\wedge^{\bullet,\bullet}X$ containing the bicomplex $C^{\bullet,\bullet}_\Gamma$. We can built it by simply adding to the bicomplex the products of its elements. Indeed, the algebra defined in this way is closed under the differentials $\partial$ and $\bar{\partial}$: $\delta(\omega\wedge\eta)=\delta\omega\wedge\eta+(-1)^{|\omega|}\omega\wedge\delta\eta$ is still a linear combination of products of elements in $C^{\bullet,\bullet}_\Gamma$ for $\delta\in\{\partial,\bar{\partial}\}$. Then, since $1\in B^{\bullet,\bullet}_\Gamma\cap \bar{B}^{\bullet,\bullet}_\Gamma$ and $B^{\bullet,\bullet}_\Gamma$ and $\bar{B}^{\bullet,\bullet}_\Gamma$ are algebras, it turns out that such cbba is $B^{\bullet,\bullet}_\Gamma\wedge \bar{B}^{\bullet,\bullet}_\Gamma$. 

We can then compute the triple ABC-Massey products of the solvmanifolds built from the Lie algebras $\g_1, \g_2^{\alpha}$ and $\g_8$ using such cbbas: by looking at the explicit expressions of the bicomplexes $C^{\bullet,\bullet}_\Gamma$ in Tables \ref{tabellaCg2}, \ref{tabellaCg823}, \ref{tabellag85}, \ref{tabellaCg8Otal} and in \cite[Tables 1, 2, 7]{AngellaKasuya_2017}, we can directly show that for a $6$-dimensional splitting Calabi-Yau manifold the cbba $B^{\bullet,\bullet}_\Gamma\wedge\bar{B}^{\bullet,\bullet}_\Gamma$ is obtained just by adding squares to the bicomplex $C^{\bullet,\bullet}_\Gamma$. Squares do not give contributions in cohomology, hence the inclusion $B^{\bullet,\bullet}_\Gamma\wedge\bar{B}^{\bullet,\bullet}_\Gamma\hookrightarrow\wedge^{\bullet,\bullet}X$ is a bigraded-quasi isomorphism and it leads to the identification of the corresponding ABC-Massey products. These computations also show that $B^{\bullet,\bullet}_\Gamma\wedge\bar{B}^{\bullet,\bullet}_\Gamma$ has squares in its decomposition if and only if $C^{\bullet,\bullet}_\Gamma\neq B^{\bullet,\bullet}_\Gamma$. Since according to \cite{milivojevic2024bigradednotionsformalityaepplibottchernmassey} a cbba is weakly formal if and only if it is connected by a chain of weak equivalences to one with no squares in its decomposition, we can conclude that for $6$-dimensional splitting Calabi-Yau manifolds, weak and strong formality are equivalent conditions.

\begin{remark}\label{remark}
	Clearly strong formality implies the vanishing of triple ABC-Massey products. Then thanks to Corollary \ref{corollarymio} we need to compute them only for those splitting Calabi-Yau manifolds not satisfying the $\partial\bar{\partial}$-Lemma.
\end{remark}

\subsection{Triple ABC-Massey products for $\g_1$}

It is shown in \cite[Example 2.3]{tardini2015geometricbottchernformalitydeformations} that solvmanifolds in case $(i)$ admit a non-zero triple ABC-Massey product. Since the complex obtained in this case is isomorphic as a bicomplex to the complex $C^{\bullet,\bullet}_\Gamma$ obtained for $b=(2k+1)\frac{\pi}{2\mathfrak{Re}A_n}$ in $\g_2^{\alpha}$ and in case \ref{caso7} of Theorem \ref{teorema}, we can conclude that these cases admit a non-trivial triple ABC-Massey product. Hence, such manifolds are not geometrically-Bott-Chern-formal. 

We give an explicit computation with our notation: the generators of the bicomplex $B^{\bullet,\bullet}_\Gamma$ for the Lie algebra $\g_1$ as constructed in paragraphs \ref{dolbeault_g2} and \ref{dolbeault_g8} for the other algebras are chosen among $\{e^{-2z_3}dz_1, e^{2z_3}dz_2, dz_3\}$. If we call $T=e^{2z_3}$, we have that the double-complex $C_\Gamma^{\bullet,\bullet}$ is as in Table \ref{tabellaCg2} and in Table \ref{tabellaCg8Otal}. So we can use the basis $\{T^{-1}dz_1, Tdz_2,dz_3\}$.

Consider Bott-Chern classes $\mathfrak{a}_{12}=[Tdz_{3\due}], \mathfrak{a}_{23}=[\barT^{-1}dz_{\uno\tre}], \mathfrak{a}_{34}=[\barT dz_{2\tre}]$. We have 
\begin{equation*}
	Tdz_{3\due}\wedge\barT^{-1}dz_{\uno\tre}=\partial\bar{\partial}(k T\barT^{-1}dz_{\uno\due})\,\, \text{ for a }k\in\C\setminus\{0\}, \quad \barT^{-1}dz_{\uno\tre}\wedge\barT dz_{2\tre}=0.
\end{equation*}
Hence $\langle\mathfrak{a}_{12},\mathfrak{a}_{23},\mathfrak{a}_{34}\rangle_{ABC}=[-kTdz_{2\uno\due\tre}]$. Further computations show that such class in non-zero in the quotient $H_{A}^{1,3}(X)/H_A^{0,2}(X)\cup\mathfrak{a}_{34}$.

For case $(ii)$, we take $\mathfrak{a}_{12}=[T^{-2}dz_{13\uno}], \mathfrak{a}_{23}=[\barT^{2}dz_{2\due\tre}], \mathfrak{a}_{34}=[dz_{\tre}]$. Then 
\begin{equation*}
	-T^{-2}dz_{13\uno}\wedge \barT^{2}dz_{2\due\tre} =\partial\bar{\partial}(k T^{-2}\barT^2dz_{12\uno\due})\,\,\text{ for a }k\in\C\setminus\{0\}, \quad  \barT^{2}dz_{2\due\tre}\wedge dz_{\tre}=0.
\end{equation*}
Hence $\langle\mathfrak{a}_{12},\mathfrak{a}_{23},\mathfrak{a}_{34}\rangle_{ABC}=[kT^{-2}\barT^2dz_{12\uno\due\tre}]$ and it's easy to check that this class is non-zero in the quotient $H_{A}^{2,3}(X)/H_A^{2,2}(X)\cup\mathfrak{a}_{34}$.

Finally, as observed in Remark \ref{remark}, case $(iii)$ has vanishing  triple ABC-Massey product since in this case the $\partial\bar{\partial}$-Lemma holds.

\subsection{Triple ABC-Massey products for $\g_2^{\alpha}$}
For $\alpha=0$, since the manifolds obtained are K\"ahler, we clearly have vanishing triple ABC-Massey products. For $\alpha_n>0$, we can proceed as for $\g_1$ since the complex $C_\Gamma^{\bullet,\bullet}$ in case  $b=(2k+1)\frac{\pi}{2\mathfrak{Re}A_n}$ is isomorphic to the one in case $(i)$ of $\g_1$ and $C_\Gamma^{\bullet,\bullet}$ in case  $b=\frac{k\pi}{\mathfrak{Re}A_n}$ is isomorphic to the one in case $(ii)$ of $\g_1$. Moreover if $b\neq k\frac{\pi}{2\mathfrak{Re}A_n}$ the $\partial$$\bar{\partial}$-Lemma holds, hence we have trivial ABC-Massey products.

\subsection{Triple ABC-Massey products for $\g_8$}

Cases \ref{caso1} and \ref{caso4} of Theorem \ref{teorema} satisfy the $\partial\bar{\partial}$-Lemma, hence thanks to Remark \ref{remark} triple ABC-Massey products vanish. 

For case \ref{caso2} we can consider Bott-Chern cohomology classes $\mathfrak{a}_{12}=[T^{-1}dz_{13}]$, $\mathfrak{a}_{23}=[\barT dz_{\due\tre}]$, $\mathfrak{a}_{34}=[dz_{\tre}]$. Then
\begin{equation*}
	T^{-1}dz_{13}\wedge\barT dz_{\due\tre}=\partial\bar{\partial}(k T^{-1}\barT dz_{1\due})\,\,\text{ for a }k\in\C\setminus\{0\}, \quad \barT dz_{\due\tre}\wedge dz_{\tre}=0
\end{equation*}
and is easy to check that $\langle\mathfrak{a}_{12},\mathfrak{a}_{23},\mathfrak{a}_{34}\rangle_{ABC}=[-kT^{-1}\barT dz_{1\due\tre}]$ is non-zero in the quotient $H_{A}^{1,2}(X)/H_A^{1,1}(X)\cup\mathfrak{a}_{34}$. 

Analogously, for case \ref{caso3} we can take Bott-Chern cohomology classes $\mathfrak{a}_{12}=[\barT^{-1}dz_{1\tre}]$, $\mathfrak{a}_{23}=[T dz_{3\due}]$, $\mathfrak{a}_{34}=[dz_3]$, compute their products
\begin{equation*}
	\barT^{-1}dz_{1\tre}\wedge Tdz_{3\due}=\partial\bar{\partial}(k \barT^{-1}T dz_{1\due})\,\,\text{ for a }k\in\C\setminus\{0\}, \quad \barT dz_{\due\tre}\wedge dz_{\tre}=0
\end{equation*}
and easily check that $\langle\mathfrak{a}_{12},\mathfrak{a}_{23},\mathfrak{a}_{34}\rangle_{ABC}=[k\barT^{-1}T dz_{13\due}]$ is non-zero in the quotient $H_{A}^{2,1}(X)/H_A^{1,1}(X)\cup\mathfrak{a}_{34}$. 

Cases \ref{caso4}, \ref{caso6} and \ref{caso7} of Theorem \ref{teorema} correspond respectively to cases $(iii)$, $(ii)$ and $ (i) $ of $\g_1$, hence we have non-trivial triple ABC-Massey products in cases \ref{caso6} and \ref{caso7} and strong formality (that is the $\partial\bar{\partial}$-Lemma) in case \ref{caso4}. Lastly, for case \ref{caso5} we take $\mathfrak{a}_{12}=[T^{-2}dz_{13\uno}]$, $\mathfrak{a}_{23}=[\barT^{-2}dz_{2\due\tre}]$ and $\mathfrak{a}_{34}=[dz_{\tre}]$, we compute the products
\begin{equation*}
	-T^{-2}dz_{13\uno}\wedge \barT^{-2}dz_{2\due\tre}=\partial\bar{\partial}(kT^{-2}\barT^2dz_{12\uno\due})\,\,\text{ for a }k\in\C\setminus\{0\}, \quad \barT^{-2}dz_{2\due\tre}\wedge dz_{\tre}=0
\end{equation*}
and we check that the triple ABC-Massey product $\langle\mathfrak{a}_{12},\mathfrak{a}_{23},\mathfrak{a}_{34}\rangle_{ABC}=[kT^{-2}\barT^2dz_{12\uno\due\tre}]$ is non-zero in the quotient $H_{A}^{2,3}(X)/H_A^{2,2}(X)\cup\mathfrak{a}_{34}$.

Summarizing our results, we have proved the following Proposition.
\begin{proposition}
	For a splitting type solvmanifold $X=\C^n\ltimes_\varphi\C^m/\Gamma$, strong formality and the $\partial\bar{\partial}$-Lemma property are equivalent conditions. If $X$ is a six-dimensional splitting Calabi-Yau manifold, these notions are also equivalent to the existence of a geometrically-Bott-Chern-formal metric $g$ on $X$ and to the notion of weak formality.
\end{proposition}

This result together with Corollary \ref{corollarymio} suggest that these properties could coincide also in general dimension, namely for the class of splitting type solvmanifolds of the form $\C^n\ltimes_\varphi \C^m/\Gamma$.

\appendix

\section{Tables}\label{tables}

This section summarizes all the results we have obtained concerning the double complexes $B^{\bullet,\bullet}_\Gamma$, $C^{\bullet,\bullet}_\Gamma$ and the corresponding cohomology groups. Explicit generators are given for each manifold (except for computations on the solvmanifolds with underlying Lie algebra $\g_1$, since they correspond to the completely solvable Nakamura manifolds, whose cohomologies are presented in \cite{kasuya2012techniquescomputationsdolbeaultcohomology}, \cite{AngellaKasuya_2017}. In this case we only give the dimensions of Dolbeault and Bott-Chern cohomology groups). 

Since we have $\bar{\partial}\equiv0$ on $B^{\bullet,\bullet}_\Gamma$ as observed in the previous sections, the isomorphism $B^{\bullet,\bullet}_\Gamma\simeq H_{\bar{\partial}}^{\bullet,\bullet}(B^{\bullet,\bullet}_\Gamma)$ holds. Then, to avoid redundancy, we will not write explicitly the Dolbeault cohomology groups, giving anyway the corresponding Hodge numbers.

For manifolds satisfying the $\partial\bar{\partial}$-Lemma we only give the generators of the double complex $B^{\bullet,\bullet}_\Gamma$ since we have the isomorphisms $B^{\bullet,\bullet}_\Gamma=C^{\bullet,\bullet}_\Gamma\simeq H_{\bar{\partial}}^{\bullet,\bullet}(C^{\bullet,\bullet}_\Gamma)\simeq H_{BC}^{\bullet,\bullet}(C^{\bullet,\bullet}_\Gamma)$.

To shorten the notation, we write $dz_{\bar{i}}$ to denote $d\bar{z}_i$ and we juxtapose indices to denote wedge products (for example $dz_{ij}$ denotes $dz_i\wedge dz_j$).

\subsection{\for{toc}{Table for $\g_1$}\except{toc}{}}\textbf{Table for $\g_1$}

\begin{table}[H]
	\begin{center}
		\begin{tabular}{c||cc||cc||cc}
			\hline
			$\g_1$&\multicolumn{2}{c||}{case $(i)$}&\multicolumn{2}{c||}{case $(ii)$}&\multicolumn{2}{c}{case $(iii)$}\\
			&$\bar{\partial}$&$BC$&$\bar{\partial}$&$BC$&$\bar{\partial}$&$BC$\\
		
			\hline
			$(1,0)$&3&1&1&1&1&1\\
			$(0,1)$&3&1&1&1&1&1\\
			\hline
			$(2,0)$&3&3&1&1&1&1\\
			$(1,1)$&9&7&5&3&3&3\\
			$(0,2)$&3&3&1&1&1&1\\
			\hline
			$(3,0)$&1&1&1&1&1&1\\
			$(2,1)$&9&9&5&5&3&3\\
			$(1,2)$&9&9&5&5&3&3\\
			$(0,3)$&1&1&1&1&1&1\\
			\hline
			$(3,1)$&3&3&1&1&1&1\\
			$(2,2)$&9&11&5&7&3&3\\
			$(1,3)$&3&3&1&1&1&1\\
			\hline
			$(3,2)$&3&5&1&1&1&1\\
			$(2,3)$&3&5&1&1&1&1\\
			\hline
			$(3,3)$&1&1&1&1&1&1\\
			\hline
			
		\end{tabular}
	\end{center}
	\caption{The dimensions of the Dolbeault and Bott-Chern cohomologies of the solvmanifolds $G_1/\Gamma$, depending on $\Gamma$, with underlying Lie algebra $\g_1$ (see \cite[Example Example 1]{kasuya2012techniquescomputationsdolbeaultcohomology} and \cite[Example 3.1]{AngellaKasuya_2017} for explicit generators).}
	\label{table1}
\end{table}

\newpage

\subsection{\for{toc}{Tables for $\g_2^{\alpha}$}\except{toc}{}}\textbf{Tables for $\g_2^\alpha$}

\begin{table}[H]
	{\footnotesize \begin{center}	
		\begin{tabular}{c||ll||ll}
			\hline
			$\g_2^0$&complex $B^{\bullet, \bullet}_\Gamma$ with $x_3=\frac{\pi}{2}\quad$ & $h^{\bullet,\bullet}_{\bar{\partial}}$& complex  $B^{\bullet, \bullet}_\Gamma$ with $ x_3\in\{\frac{\pi}{3}, \frac{\pi}{4}, \frac{\pi}{6}\} \quad$&$h^{\bullet,\bullet}_{\bar{\partial}}$\\
			\hline
			
			$(1,0)$&$ \C\langle dz_3\rangle $&1&$\C\langle dz_3\rangle$&1\\
			$(0,1)$&$ \C\langle dz_{\tre}\rangle $&1&$\C\langle dz_{\tre}\rangle$&1\\
			\hline
			$(2,0)$&$\C\langle dz_{12}\rangle$&1&$\C\langle dz_{12}\rangle$&1\\	
			$(1,1)$&$\C\langle dz_{1\uno}, dz_{1\due}, dz_{2\uno}, dz_{2\due}, dz_{3\tre}\rangle$&5&$\C\langle dz_{1\uno}, dz_{2\due}, dz_{3\tre}\rangle$&3\\
			$(0,2)$&$\C\langle dz_{\uno\due}\rangle$&1&$\C\langle dz_{\uno\due}\rangle$&1\\	
			\hline
			$(3,0)$&$\C\langle dz_{123}\rangle$&1&$\C\langle dz_{123}\rangle$&1\\
			$(2,1)$&$\C\langle dz_{12\tre}, dz_{13\uno}, dz_{13\due}, dz_{23\uno}, dz_{23\due}\rangle$&5&$\C\langle dz_{12\tre},  dz_{13\uno}, dz_{23\due}\rangle$&3\\
			$(1,2)$&$\C\langle dz_{1\uno\tre}, dz_{1\due\tre}, dz_{2\uno\tre}, dz_{2\due\tre}, dz_{3\uno\due}\rangle$&5&$\C\langle  dz_{1\uno\tre}, dz_{2\due\tre},  dz_{3\uno\due}\rangle$&3\\
			$(0,3)$&$\C\langle dz_{\uno\due\tre}\rangle$&1&$\C\langle dz_{\uno\due\tre}\rangle$&1\\
			\hline
			$(3,1)$&$\C\langle dz_{123\tre}\rangle$&1&$\C\langle dz_{123\tre}\rangle$&1\\
			$(2,2)$&$\C\langle dz_{12\uno\due}, dz_{13\uno\tre}, dz_{13\due\tre}, dz_{23\uno\tre}, dz_{23\due\tre}\rangle$&5&$\C\langle dz_{12\uno\due},  dz_{13\uno\tre}, dz_{23\due\tre}\rangle$&3\\
			$(1,3)$&$\C\langle dz_{3\uno\due\tre}\rangle$&1&$\C\langle dz_{3\uno\due\tre}\rangle$&1\\
			\hline
			$(3,2)$&$\C\langle dz_{123\uno\due}\rangle$&1&$\C\langle dz_{123\uno\due}\rangle$&1\\
			$(2,3)$&$\C\langle dz_{12\uno\due\tre}\rangle$&1&$\C\langle dz_{12\uno\due\tre}\rangle$&1\\
			\hline
			$(3,3)$&$\C\langle dz_{123\uno\due\tre}\rangle$&1&$\C\langle dz_{123\uno\due\tre}\rangle$&1\\
			\hline
		\end{tabular}
	\caption{The double complex $B_\Gamma^{\bullet,\bullet}$ for computing the Dolbeault cohomology of the solvmanifolds $G_2^0/\Gamma$ with underlying Lie algebra $\g_2^0$ and the corresponding Hodge numbers (computed in \cite[Table 4.2]{otalthesis}), depending on the value $x_3$ defining the lattice $\Gamma$. These manifolds satisfy the $\partial\bar{\partial}$-Lemma.}
	\label{tabellag02}
\end{center}}
\end{table}

\begin{landscape}
	
\begin{table}[H]
	{\footnotesize \begin{center}	
		\begin{tabularx}{22.8cm}{c || >{\setlength{\baselineskip}{1.15\baselineskip}}X c || >{\setlength{\baselineskip}{1.15\baselineskip}}X c || >{\setlength{\baselineskip}{1.15\baselineskip}}X c}
			\hline
			
			$\g_2^{\alpha_n}$&$B^{\bullet, \bullet}_\Gamma$, $\Gamma'=\frac{\pi}{2\mathfrak{Im}A_n}\Z\oplus\frac{i(2k+1)\pi}{2\mathfrak{Re}A_n}\Z$& $h^{\bullet,\bullet}_{\bar{\partial}}$&$B^{\bullet, \bullet}_\Gamma$, $\Gamma'=\frac{\pi}{2\mathfrak{Im}A_n}\Z\oplus\frac{2ik\pi}{2\mathfrak{Re}A_n}\Z$&$h^{\bullet,\bullet}_{\bar{\partial}}$&$B^{\bullet, \bullet}_\Gamma$, $\Gamma'=\frac{\pi}{2\mathfrak{Im}A_n}\Z\oplus ib\Z, b\neq \frac{k\pi}{2\mathfrak{Re}A_n}$&$h^{\bullet,\bullet}_{\bar{\partial}}$\\
			\hline
			

			$(1,0)$&$ \C\langle T^{-1}dz_1, Tdz_2, dz_3\rangle $&3&$\C\langle dz_3\rangle$&$1$&$\C\langle dz_3\rangle$&$1$\\
			
			$(0,1)$&$ \C\langle T^{-1}dz_{\uno}, Tdz_{\due}, dz_{\tre}\rangle $&3&$\C\langle dz_{\tre}\rangle$&$1$&$\C\langle dz_{\tre}\rangle$&$1$\\
			\hline
			
			$(2,0)$&$\C\langle dz_{12}, T^{-1}dz_{13}, Tdz_{23}\rangle$&3&$\C\langle dz_{12}\rangle$&$1$&$\C\langle dz_{12}\rangle$&$1$\\	
			
			$(1,1)$&$\C\langle T^{-2}dz_{1\uno}, dz_{1\due}, T^{-1}dz_{1\tre}, dz_{2\uno},T^2dz_{2\due},$ $ T dz_{2\tre}, T^{-1}dz_{3\uno}, Tdz_{3\due}, dz_{3\tre}\rangle  $&9&$\C\langle T^{-2}dz_{1\uno}, dz_{1\due}, dz_{2\uno}, T^2dz_{2\due}, dz_{3\tre}\rangle$&$5$&$\C\langle dz_{1\due}, dz_{2\uno}, dz_{3\tre}\rangle$&$3$\\
			
			$(0,2)$&$\C\langle dz_{\uno\due}, T^{-1}dz_{\uno\tre}, Tdz_{\due\tre} \rangle$&3&$\C\langle dz_{\uno\due}\rangle$&1&$\C\langle dz_{\uno\due}\rangle$&1\\	
			\hline
			
			$(3,0)$&$\C\langle dz_{123}\rangle$&1&$\C\langle dz_{123}\rangle$&1&$\C\langle dz_{123}\rangle$&1\\
			
			$(2,1)$&$\C\langle T^{-1}dz_{12\uno}, Tdz_{12\due}, dz_{12\tre}, T^{-2}dz_{13\uno}, $ $dz_{13\due}, T^{-1}dz_{13\tre}, dz_{23\uno}, T^2dz_{23\due}, Tdz_{23\tre} \rangle$&9&$\C\langle dz_{12\tre}, T^{-2}dz_{13\uno}, dz_{13\due}, dz_{23\uno}, T^2dz_{23\due}\rangle$&5&$\C\langle dz_{12\tre},  dz_{13\due}, dz_{23\uno}\rangle$&3\\
			
			$(1,2)$&$\C\langle T^{-1}dz_{1\uno\due}, T^{-2}dz_{1\uno\tre}, dz_{1\due\tre}, Tdz_{2\uno\due}, $ $dz_{2\uno\tre}, T^2dz_{2\due\tre}, dz_{3\uno\due}, T^{-1}dz_{3\uno\tre}, Tdz_{3\due\tre} \rangle$&9&$\C\langle T^{-2}dz_{1\uno\tre}, dz_{1\due\tre}, dz_{2\uno\tre}, T^2dz_{2\due\tre}, dz_{3\uno\due}\rangle$&5&$\C\langle  dz_{1\due\tre}, dz_{2\uno\tre},  dz_{3\uno\due}\rangle$&3\\
			
			$(0,3)$&$\C\langle dz_{\uno\due\tre}\rangle$&1&$\C\langle dz_{\uno\due\tre}\rangle$&1&$\C\langle dz_{\uno\due\tre}\rangle$&1\\
			\hline
			
			$(3,1)$&$\C\langle T^{-1}dz_{123\uno}, Tdz_{123\due}, dz_{123\tre} \rangle$&3&$\C\langle dz_{123\tre}\rangle$&1&$\C\langle dz_{123\tre}\rangle$&1\\
			
			$(2,2)$&$\C\langle dz_{12\uno\due}, T^{-1}dz_{12\uno\tre}, Tdz_{12\due\tre}, T^{-1}dz_{13\uno\due},$ $T^{-2}dz_{13\uno\tre}, dz_{13\due\tre}, Tdz_{23\uno\due}, dz_{23\uno\tre}, T^2dz_{23\due\tre} \rangle$&9&$\C\langle dz_{12\uno\due}, T^{-2}dz_{13\uno\tre}, dz_{13\due\tre}, dz_{23\uno\tre},$ $ T^2dz_{23\due\tre}\rangle$&5&$\C\langle dz_{12\uno\due},  dz_{13\due\tre}, dz_{23\uno\tre}\rangle$&3\\
			
			$(1,3)$&$\C\langle T^{-1}dz_1\uno\due\tre, Tdz_{2\uno\due\tre}, dz_{3\uno\due\tre} \rangle$&3&$\C\langle dz_{3\uno\due\tre}\rangle$&1&$\C\langle dz_{3\uno\due\tre}\rangle$&1\\
			\hline
			
			$(3,2)$&$\C\langle dz_{123\uno\due}, T^{-1}dz_{123\uno\tre}, Tdz_{123\due\tre} \rangle$&3&$\C\langle dz_{123\uno\due}\rangle$&1&$\C\langle dz_{123\uno\due}\rangle$&1\\
			
			$(2,3)$&$\C\langle dz_{12\uno\due\tre}, T^{-1}dz_{13\uno\due\tre}, Tdz_{23\uno\due\tre} \rangle$&3&$\C\langle dz_{12\uno\due\tre}\rangle$&1&$\C\langle dz_{12\uno\due\tre}\rangle$&1\\
			\hline
			
			$(3,3)$&$\C\langle dz_{123\uno\due\tre}\rangle$&1&$\C\langle dz_{123\uno\due\tre}\rangle$&1&$\C\langle dz_{123\uno\due\tre}\rangle$&1\\
			\hline

		\end{tabularx}
		\caption{The double complex $B_\Gamma^{\bullet,\bullet}$ for computing the Dolbeault cohomology of the solvmanifolds $G_2^{\alpha_n}/\Gamma$ with underlying Lie algebra $\g_2^{\alpha_n}, \alpha_n=\frac{1}{\pi}|\log(|\frac{n+\sqrt{n^2-4}}{2}|)|$ (here, $T=e^{(A_n+\bar{A}_n)z_3}$) and the corresponding Hodge numbers (computed in \cite[Table 4.2]{otalthesis}), depending on the lattice $\Gamma'$. The third case satisfies the $\partial\bar{\partial}$-Lemma.}
		\label{tabellaBg2}
	\end{center}}
\end{table}

\begin{table}[H]
	{\footnotesize \begin{center}	
		\begin{tabularx}{22.8cm}{c || >{\setlength{\baselineskip}{1.2\baselineskip}}X || >{\setlength{\baselineskip}{1.2\baselineskip}}X }
			\hline
				
			$\g_2^{\alpha_n}$&$C^{\bullet, \bullet}_\Gamma$, $\Gamma'=\frac{\pi}{2\mathfrak{Im}A_n}\Z\oplus\frac{i(2k+1)\pi}{2\mathfrak{Re}A_n}\Z$&$C^{\bullet, \bullet}_\Gamma$, $\Gamma'=\frac{\pi}{2\mathfrak{Im}A_n}\Z\oplus\frac{2ik\pi}{2\mathfrak{Re}A_n}\Z$\\
			\hline
				

			$(1,0)$&$ \C\langle T^{-1}dz_1, Tdz_2, dz_3, \barT^{-1}dz_1, \barT dz_2 \rangle $&$\C\langle dz_3\rangle$\\
				
			$(0,1)$&$ \C\langle T^{-1}dz_{\uno}, Tdz_{\due}, dz_{\tre}, \barT^{-1}dz_{\uno}, \barT dz_{\due} \rangle $&$\C\langle dz_{\tre}\rangle$\\
			\hline
			
			$(2,0)$&$\C\langle dz_{12}, T^{-1}dz_{13}, Tdz_{23}, \barT^{-1}dz_{13}, \barT dz_{23} \rangle$&$\C\langle dz_{12}\rangle$\\	
			
			$(1,1)$&$\C\langle T^{-2}dz_{1\uno}, dz_{1\due}, T^{-1}dz_{1\tre}, dz_{2\uno},T^2dz_{2\due}, T dz_{2\tre},T^{-1}dz_{3\uno}, Tdz_{3\due},dz_{3\tre},$ $   \barT^{-2}dz_{1\uno}, \barT^{-1} dz_{1\tre}, \barT^2dz_{2\due}, $ $\barT dz_{2\tre}, \barT^{-1}dz_{3\uno}, \barT dz_{3\due} \rangle  $&$\C\langle T^{-2}dz_{1\uno}, dz_{1\due}, dz_{2\uno}, T^2dz_{2\due}, dz_{3\tre}, \barT^{-2}dz_{1\uno}, \barT^2dz_{2\due}\rangle$\\
			
			$(0,2)$&$\C\langle dz_{\uno\due}, T^{-1}dz_{\uno\tre}, Tdz_{\due\tre}, \barT^{-1}dz_{\uno\tre}, \barT dz_{\due\tre} \rangle$&$\C\langle dz_{\uno\due}\rangle$\\	
			\hline
			
			$(3,0)$&$\C\langle dz_{123}\rangle$&$\C\langle dz_{123}\rangle$\\
			
			$(2,1)$&$\C\langle T^{-1}dz_{12\uno}, Tdz_{12\due}, dz_{12\tre}, T^{-2}dz_{13\uno}, dz_{13\due}, T^{-1}dz_{13\tre}, dz_{23\uno}, T^2dz_{23\due},$ $ Tdz_{23\tre}, \barT^{-1}dz_{12\uno}, \barT dz_{12\due} \barT^{-2}dz_{13\uno}, \barT^{-1}dz_{13\tre}, \barT^2dz_{23\due}, \barT dz_{23\tre} \rangle$&$\C\langle dz_{12\tre}, T^{-2}dz_{13\uno}, dz_{13\due}, dz_{23\uno}, T^2dz_{23\due}, \barT^{-2}dz_{13\uno}, \barT^2dz_{23\due} \rangle$\\
			
			$(1,2)$&$\C\langle T^{-1}dz_{1\uno\due}, T^{-2}dz_{1\uno\tre}, dz_{1\due\tre}, Tdz_{2\uno\due}, dz_{2\uno\tre}, T^2dz_{2\due\tre}, dz_{3\uno\due}, T^{-1}dz_{3\uno\tre},$ $ Tdz_{3\due\tre}, \barT^{-1}dz_{1\uno\due}, \barT^{-2}dz_{1\uno\tre}, \barT dz_{2\uno\due}, \barT^2dz_{2\due\tre}, \barT^{-1}dz_{3\uno\tre}, \barT dz_{3\due\tre} \rangle$&$\C\langle T^{-2}dz_{1\uno\tre}, dz_{1\due\tre}, dz_{2\uno\tre}, T^2dz_{2\due\tre}, dz_{3\uno\due}, \barT^{-2}dz_{1\uno\tre}, \barT^2dz_{2\due\tre} \rangle$\\
			
			$(0,3)$&$\C\langle dz_{\uno\due\tre}\rangle$&$\C\langle dz_{\uno\due\tre}\rangle$\\
			\hline
			
			$(3,1)$&$\C\langle T^{-1}dz_{123\uno}, Tdz_{123\due}, dz_{123\tre}, \barT^{-1}dz_{123\uno}, \barT dz_{123\due} \rangle$&$\C\langle dz_{123\tre}\rangle$\\
			
			$(2,2)$&$\C\langle dz_{12\uno\due}, T^{-1}dz_{12\uno\tre}, Tdz_{12\due\tre}, T^{-1}dz_{13\uno\due}, T^{-2}dz_{13\uno\tre}, dz_{13\due\tre}, Tdz_{23\uno\due}, dz_{23\uno\tre},$ $ T^2dz_{23\due\tre}, \barT^{-1}dz_{12\uno\tre}, \barT dz_{12\due\tre}, \barT^{-1}dz_{13\uno\due}, \barT^{-2}dz_{13\uno\tre}, \barT dz_{23\uno\due}, \barT^2dz_{23\due\tre} \rangle$&$\C\langle dz_{12\uno\due}, T^{-2}dz_{13\uno\tre}, dz_{13\due\tre}, dz_{23\uno\tre}, T^2dz_{23\due\tre}, \barT^{-2}dz_{13\uno\tre}, \barT^2dz_{23\due\tre} \rangle$\\
			
			$(1,3)$&$\C\langle T^{-1}dz_{1\uno\due\tre}, Tdz_{2\uno\due\tre}, dz_{3\uno\due\tre}, \barT^{-1}dz_{1\uno\due\tre}, \barT dz_{2\uno\due\tre} \rangle$&$\C\langle dz_{3\uno\due\tre}\rangle$\\
			\hline
			
			$(3,2)$&$\C\langle dz_{123\uno\due}, T^{-1}dz_{123\uno\tre}, Tdz_{123\due\tre}, \barT^{-1}dz_{123\uno\tre}, \barT dz_{123\due\tre} \rangle$&$\C\langle dz_{123\uno\due}\rangle$\\
			
			$(2,3)$&$\C\langle dz_{12\uno\due\tre}, T^{-1}dz_{13\uno\due\tre}, Tdz_{23\uno\due\tre}, \barT^{-1}dz_{13\uno\due\tre}, \barT dz_{23\uno\due\tre} \rangle$&$\C\langle dz_{12\uno\due\tre}\rangle$\\
			\hline
			
			$(3,3)$&$\C\langle dz_{123\uno\due\tre}\rangle$&$\C\langle dz_{123\uno\due\tre}\rangle$\\
			\hline

		\end{tabularx}
		\caption{The double complex $C_\Gamma^{\bullet,\bullet}=B_\Gamma^{\bullet,\bullet}+\bar{B}_\Gamma^{\bullet,\bullet}$ for computing the Bott-Chern cohomology of the solvmanifolds $G_2^{\alpha_n}/\Gamma$ with underlying Lie algebra $\g_2^{\alpha_n}, \alpha_n=\frac{1}{\pi}|\log(|\frac{n+\sqrt{n^2-4}}{2}|)|$ (here, $T=e^{(A_n+\bar{A}_n)z_3}$), depending on the lattice $\Gamma'$.}
		\label{tabellaCg2}
	\end{center}}
\end{table}

\begin{table}[H]
	{\footnotesize \begin{center}	
		\begin{tabularx}{22.8cm}{c || >{\setlength{\baselineskip}{1.2\baselineskip}}X c|| >{\setlength{\baselineskip}{1.2\baselineskip}}X c}
			\hline
			
			$\g_2^{\alpha_n}$&$H_{BC}^{\bullet, \bullet}(G^{\alpha_n}_2/\Gamma)$, $\Gamma'=\frac{\pi}{2\mathfrak{Im}A_n}\Z\oplus\frac{i(2k+1)\pi}{2\mathfrak{Re}A_n}\Z$&dim&$H_{BC}^{\bullet, \bullet}(G^{\alpha_n}_2/\Gamma)$, $\Gamma'=\frac{\pi}{2\mathfrak{Im}A_n}\Z\oplus\frac{2ik\pi}{2\mathfrak{Re}A_n}\Z$&dim\\
			\hline
			

			$(1,0)$&$ \C\langle [dz_3] \rangle $&1&$\C\langle [dz_3]\rangle$&1\\
			
			$(0,1)$&$ \C\langle [dz_{\tre}] \rangle $&1&$\C\langle [dz_{\tre}]\rangle$&1\\
			\hline
			
			$(2,0)$&$\C\langle [dz_{12}], [T^{-1}dz_{13}], [Tdz_{23}] \rangle$&3&$\C\langle [dz_{12}]\rangle$&1\\	
			
			$(1,1)$&$\C\langle [dz_{1\due}], [ dz_{2\uno}], [T^{-1}dz_{3\uno}], [Tdz_{3\due}], [dz_{3\tre}],   [\barT^{-1} dz_{1\tre}], [\barT dz_{2\tre}] \rangle  $&7&$\C\langle [ dz_{1\due}], [dz_{2\uno}], [dz_{3\tre}] \rangle$&3\\
			
			$(0,2)$&$\C\langle [dz_{\uno\due}], [\barT^{-1}dz_{\uno\tre}], [\barT dz_{\due\tre}] \rangle$&3&$\C\langle [dz_{\uno\due}]\rangle$&1\\	
			\hline
			
			$(3,0)$&$\C\langle [dz_{123}]\rangle$&1&$\C\langle [dz_{123}]\rangle$&1\\
			
			$(2,1)$&$\C\langle [dz_{12\tre}], [T^{-2}dz_{13\uno}], [dz_{13\due}], [T^{-1}dz_{13\tre}], [dz_{23\uno}], [T^2dz_{23\due}], [Tdz_{23\tre}],$ $ [\barT^{-1}dz_{13\tre}], [\barT dz_{23\tre}] \rangle$&9&$\C\langle [dz_{12\tre}], [T^{-2}dz_{13\uno}], [dz_{13\due}], [dz_{23\uno}], [T^2dz_{23\due}] \rangle$&5\\
			
			$(1,2)$&$\C\langle  [dz_{1\due\tre}],  [dz_{2\uno\tre}], [dz_{3\uno\due}], [T^{-1}dz_{3\uno\tre}], [Tdz_{3\due\tre}], [\barT^{-2}dz_{1\uno\tre}],  [\barT^2dz_{2\due\tre}],$ $ [\barT^{-1}dz_{3\uno\tre}], [\barT dz_{3\due\tre}] \rangle$&9&$\C\langle [ dz_{1\due\tre}], [dz_{2\uno\tre}], [dz_{3\uno\due}], [\barT^{-2}dz_{1\uno\tre}], [\barT^2dz_{2\due\tre}] \rangle$&5\\
			
			$(0,3)$&$\C\langle [dz_{\uno\due\tre}]\rangle$&1&$\C\langle [dz_{\uno\due\tre}]\rangle$&1\\
			\hline
			
			$(3,1)$&$\C\langle [T^{-1}dz_{123\uno}], [Tdz_{123\due}], [dz_{123\tre}] \rangle$&3&$\C\langle [dz_{123\tre}]\rangle$&1\\
			
			$(2,2)$&$\C\langle [dz_{12\uno\due}],  [T^{-1}dz_{13\uno\due}], [T^{-2}dz_{13\uno\tre}], [dz_{13\due\tre}], [Tdz_{23\uno\due}], [dz_{23\uno\tre}],$ $ [T^2dz_{23\due\tre}], [\barT^{-1}dz_{12\uno\tre}], [\barT dz_{12\due\tre}], [\barT^{-2}dz_{13\uno\tre}], [\barT^2dz_{23\due\tre}] \rangle$&11&$\C\langle [dz_{12\uno\due}], [T^{-2}dz_{13\uno\tre}], [dz_{13\due\tre}], [dz_{23\uno\tre}], [T^2dz_{23\due\tre}], [\barT^{-2}dz_{13\uno\tre}],$ $ [\barT^2dz_{23\due\tre}] \rangle$&7\\
			
			$(1,3)$&$\C\langle [dz_{3\uno\due\tre}], [\barT^{-1}dz_{1\uno\due\tre}], [\barT dz_{2\uno\due\tre}] \rangle$&3&$\C\langle [dz_{3\uno\due\tre}]\rangle$&1\\
			\hline
			
			$(3,2)$&$\C\langle [dz_{123\uno\due}], [T^{-1}dz_{123\uno\tre}], [Tdz_{123\due\tre}], [\barT^{-1}dz_{123\uno\tre}], [\barT dz_{123\due\tre}] \rangle$&5&$\C\langle [dz_{123\uno\due}]\rangle$&1\\
			
			$(2,3)$&$\C\langle [dz_{12\uno\due\tre}], [T^{-1}dz_{13\uno\due\tre}], [Tdz_{23\uno\due\tre}], [\barT^{-1}dz_{13\uno\due\tre}], [\barT dz_{23\uno\due\tre}] \rangle$&5&$\C\langle [dz_{12\uno\due\tre}]\rangle$&1\\
			\hline
			
			$(3,3)$&$\C\langle [dz_{123\uno\due\tre}]\rangle$&1&$\C\langle [dz_{123\uno\due\tre}]\rangle$&1\\
			\hline
			\end{tabularx}
		\caption{The Bott-Chern cohomology groups and the relative dimensions for the solvmanifolds $G_2^{\alpha_n}/\Gamma$ with underlying Lie algebra $\g_2^{\alpha_n}$ with $ \alpha_n=\frac{1}{\pi}|\log(|\frac{n+\sqrt{n^2-4}}{2}|)|$ (here, $T=e^{(A_n+\bar{A}_n)z_3}$), depending on the lattice $\Gamma'$.}
		\label{tabellaBCg2}
	\end{center}}
\end{table}
\end{landscape}

\subsection{\for{toc}{Tables for $\g_8$}\except{toc}{}}\textbf{Tables for $\g_8$}

\begin{table}[H]
	{\footnotesize \begin{center}	
			\begin{tabularx}{16cm}{c || l c || l c || l c}
				\hline
				
				$\g_8$&$B^{\bullet, \bullet}_\Gamma$, case \ref{caso1}& $h^{\bullet,\bullet}_{\bar{\partial}}$&$B^{\bullet, \bullet}_\Gamma$, case \ref{caso2}&$h^{\bullet,\bullet}_{\bar{\partial}}$&$B^{\bullet, \bullet}_\Gamma$, case \ref{caso3}&$h^{\bullet,\bullet}_{\bar{\partial}}$\\
				\hline
				

				$(1,0)$&$\C\langle dz_3\rangle$&$1$&$\C\langle T^{-1}dz_1, Tdz_2, dz_3\rangle$&3&$ \C\langle dz_3\rangle $&1\\
				
				$(0,1)$&$\C\langle dz_{\tre}\rangle$&$1$&$\C\langle dz_{\tre}\rangle$&$1$&$ \C\langle T^{-1}dz_{\uno}, Tdz_{\due}, dz_{\tre}\rangle $&3\\
				\hline
				
				$(2,0)$&$\C\langle dz_{12}\rangle$&1&$\C\langle dz_{12}, T^{-1}dz_{13}, Tdz_{23}\rangle$&3&$\C\langle dz_{12} \rangle$&1\\	
				
				$(1,1)$&$\C\langle dz_{3\tre}\rangle$&1&$\C\langle T^{-1}dz_{1\tre}, Tdz_{2\tre}, dz_{3\tre}\rangle$&3&$\C\langle T^{-1}dz_{3\uno}, Tdz_{3\due}, dz_{3\tre}\rangle  $&3\\
				
				$(0,2)$&$\C\langle dz_{\uno\due}\rangle$&1&$\C\langle dz_{\uno\due}\rangle$&1&$\C\langle dz_{\uno\due}, T^{-1}dz_{\uno\tre}, Tdz_{\due\tre} \rangle$&3\\	
				\hline
				
				$(3,0)$&$\C\langle dz_{123}\rangle$&1&$\C\langle dz_{123}\rangle$&1&$\C\langle dz_{123}\rangle$&1\\
				
				$(2,1)$&$\C\langle dz_{12\tre} \rangle$&1&$\C\langle dz_{12\tre}, T^{-1}dz_{13\tre}, Tdz_{23\tre} \rangle$&3&$\C\langle T^{-1}dz_{12\uno}, Tdz_{12\due}, dz_{12\tre} \rangle$&3\\
				
				$(1,2)$&$\C\langle  dz_{3\uno\due}\rangle$&1&$\C\langle T^{-1}dz_{1\uno\due}, Tdz_{2\uno\due}, dz_{3\uno\due}\rangle$&3&$\C\langle dz_{3\uno\due}, T^{-1}dz_{3\uno\tre}, Tdz_{3\due\tre} \rangle$&3\\
				
				$(0,3)$&$\C\langle dz_{\uno\due\tre}\rangle$&1&$\C\langle dz_{\uno\due\tre}\rangle$&1&$\C\langle dz_{\uno\due\tre}\rangle$&1\\
				\hline
				
				$(3,1)$&$\C\langle dz_{123\tre}\rangle$&1&$\C\langle dz_{123\tre}\rangle$&1&$\C\langle T^{-1}dz_{123\uno}, Tdz_{123\due}, dz_{123\tre} \rangle$&3\\
				
				$(2,2)$&$\C\langle dz_{12\uno\due},  \rangle$&1&$\C\langle dz_{12\uno\due}, T^{-1}dz_{13\uno\due}, Tdz_{23\uno\due}\rangle$&3&$\C\langle dz_{12\uno\due}, T^{-1}dz_{12\uno\tre}, Tdz_{12\due\tre} \rangle$&3\\
				
				$(1,3)$&$\C\langle dz_{3\uno\due\tre}\rangle$&1&$\C\langle T^{-1}dz_{1\uno\due\tre}, Tdz_{2\uno\due\tre}, dz_{3\uno\due\tre}\rangle$&3&$\C\langle  dz_{3\uno\due\tre} \rangle$&1\\
				\hline
				
				$(3,2)$&$\C\langle dz_{123\uno\due}\rangle$&1&$\C\langle dz_{123\uno\due}\rangle$&1&$\C\langle dz_{123\uno\due}, T^{-1}dz_{123\uno\tre}, Tdz_{123\due\tre} \rangle$&3\\
				
				$(2,3)$&$\C\langle dz_{12\uno\due\tre}\rangle$&1&$\C\langle dz_{12\uno\due\tre}, T^{-1}dz_{13\uno\due\tre}, Tdz_{23\uno\due\tre}\rangle$&3&$\C\langle dz_{12\uno\due\tre} \rangle$&1\\
				\hline
				
				$(3,3)$&$\C\langle dz_{123\uno\due\tre}\rangle$&1&$\C\langle dz_{123\uno\due\tre}\rangle$&1&$\C\langle dz_{123\uno\due\tre}\rangle$&1\\
				\hline

			\end{tabularx}
			\caption{The double complex $B_\Gamma^{\bullet,\bullet}$ for computing the Dolbeault cohomology of the solvmanifolds $G_8/\Gamma$ with underlying Lie algebra $\g_8$ in cases \ref{caso1}, \ref{caso2} (computed in \cite[Example 2]{kasuya2012techniquescomputationsdolbeaultcohomology}) and \ref{caso3} of Theorem \ref{teorema} (here, $T=e^{-(A+\bar{A}-2i)z_3}$ with $A$ as in \eqref{struttura8}) and the relative Hodge numbers. Case \ref{caso1} satisfies the $\partial\bar{\partial}$-Lemma.}
			\label{tabellaBg8123}
	\end{center}}
\end{table}

\begin{landscape}

\begin{table}[H]
	{\footnotesize \begin{center}	
			\begin{tabularx}{22.8cm}{c || >{\setlength{\baselineskip}{1.2\baselineskip}}X || >{\setlength{\baselineskip}{1.2\baselineskip}}X }
				\hline
				
				$\g_8$&$C^{\bullet, \bullet}_\Gamma$, case \ref{caso2}&$C^{\bullet, \bullet}_\Gamma$, case \ref{caso3}\\
				\hline
				

				$(1,0)$&$\C\langle T^{-1}dz_1, Tdz_2, dz_3\rangle$&$ \C\langle dz_3, \barT^{-1}dz_1, \barT dz_2 \rangle $\\
				
				$(0,1)$&$\C\langle dz_{\tre}, \barT^{-1}dz_{\uno}, \barT dz_{\due}\rangle$&$ \C\langle T^{-1}dz_{\uno}, Tdz_{\due}, dz_{\tre} \rangle $\\
				\hline
				
				$(2,0)$&$\C\langle dz_{12}, T^{-1}dz_{13}, Tdz_{23} \rangle$&$\C\langle dz_{12}, \barT^{-1}dz_{13}, \barT dz_{23} \rangle$\\	
				
				$(1,1)$&$\C\langle T^{-1}dz_{1\tre},  Tdz_{2\tre}, dz_{3\tre}, \barT^{-1}dz_{3\uno}, \barT dz_{3\due}\rangle$&$\C\langle T^{-1}dz_{3\uno}, Tdz_{3\due},dz_{3\tre},  \barT^{-1} dz_{1\tre},\barT dz_{2\tre} \rangle  $\\
				
				$(0,2)$&$\C\langle dz_{\uno\due}, \barT^{-1}dz_{\uno\tre}, \barT dz_{\due\tre} \rangle$&$\C\langle dz_{\uno\due}, T^{-1}dz_{\uno\tre}, Tdz_{\due\tre} \rangle$\\	
				\hline
				
				$(3,0)$&$\C\langle dz_{123}\rangle$&$\C\langle dz_{123}\rangle$\\
				
				$(2,1)$&$\C\langle dz_{12\tre}, T^{-1}dz_{13\tre}, Tdz_{23\tre}, \barT^{-1}dz_{12\uno}, \barT dz_{12\due} \rangle$&$\C\langle T^{-1}dz_{12\uno}, Tdz_{12\due}, dz_{12\tre}, \barT^{-1}dz_{13\tre}, \barT dz_{23\tre} \rangle$\\
				
				$(1,2)$&$\C\langle T^{-1}dz_{1\uno\due}, Tdz_{2\uno\due}, dz_{3\uno\due}, \barT^{-1}dz_{3\uno\tre}, \barT dz_{3\due\tre} \rangle$&$\C\langle dz_{3\uno\due}, T^{-1}dz_{3\uno\tre}, Tdz_{3\due\tre}, \barT^{-1}dz_{1\uno\due}, \barT dz_{2\uno\due} \rangle$\\
				
				$(0,3)$&$\C\langle dz_{\uno\due\tre}\rangle$&$\C\langle dz_{\uno\due\tre}\rangle$\\
				\hline
				
				$(3,1)$&$\C\langle dz_{123\tre}, \barT^{-1}dz_{123\uno}, \barT dz_{123\due} \rangle$&$\C\langle T^{-1}dz_{123\uno}, Tdz_{123\due}, dz_{123\tre} \rangle$\\
				
				$(2,2)$&$\C\langle dz_{12\uno\due}, T^{-1}dz_{13\uno\due}, Tdz_{23\uno\due}, \barT^{-1}dz_{12\uno\tre}, \barT dz_{12\due\tre} \rangle$&$\C\langle dz_{12\uno\due}, T^{-1}dz_{12\uno\tre}, Tdz_{12\due\tre}, \barT^{-1}dz_{13\uno\due}, \barT dz_{23\uno\due} \rangle$\\
				
				$(1,3)$&$\C\langle T^{-1}dz_{1\uno\due\tre}, Tdz_{2\uno\due\tre}, dz_{3\uno\due\tre}\rangle$&$\C\langle dz_{3\uno\due\tre}, \barT^{-1}dz_{1\uno\due\tre}, \barT dz_{2\uno\due\tre} \rangle$\\
				\hline
				
				$(3,2)$&$\C\langle dz_{123\uno\due}, \barT^{-1}dz_{123\uno\tre}, \barT dz_{123\due\tre} \rangle$&$\C\langle dz_{123\uno\due}, T^{-1}dz_{123\uno\tre}, Tdz_{123\due\tre} \rangle$\\
				
				$(2,3)$&$\C\langle dz_{12\uno\due\tre}, T^{-1}dz_{13\uno\due\tre}, Tdz_{23\uno\due\tre} \rangle$&$\C\langle dz_{12\uno\due\tre}, \barT^{-1}dz_{13\uno\due\tre}, \barT dz_{23\uno\due\tre} \rangle$\\
				\hline
				
				$(3,3)$&$\C\langle dz_{123\uno\due\tre}\rangle$&$\C\langle dz_{123\uno\due\tre}\rangle$\\
				\hline

			\end{tabularx}
			\caption{The double complex $C_\Gamma^{\bullet,\bullet}=B_\Gamma^{\bullet,\bullet}+\bar{B}_\Gamma^{\bullet,\bullet}$ for computing the Bott-Chern cohomology of the solvmanifolds $G_8/\Gamma$ with underlying Lie algebra $\g_8$ in cases \ref{caso2} (computed in \cite[Table 10]{AngellaKasuya_2017}) and \ref{caso3} of Theorem \ref{teorema} (here, $T=e^{-(A+\bar{A}-2i)z_3}$ with $A$ as in \eqref{struttura8}).}
			\label{tabellaCg823}
	\end{center}}
\end{table}

\begin{table}[H]
	{\footnotesize \begin{center}	
			\begin{tabularx}{22.8cm}{c || >{\setlength{\baselineskip}{1.2\baselineskip}}X c|| >{\setlength{\baselineskip}{1.2\baselineskip}}X c}
				\hline
				
				$\g_8$&$H_{BC}^{\bullet, \bullet}(G_8/\Gamma)$, case \ref{caso2}&dim&$H_{BC}^{\bullet, \bullet}(G_8/\Gamma)$, case \ref{caso3}&dim\\
				\hline
				

				$(1,0)$&$ \C\langle [dz_3] \rangle $&1&$\C\langle [dz_3]\rangle$&1\\
				
				$(0,1)$&$ \C\langle [dz_{\tre}] \rangle $&1&$\C\langle [dz_{\tre}]\rangle$&1\\
				\hline
				
				$(2,0)$&$\C\langle [dz_{12}], [T^{-1}dz_{13}], [Tdz_{23}] \rangle$&3&$\C\langle [dz_{12}] \rangle$&1\\	
				
				$(1,1)$&$\C\langle [dz_{3\tre}] \rangle$&1&$\C\langle [T^{-1}dz_{3\uno}], [Tdz_{3\due}], [dz_{3\tre}], [\barT^{-1} dz_{1\tre}], [\barT dz_{2\tre}] \rangle  $&5\\
				
				$(0,2)$&$\C\langle [dz_{\uno\due}] [\barT^{-1}dz_{\uno\tre}], [\barT dz_{\due\tre}] \rangle$&3&$\C\langle [dz_{\uno\due}] \rangle$&1\\	
				\hline
				
				$(3,0)$&$\C\langle [dz_{123}]\rangle$&1&$\C\langle [dz_{123}]\rangle$&1\\
				
				$(2,1)$&$\C\langle [dz_{12\tre}], [T^{-1}dz_{13\tre}], [Tdz_{23\tre}] \rangle$&3&$\C\langle [dz_{12\tre}], [\barT^{-1}dz_{13\tre}], [\barT dz_{23\tre}] \rangle$&3\\
				
				$(1,2)$&$\C\langle [dz_{3\uno\due}], [\barT^{-1}dz_{3\uno\tre}], [\barT dz_{3\due\tre}] \rangle$&3&$\C\langle [dz_{3\uno\due}], [T^{-1}dz_{3\uno\tre}], [Tdz_{3\due\tre}] \rangle$&3\\
				
				$(0,3)$&$\C\langle [dz_{\uno\due\tre}]\rangle$&1&$\C\langle [dz_{\uno\due\tre}]\rangle$&1\\
				\hline
				
				$(3,1)$&$\C\langle [dz_{123\tre}]\rangle$&1&$\C\langle [T^{-1}dz_{123\uno}], [Tdz_{123\due}], [dz_{123\tre}] \rangle$&3\\
				
				$(2,2)$&$\C\langle [dz_{12\uno\due}], [T^{-1}dz_{13\uno\due}], [Tdz_{23\uno\due}], [\barT^{-1}dz_{12\uno\tre}], [\barT dz_{12\due\tre}] \rangle$&5&$\C\langle [dz_{12\uno\due}] \rangle$&1\\
				
				$(1,3)$&$\C\langle [dz_{3\uno\due\tre}]\rangle$&1&$\C\langle [dz_{3\uno\due\tre}], [\barT^{-1}dz_{1\uno\due\tre}], [\barT dz_{2\uno\due\tre}] \rangle$&3\\
				\hline
				
				$(3,2)$&$\C\langle [dz_{123\uno\due}], [\barT^{-1}dz_{123\uno\tre}], [\barT dz_{123\due\tre}]\rangle$&3&$\C\langle [dz_{123\uno\due}], [T^{-1}dz_{123\uno\tre}], [Tdz_{123\due\tre}] \rangle$&3\\
				
				$(2,3)$&$\C\langle [dz_{12\uno\due\tre}], [T^{-1}dz_{13\uno\due\tre}], [Tdz_{23\uno\due\tre}] \rangle$&3&$\C\langle [dz_{12\uno\due\tre}], [\barT^{-1}dz_{13\uno\due\tre}], [\barT dz_{23\uno\due\tre}] \rangle$&3\\
				\hline
				
				$(3,3)$&$\C\langle [dz_{123\uno\due\tre}]\rangle$&1&$\C\langle [dz_{123\uno\due\tre}]\rangle$&1\\
				\hline
			\end{tabularx}
			\caption{The Bott-Chern cohomology groups and the relative dimensions for the solvmanifolds $G_8/\Gamma$ with underlying Lie algebra $\g_8$ in cases \ref{caso2} (computed in \cite[Table 11]{AngellaKasuya_2017}) and \ref{caso3} of Theorem \ref{teorema} (here, $T=e^{-(A+\bar{A}-2i)z_3}$ with $A$ as in \eqref{struttura8}).}
			\label{tabellaBCg823}
	\end{center}}
\end{table}
\end{landscape}

\begin{landscape}
	\begin{table}[H]
		{\footnotesize \begin{center}	
				\begin{tabularx}{22.8cm}{c || >{\setlength{\baselineskip}{1.15\baselineskip}}X c || >{\setlength{\baselineskip}{1.15\baselineskip}}X || >{\setlength{\baselineskip}{1.15\baselineskip}}X c}
					\hline
					
					$\g_8$&$B^{\bullet, \bullet}_\Gamma$, case \ref{caso5}& $h^{\bullet,\bullet}_{\bar{\partial}}$&$C^{\bullet, \bullet}_\Gamma$, case \ref{caso5}&$H_{BC}^{\bullet, \bullet}(G_8/\Gamma)$, case \ref{caso5}&dim\\
					\hline
					

					$(1,0)$&$\C\langle dz_3\rangle$&$1$&$\C\langle dz_3\rangle$&$ \C\langle [dz_3]\rangle $&1\\
					
					$(0,1)$&$\C\langle dz_{\tre}\rangle$&$1$&$\C\langle dz_{\tre}\rangle$&$ \C\langle [dz_{\tre}]\rangle $&1\\
					\hline
					
					$(2,0)$&$\C\langle dz_{12}\rangle$&$1$&$\C\langle dz_{12}\rangle$&$\C\langle [dz_{12}]\rangle$&1\\	
					
					$(1,1)$&$\C\langle T^{-2}dz_{1\uno}, T^2dz_{2\due}, dz_{3\tre}\rangle$&$3$&$\C\langle T^{-2}dz_{1\uno}, T^2dz_{2\due}, dz_{3\tre}, \barT^{-2}dz_{1\uno}, \barT^2dz_{2\due}  \rangle$&$\C\langle [dz_{3\tre}]\rangle  $&1\\
					
					$(0,2)$&$\C\langle dz_{\uno\due}\rangle$&1&$\C\langle dz_{\uno\due}\rangle$&$\C\langle [dz_{\uno\due}]\rangle$&1\\	
					\hline
					
					$(3,0)$&$\C\langle dz_{123}\rangle$&1&$\C\langle dz_{123}\rangle$&$\C\langle [dz_{123}]\rangle$&1\\
					
					$(2,1)$&$\C\langle dz_{12\tre},  T^{-2}dz_{13\uno}, T^2dz_{23\due}\rangle$&3&$\C\langle dz_{12\tre}, T^{-2}dz_{13\uno}, T^2dz_{23\due}, \barT^{-2}dz_{13\uno},$ $ \barT^2dz_{23\due}\rangle$&$\C\langle [dz_{12\tre}], [T^{-2}dz_{13\uno}], [T^2dz_{23\due}]\rangle$&3\\
					
					$(1,2)$&$\C\langle  T^{-2}dz_{1\uno\tre}, T^2dz_{2\due\tre},  dz_{3\uno\due}\rangle$&3&$\C\langle T^{-2}dz_{1\uno\tre}, T^2dz_{2\due\tre}, dz_{3\uno\due}, \barT^{-2}dz_{1\uno\tre}, $ $\barT^2dz_{2\due\tre}\rangle$&$\C\langle [dz_{3\uno\due}], [\barT^{-2}dz_{1\uno\tre}], [\barT^2dz_{2\due\tre}] \rangle$&3\\
					
					$(0,3)$&$\C\langle dz_{\uno\due\tre}\rangle$&1&$\C\langle dz_{\uno\due\tre}\rangle$&$\C\langle [dz_{\uno\due\tre}]\rangle$&1\\
					\hline
					
					$(3,1)$&$\C\langle dz_{123\tre}\rangle$&1&$\C\langle dz_{123\tre}\rangle$&$\C\langle [dz_{123\tre}] \rangle$&1\\
					
					$(2,2)$&$\C\langle dz_{12\uno\due},  T^{-2}dz_{13\uno\tre}, T^2dz_{23\due\tre}\rangle$&3&$\C\langle dz_{12\uno\due}, T^{-2}dz_{13\uno\tre}, T^2dz_{23\due\tre}, \barT^{-2}dz_{13\uno\tre}, $ $\barT^2dz_{23\due\tre}\rangle$&$\C\langle [dz_{12\uno\due}], [T^{-2}dz_{13\uno\tre}], [T^2dz_{23\due\tre}], [\barT^{-2}dz_{13\uno\tre}],$ $ [\barT^2dz_{23\due\tre}] \rangle$&5\\
					
					$(1,3)$&$\C\langle dz_{3\uno\due\tre}\rangle$&1&$\C\langle dz_{3\uno\due\tre}\rangle$&$\C\langle [dz_{3\uno\due\tre}] \rangle$&1\\
					\hline
					
					$(3,2)$&$\C\langle dz_{123\uno\due}\rangle$&1&$\C\langle dz_{123\uno\due}\rangle$&$\C\langle [dz_{123\uno\due}] \rangle$&1\\
					
					$(2,3)$&$\C\langle dz_{12\uno\due\tre}\rangle$&1&$\C\langle dz_{12\uno\due\tre}\rangle$&$\C\langle [dz_{12\uno\due\tre}] \rangle$&1\\
					\hline
					
					$(3,3)$&$\C\langle dz_{123\uno\due\tre}\rangle$&1&$\C\langle dz_{123\uno\due\tre}\rangle$&$\C\langle [dz_{123\uno\due\tre}]\rangle$&1\\
					\hline

				\end{tabularx}
				\caption{The double complexes $B_\Gamma^{\bullet,\bullet}$ and $C_\Gamma^{\bullet,\bullet}$ and Bott-Chern cohomology groups of the solvmanifolds $G_8/\Gamma$ with underlying Lie algebra $\g_8$ in case \ref{caso5} of Theorem \ref{teorema} (here, $T=e^{-(A+\bar{A}-2i)z_3}$ with $A$ as in \eqref{struttura8}), the relative Hodge numbers and Bott-Chern groups dimension.}
				\label{tabellag85}
		\end{center}}
	\end{table}

\end{landscape}

\begin{landscape}

\begin{table}[H]
		{\footnotesize\centering
			
			\begin{tabularx}{22.8cm}{c || l c || >{\setlength{\baselineskip}{1.15\baselineskip}}X c || >{\setlength{\baselineskip}{1.15\baselineskip}}X c}
				\hline
				
				$\g_8$&$B^{\bullet, \bullet}_\Gamma$, case \ref{caso4}& $h^{\bullet,\bullet}_{\bar{\partial}}$&$B^{\bullet, \bullet}_\Gamma$, case \ref{caso6}&$h^{\bullet,\bullet}_{\bar{\partial}}$&$B^{\bullet, \bullet}_\Gamma$, case \ref{caso7}&$h^{\bullet,\bullet}_{\bar{\partial}}$\\
				\hline
				

				$(1,0)$&$\C\langle dz_3\rangle$&$1$&$\C\langle dz_3\rangle$&$1$&$ \C\langle T^{-1}dz_1, Tdz_2, dz_3\rangle $&3\\
				
				$(0,1)$&$\C\langle dz_{\tre}\rangle$&$1$&$\C\langle dz_{\tre}\rangle$&$1$&$ \C\langle T^{-1}dz_{\uno}, Tdz_{\due}, dz_{\tre}\rangle $&3\\
				\hline
				
				$(2,0)$&$\C\langle dz_{12}\rangle$&$1$&$\C\langle dz_{12}\rangle$&$1$&$\C\langle dz_{12}, T^{-1}dz_{13}, Tdz_{23}\rangle$&3\\	
				
				$(1,1)$&$\C\langle dz_{1\due}, dz_{2\uno}, dz_{3\tre}\rangle$&$3$&$\C\langle T^{-2}dz_{1\uno}, dz_{1\due}, dz_{2\uno}, T^2dz_{2\due}, dz_{3\tre}\rangle$&$5$&$\C\langle T^{-2}dz_{1\uno}, dz_{1\due}, T^{-1}dz_{1\tre}, dz_{2\uno},T^2dz_{2\due}, T dz_{2\tre},$ $ T^{-1}dz_{3\uno}, Tdz_{3\due}, dz_{3\tre}\rangle  $&9\\
				
				$(0,2)$&$\C\langle dz_{\uno\due}\rangle$&1&$\C\langle dz_{\uno\due}\rangle$&1&$\C\langle dz_{\uno\due}, T^{-1}dz_{\uno\tre}, Tdz_{\due\tre} \rangle$&3\\	
				\hline
				
				$(3,0)$&$\C\langle dz_{123}\rangle$&1&$\C\langle dz_{123}\rangle$&1&$\C\langle dz_{123}\rangle$&1\\
				
				$(2,1)$&$\C\langle dz_{12\tre},  dz_{13\due}, dz_{23\uno}\rangle$&3&$\C\langle dz_{12\tre}, T^{-2}dz_{13\uno}, dz_{13\due}, dz_{23\uno}, T^2dz_{23\due}\rangle$&5&$\C\langle T^{-1}dz_{12\uno}, Tdz_{12\due}, dz_{12\tre}, T^{-2}dz_{13\uno}, dz_{13\due},$ $ T^{-1}dz_{13\tre}, dz_{23\uno}, T^2dz_{23\due}, Tdz_{23\tre} \rangle$&9\\
				
				$(1,2)$&$\C\langle  dz_{1\due\tre}, dz_{2\uno\tre},  dz_{3\uno\due}\rangle$&3&$\C\langle T^{-2}dz_{1\uno\tre}, dz_{1\due\tre}, dz_{2\uno\tre}, T^2dz_{2\due\tre}, dz_{3\uno\due}\rangle$&5&$\C\langle T^{-1}dz_{1\uno\due}, T^{-2}dz_{1\uno\tre}, dz_{1\due\tre}, Tdz_{2\uno\due}, dz_{2\uno\tre},$ $ T^2dz_{2\due\tre}, dz_{3\uno\due}, T^{-1}dz_{3\uno\tre}, Tdz_{3\due\tre} \rangle$&9\\
				
				$(0,3)$&$\C\langle dz_{\uno\due\tre}\rangle$&1&$\C\langle dz_{\uno\due\tre}\rangle$&1&$\C\langle dz_{\uno\due\tre}\rangle$&1\\
				\hline
				
				$(3,1)$&$\C\langle dz_{123\tre}\rangle$&1&$\C\langle dz_{123\tre}\rangle$&1&$\C\langle T^{-1}dz_{123\uno}, Tdz_{123\due}, dz_{123\tre} \rangle$&3\\
				
				$(2,2)$&$\C\langle dz_{12\uno\due},  dz_{13\due\tre}, dz_{23\uno\tre}\rangle$&3&$\C\langle dz_{12\uno\due}, T^{-2}dz_{13\uno\tre}, dz_{13\due\tre}, dz_{23\uno\tre},$ $ T^2dz_{23\due\tre}\rangle$&5&$\C\langle dz_{12\uno\due}, T^{-1}dz_{12\uno\tre}, Tdz_{12\due\tre}, T^{-1}dz_{13\uno\due}, T^{-2}dz_{13\uno\tre},$ $ dz_{13\due\tre}, Tdz_{23\uno\due}, dz_{23\uno\tre}, T^2dz_{23\due\tre} \rangle$&9\\
				
				$(1,3)$&$\C\langle dz_{3\uno\due\tre}\rangle$&1&$\C\langle dz_{3\uno\due\tre}\rangle$&1&$\C\langle T^{-1}dz_1\uno\due\tre, Tdz_{2\uno\due\tre}, dz_{3\uno\due\tre} \rangle$&3\\
				\hline
				
				$(3,2)$&$\C\langle dz_{123\uno\due}\rangle$&1&$\C\langle dz_{123\uno\due}\rangle$&1&$\C\langle dz_{123\uno\due}, T^{-1}dz_{123\uno\tre}, Tdz_{123\due\tre} \rangle$&3\\
				
				$(2,3)$&$\C\langle dz_{12\uno\due\tre}\rangle$&1&$\C\langle dz_{12\uno\due\tre}\rangle$&1&$\C\langle dz_{12\uno\due\tre}, T^{-1}dz_{13\uno\due\tre}, Tdz_{23\uno\due\tre} \rangle$&3\\
				\hline
				
				$(3,3)$&$\C\langle dz_{123\uno\due\tre}\rangle$&1&$\C\langle dz_{123\uno\due\tre}\rangle$&1&$\C\langle dz_{123\uno\due\tre}\rangle$&1\\
				\hline

			\end{tabularx}
			\caption{The double complex $B_\Gamma^{\bullet,\bullet}$ for computing the Dolbeault cohomology of the solvmanifolds $G_8/\Gamma$ with underlying Lie algebra $\g_8$ in cases \ref{caso4}, \ref{caso6} and \ref{caso7} of Theorem \ref{teorema} (here, $T=e^{-(A+\bar{A}-2i)z_3}$ with $A$ as in \eqref{struttura8}) and the corresponding Hodge numbers, computed in \cite[Table 4.5]{otalthesis}. Case \ref{caso4} satisfies the $\partial\bar{\partial}$-Lemma.}
			\label{tabellaBg8Otal}}
	
\end{table}

\end{landscape}

\begin{landscape}

\begin{table}[H]
	{\footnotesize \begin{center}	
			\begin{tabularx}{22.8cm}{c || >{\setlength{\baselineskip}{1.2\baselineskip}}X || >{\setlength{\baselineskip}{1.2\baselineskip}}X }
				\hline
				
				$\g_8$&$C^{\bullet, \bullet}_\Gamma$, case \ref{caso6}&$C^{\bullet, \bullet}_\Gamma$, case \ref{caso7}\\
				\hline
				

				$(1,0)$&$\C\langle dz_3\rangle$&$ \C\langle T^{-1}dz_1, Tdz_2, dz_3, \barT^{-1}dz_1, \barT dz_2 \rangle $\\
				
				$(0,1)$&$\C\langle dz_{\tre}\rangle$&$ \C\langle T^{-1}dz_{\uno}, Tdz_{\due}, dz_{\tre}, \barT^{-1}dz_{\uno}, \barT dz_{\due} \rangle $\\
				\hline
				
				$(2,0)$&$\C\langle dz_{12}\rangle$&$\C\langle dz_{12}, T^{-1}dz_{13}, Tdz_{23}, \barT^{-1}dz_{13}, \barT dz_{23} \rangle$\\	
				
				$(1,1)$&$\C\langle T^{-2}dz_{1\uno}, dz_{1\due}, dz_{2\uno}, T^2dz_{2\due}, dz_{3\tre}, \barT^{-2}dz_{1\uno}, \barT^2dz_{2\due}\rangle$&$\C\langle T^{-2}dz_{1\uno}, dz_{1\due}, T^{-1}dz_{1\tre}, dz_{2\uno},T^2dz_{2\due}, T dz_{2\tre},T^{-1}dz_{3\uno}, Tdz_{3\due},dz_{3\tre},$ $   \barT^{-2}dz_{1\uno}, \barT^{-1} dz_{1\tre}, \barT^2dz_{2\due}, $ $\barT dz_{2\tre}, \barT^{-1}dz_{3\uno}, \barT dz_{3\due} \rangle  $\\
				
				$(0,2)$&$\C\langle dz_{\uno\due}\rangle$&$\C\langle dz_{\uno\due}, T^{-1}dz_{\uno\tre}, Tdz_{\due\tre}, \barT^{-1}dz_{\uno\tre}, \barT dz_{\due\tre} \rangle$\\	
				\hline
				
				$(3,0)$&$\C\langle dz_{123}\rangle$&$\C\langle dz_{123}\rangle$\\
				
				$(2,1)$&$\C\langle dz_{12\tre}, T^{-2}dz_{13\uno}, dz_{13\due}, dz_{23\uno}, T^2dz_{23\due}, \barT^{-2}dz_{13\uno}, \barT^2dz_{23\due} \rangle$&$\C\langle T^{-1}dz_{12\uno}, Tdz_{12\due}, dz_{12\tre}, T^{-2}dz_{13\uno}, dz_{13\due}, T^{-1}dz_{13\tre}, dz_{23\uno}, T^2dz_{23\due},$ $ Tdz_{23\tre}, \barT^{-1}dz_{12\uno}, \barT dz_{12\due} \barT^{-2}dz_{13\uno}, \barT^{-1}dz_{13\tre}, \barT^2dz_{23\due}, \barT dz_{23\tre} \rangle$\\
				
				$(1,2)$&$\C\langle T^{-2}dz_{1\uno\tre}, dz_{1\due\tre}, dz_{2\uno\tre}, T^2dz_{2\due\tre}, dz_{3\uno\due}, \barT^{-2}dz_{1\uno\tre}, \barT^2dz_{2\due\tre} \rangle$&$\C\langle T^{-1}dz_{1\uno\due}, T^{-2}dz_{1\uno\tre}, dz_{1\due\tre}, Tdz_{2\uno\due}, dz_{2\uno\tre}, T^2dz_{2\due\tre}, dz_{3\uno\due}, T^{-1}dz_{3\uno\tre},$ $ Tdz_{3\due\tre}, \barT^{-1}dz_{1\uno\due}, \barT^{-2}dz_{1\uno\tre}, \barT dz_{2\uno\due}, \barT^2dz_{2\due\tre}, \barT^{-1}dz_{3\uno\tre}, \barT dz_{3\due\tre} \rangle$\\
				
				$(0,3)$&$\C\langle dz_{\uno\due\tre}\rangle$&$\C\langle dz_{\uno\due\tre}\rangle$\\
				\hline
				
				$(3,1)$&$\C\langle dz_{123\tre}\rangle$&$\C\langle T^{-1}dz_{123\uno}, Tdz_{123\due}, dz_{123\tre}, \barT^{-1}dz_{123\uno}, \barT dz_{123\due} \rangle$\\
				
				$(2,2)$&$\C\langle dz_{12\uno\due}, T^{-2}dz_{13\uno\tre}, dz_{13\due\tre}, dz_{23\uno\tre}, T^2dz_{23\due\tre}, \barT^{-2}dz_{13\uno\tre}, \barT^2dz_{23\due\tre} \rangle$&$\C\langle dz_{12\uno\due}, T^{-1}dz_{12\uno\tre}, Tdz_{12\due\tre}, T^{-1}dz_{13\uno\due}, T^{-2}dz_{13\uno\tre}, dz_{13\due\tre}, Tdz_{23\uno\due}, dz_{23\uno\tre},$ $ T^2dz_{23\due\tre}, \barT^{-1}dz_{12\uno\tre}, \barT dz_{12\due\tre}, \barT^{-1}dz_{13\uno\due}, \barT^{-2}dz_{13\uno\tre}, \barT dz_{23\uno\due}, \barT^2dz_{23\due\tre} \rangle$\\
				
				$(1,3)$&$\C\langle dz_{3\uno\due\tre}\rangle$&$\C\langle T^{-1}dz_{1\uno\due\tre}, Tdz_{2\uno\due\tre}, dz_{3\uno\due\tre}, \barT^{-1}dz_{1\uno\due\tre}, \barT dz_{2\uno\due\tre} \rangle$\\
				\hline
				
				$(3,2)$&$\C\langle dz_{123\uno\due}\rangle$&$\C\langle dz_{123\uno\due}, T^{-1}dz_{123\uno\tre}, Tdz_{123\due\tre}, \barT^{-1}dz_{123\uno\tre}, \barT dz_{123\due\tre} \rangle$\\
				
				$(2,3)$&$\C\langle dz_{12\uno\due\tre}\rangle$&$\C\langle dz_{12\uno\due\tre}, T^{-1}dz_{13\uno\due\tre}, Tdz_{23\uno\due\tre}, \barT^{-1}dz_{13\uno\due\tre}, \barT dz_{23\uno\due\tre} \rangle$\\
				\hline
				
				$(3,3)$&$\C\langle dz_{123\uno\due\tre}\rangle$&$\C\langle dz_{123\uno\due\tre}\rangle$\\
				\hline

			\end{tabularx}
			\caption{The double complex $C_\Gamma^{\bullet,\bullet}=B_\Gamma^{\bullet,\bullet}+\bar{B}_\Gamma^{\bullet,\bullet}$ for computing the Bott-Chern cohomology of the solvmanifolds $G_8/\Gamma$ with underlying Lie algebra $\g_8$ in case \ref{caso6} and in case \ref{caso7} (computed in \cite[Table 7]{AngellaKasuya_2017}) of Theorem \ref{teorema} (here, $T=e^{-(A+\bar{A}-2i)z_3}$ with $A$ as in \eqref{struttura8}).}
			\label{tabellaCg8Otal}
	\end{center}}
\end{table}

\begin{table}[H]
	{\footnotesize \begin{center}	
			\begin{tabularx}{22.8cm}{c || >{\setlength{\baselineskip}{1.2\baselineskip}}X c|| >{\setlength{\baselineskip}{1.2\baselineskip}}X c}
				\hline
				
				$\g_8$&$H_{BC}^{\bullet, \bullet}(G_8/\Gamma)$, case \ref{caso6}&dim&$H_{BC}^{\bullet, \bullet}(G_8/\Gamma)$, case \ref{caso7}&dim\\
				\hline
				

				$(1,0)$&$ \C\langle [dz_3] \rangle $&1&$\C\langle [dz_3]\rangle$&1\\
				
				$(0,1)$&$ \C\langle [dz_{\tre}] \rangle $&1&$\C\langle [dz_{\tre}]\rangle$&1\\
				\hline
				
				$(2,0)$&$\C\langle [dz_{12}]\rangle$&1&$\C\langle [dz_{12}], [T^{-1}dz_{13}], [Tdz_{23}] \rangle$&3\\	
				
				$(1,1)$&$\C\langle [ dz_{1\due}], [dz_{2\uno}], [dz_{3\tre}] \rangle$&3&$\C\langle [dz_{1\due}], [ dz_{2\uno}], [T^{-1}dz_{3\uno}], [Tdz_{3\due}], [dz_{3\tre}],   [\barT^{-1} dz_{1\tre}], [\barT dz_{2\tre}] \rangle  $&7\\
				
				$(0,2)$&$\C\langle [dz_{\uno\due}]\rangle$&1&$\C\langle [dz_{\uno\due}], [\barT^{-1}dz_{\uno\tre}], [\barT dz_{\due\tre}] \rangle$&3\\	
				\hline
				
				$(3,0)$&$\C\langle [dz_{123}]\rangle$&1&$\C\langle [dz_{123}]\rangle$&1\\
				
				$(2,1)$&$\C\langle [dz_{12\tre}], [T^{-2}dz_{13\uno}], [dz_{13\due}], [dz_{23\uno}], [T^2dz_{23\due}] \rangle$&5&$\C\langle [dz_{12\tre}], [T^{-2}dz_{13\uno}], [dz_{13\due}], [T^{-1}dz_{13\tre}], [dz_{23\uno}], [T^2dz_{23\due}], [Tdz_{23\tre}],$ $ [\barT^{-1}dz_{13\tre}], [\barT dz_{23\tre}] \rangle$&9\\
				
				$(1,2)$&$\C\langle [ dz_{1\due\tre}], [dz_{2\uno\tre}], [dz_{3\uno\due}], [\barT^{-2}dz_{1\uno\tre}], [\barT^2dz_{2\due\tre}] \rangle$&5&$\C\langle  [dz_{1\due\tre}],  [dz_{2\uno\tre}], [dz_{3\uno\due}], [T^{-1}dz_{3\uno\tre}], [Tdz_{3\due\tre}], [\barT^{-2}dz_{1\uno\tre}],  [\barT^2dz_{2\due\tre}],$ $ [\barT^{-1}dz_{3\uno\tre}], [\barT dz_{3\due\tre}] \rangle$&9\\
				
				$(0,3)$&$\C\langle [dz_{\uno\due\tre}]\rangle$&1&$\C\langle [dz_{\uno\due\tre}]\rangle$&1\\
				\hline
				
				$(3,1)$&$\C\langle [dz_{123\tre}]\rangle$&1&$\C\langle [T^{-1}dz_{123\uno}], [Tdz_{123\due}], [dz_{123\tre}] \rangle$&3\\
				
				$(2,2)$&$\C\langle [dz_{12\uno\due}], [T^{-2}dz_{13\uno\tre}], [dz_{13\due\tre}], [dz_{23\uno\tre}], [T^2dz_{23\due\tre}], [\barT^{-2}dz_{13\uno\tre}],$ $ [\barT^2dz_{23\due\tre}] \rangle$&7&$\C\langle [dz_{12\uno\due}],  [T^{-1}dz_{13\uno\due}], [T^{-2}dz_{13\uno\tre}], [dz_{13\due\tre}], [Tdz_{23\uno\due}], [dz_{23\uno\tre}],$ $ [T^2dz_{23\due\tre}], [\barT^{-1}dz_{12\uno\tre}], [\barT dz_{12\due\tre}], [\barT^{-2}dz_{13\uno\tre}], [\barT^2dz_{23\due\tre}] \rangle$&11\\
				
				$(1,3)$&$\C\langle [dz_{3\uno\due\tre}]\rangle$&1&$\C\langle [dz_{3\uno\due\tre}], [\barT^{-1}dz_{1\uno\due\tre}], [\barT dz_{2\uno\due\tre}] \rangle$&3\\
				\hline
				
				$(3,2)$&$\C\langle [dz_{123\uno\due}]\rangle$&1&$\C\langle [dz_{123\uno\due}], [T^{-1}dz_{123\uno\tre}], [Tdz_{123\due\tre}], [\barT^{-1}dz_{123\uno\tre}], [\barT dz_{123\due\tre}] \rangle$&5\\
				
				$(2,3)$&$\C\langle [dz_{12\uno\due\tre}]\rangle$&1&$\C\langle [dz_{12\uno\due\tre}], [T^{-1}dz_{13\uno\due\tre}], [Tdz_{23\uno\due\tre}], [\barT^{-1}dz_{13\uno\due\tre}], [\barT dz_{23\uno\due\tre}] \rangle$&5\\
				\hline
				
				$(3,3)$&$\C\langle [dz_{123\uno\due\tre}]\rangle$&1&$\C\langle [dz_{123\uno\due\tre}]\rangle$&1\\
				\hline
			\end{tabularx}
			\caption{The Bott-Chern cohomology groups and the corresponding dimensions for the solvmanifolds $G_8/\Gamma$ with underlying Lie algebra $\g_8$ in case \ref{caso6} and in case \ref{caso7} (computed in \cite[Table 8]{AngellaKasuya_2017}) of Theorem \ref{teorema} (here, $T=e^{-(A+\bar{A}-2i)z_3}$ with $A$ as in \eqref{struttura8}).}
			\label{tabellaBCg8Otal}
	\end{center}}
\end{table}
\end{landscape}

\section{Double complexes decomposition}\label{decompositions}

This section contains the decompositions in irreducible double-complexes for $C_\Gamma^{\bullet,\bullet}$. For the explicit decomposition, we denote a dot in position $(p,q)$ by $D^{p,q}$, a zigzag with shape $\{(p,q),(p+1,q)\}$ (horizontal lines starting from the position $(p,q)$) by $S_h^{p,q}$ and a zigzag with shape $\{(p,q),(p,q+1)\}$ (vertical lines starting from the position $(p,q)$) by $S_v^{p,q}$. We divide the decompositions in dots, horizontal lines and vertical lines to lighten the direct sums. In the graphic representations we only give the dimensions of the involved vector spaces, melting together zigzags with the same shape (for example, we write $\C^2\to\C^2$ to represent the direct sum of two zigzags $\C\to\C$).

\subsection{\for{toc}{Decompositions for $\g_1$}\except{toc}{}}\textbf{Decompositions for $\g_1$:} distinguishing the three cases involved, we have:
\begin{description}[itemsep=5pt]
	\item[case $(i)$] \begin{description}[itemsep=5pt]
		\item[dots] $\bigoplus_{j=0}^3 D^{j,0} \oplus \bigoplus_{j=0}^3 D^{0,j}\oplus (D^{1,1})^{\oplus 3}\oplus (D^{2,1})^{\oplus 3}\oplus (D^{1,2})^{\oplus 3}\oplus (D^{2,2})^{\oplus 3}$;
		\item[horizontal lines] $(S_h^{1,0})^{\oplus 2}\oplus (S_h^{0,1})^{\oplus 2}\oplus (S_h^{1,1})^{\oplus 4}\oplus (S_h^{0,2})^{\oplus 2}\oplus (S_h^{2,1})^{\oplus 2}\oplus (S_h^{1,2})^{\oplus 4}\oplus (S_h^{2,2})^{\oplus 2}\oplus (S_h^{1,3})^{\oplus 2}$;
		\item[vertical lines] $(S_v^{1,0})^{\oplus 2}\oplus (S_v^{0,1})^{\oplus 2}\oplus (S_v^{2,0})^{\oplus 2}\oplus (S_v^{1,1})^{\oplus 4}\oplus (S_v^{2,1})^{\oplus 4}\oplus (S_v^{1,2})^{\oplus 2}\oplus (S_v^{3,1})^{\oplus 2}\oplus (S_v^{2,2})^{\oplus 2}$.
		
	\end{description}
	
	\item[case $(ii)$] \begin{description}[itemsep=5pt]
		\item[dots] $\bigoplus_{j=0}^3 D^{j,0} \oplus \bigoplus_{j=0}^3 D^{0,j}\oplus (D^{1,1})^{\oplus 3}\oplus (D^{2,1})^{\oplus 3}\oplus (D^{1,2})^{\oplus 3}\oplus (D^{2,2})^{\oplus 3}$;
		
		\item[horizontal lines] $(S_h^{1,1})^{\oplus 2}\oplus (S_h^{1,2})^{\oplus 2}$;
		
		\item[vertical lines] $(S_v^{1,1})^{\oplus 2}\oplus (S_v^{2,1})^{\oplus 2}$.
		
	\end{description}
	
	\item[case $(iii)$] \begin{description}[itemsep=5pt]
		\item[dots] $\bigoplus_{j=0}^3 D^{j,0} \oplus \bigoplus_{j=0}^3 D^{0,j}\oplus (D^{1,1})^{\oplus 3}\oplus (D^{2,1})^{\oplus 3}\oplus (D^{1,2})^{\oplus 3}\oplus (D^{2,2})^{\oplus 3}$;
		
		\item[horizontal lines] $\emptyset$;
		
		\item[vertical lines] $\emptyset$.
		
	\end{description}
\end{description}

\begin{figure}[H]
	
	\begin{minipage}{.45\textwidth}
		
	\centering
	\captionsetup{width=\linewidth}
\scalebox{.7}{\begin{tikzpicture}
	\draw [step=2cm] (0,0) grid (8,8);
	\draw[thick,->](0,0) -- (8.25,0);
	\draw[thick,->](0,0) -- (0,8.25);
	\node at (8.5,0) {$q$};
	\node at (0,8.5) {$p$};
	\node at (-0.3,1) {$0$};
	\node at (-0.3,3) {$1$};
	\node at (-0.3,5) {$2$};
	\node at (-0.3,7) {$3$};
	\node at (1,-0.3) {$0$};
	\node at (3,-0.3) {$1$};
	\node at (5,-0.3) {$2$};
	\node at (7,-0.3) {$3$};
	
	\node at (1,1) {$\C$};
	
	\node at (2.5,0.5) {$\C$};
	\node at (3.5,1) {$\C^2$};
	\node at (3,1.5) {$\C^2$};
	\draw[thin, ->] (3.75,1) -- (4.25,1);
	\draw[thin, ->] (3,1.75) -- (3,2.15);
	
	\node at (0.5,2.5) {$\C$};
	\node at (1,3.5) {$\C^2$};
	\node at (1.5,3) {$\C^2$};
	\draw[thin, ->] (1,3.75) -- (1,4.25);
	\draw[thin, ->] (1.75,3) -- (2.15,3);
	
	\node at (5.5,0.5) {$\C$};
	\node at (4.5,1) {$\C^2$};
	\node at (5,1.5) {$\C^2$};
	\draw[thin, ->] (5,1.75) -- (5,2.15);
	
	\node at (3,3) {$\C^3$};
	\node at (3,2.40) {$\C^2$};
	\node at (2.4,3) {$\C^2$};
	\node at (3,3.60) {$\C^4$};
	\node at (3.60,3) {$\C^4$};
	\draw[thin,->] (3,3.85) -- (3,4.15);
	\draw[thin,->] (3.85,3) -- (4.15,3);
	
	\node at (0.5,5.5) {$\C$};
	\node at (1,4.5) {$\C^2$};
	\node at (1.5,5) {$\C^2$};
	\draw[thin, ->] (1.75,5) -- (2.15,5);
	
	\node at (7,1) {$\C$};
	
	\node at (5,3) {$\C^3$};
	\node at (5,2.40) {$\C^2$};
	\node at (4.4,3) {$\C^4$};
	\node at (5,3.60) {$\C^4$};
	\node at (5.60,3) {$\C^2$};
	\draw[thin,->] (5,3.85) -- (5,4.15);
	\draw[thin,->] (5.85,3) -- (6.25,3);
	
	\node at (3,5) {$\C^3$};
	\node at (3,4.40) {$\C^4$};
	\node at (2.4,5) {$\C^2$};
	\node at (3,5.60) {$\C^2$};
	\node at (3.60,5) {$\C^4$};
	\draw[thin,->] (3,5.85) -- (3,6.25);
	\draw[thin,->] (3.85,5) -- (4.15,5);
	
	\node at (1,7) {$\C$};
	
	\node at (7.5,2.5) {$\C$};
	\node at (7,3.5) {$\C^2$};
	\node at (6.5,3) {$\C^2$};
	\draw[thin, ->] (7,3.75) -- (7,4.25);
	
	\node at (5,5) {$\C^3$};
	\node at (5,4.40) {$\C^4$};
	\node at (4.4,5) {$\C^4$};
	\node at (5,5.60) {$\C^2$};
	\node at (5.60,5) {$\C^2$};
	\draw[thin,->] (5,5.85) -- (5,6.25);
	\draw[thin,->] (5.85,5) -- (6.25,5);
	
	\node at (2.5,7.5) {$\C$};
	\node at (3.5,7) {$\C^2$};
	\node at (3,6.5) {$\C^2$};
	\draw[thin, ->] (3.75,7) -- (4.25,7);
	
	\node at (7.5,5.5) {$\C$};
	\node at (7,4.5) {$\C^2$};
	\node at (6.5,5) {$\C^2$};
	
	\node at (5.5,7.5) {$\C$};
	\node at (5,6.5) {$\C^2$};
	\node at (4.5,7) {$\C^2$};
	
	\node at (7,7) {$\C$};

\end{tikzpicture}
}\caption{Graphic representation of the decomposition stated in Theorem \ref{decomposition} related to the solvmanifold $G_1/\Gamma$ with underlying Lie algebra $\g_1$ in case $ (i) $.}
\label{figure1}
\end{minipage}
\hspace{1.2cm}
\begin{minipage}{.45\textwidth}
		\centering
		\captionsetup{width=\linewidth}
	\scalebox{.7}{\begin{tikzpicture}
	\draw [step=2cm] (0,0) grid (8,8);
	\draw[thick,->](0,0) -- (8.25,0);
	\draw[thick,->](0,0) -- (0,8.25);
	\node at (8.5,0) {$q$};
	\node at (0,8.5) {$p$};
	\node at (-0.3,1) {$0$};
	\node at (-0.3,3) {$1$};
	\node at (-0.3,5) {$2$};
	\node at (-0.3,7) {$3$};
	\node at (1,-0.3) {$0$};
	\node at (3,-0.3) {$1$};
	\node at (5,-0.3) {$2$};
	\node at (7,-0.3) {$3$};
	
	\node at (1,1) {$\C$};
	
	\node at (3,1) {$\C$};
	
	\node at (1,3) {$\C$};
	
	\node at (5,1) {$\C$};
	
	\node at (2.5,2.5) {$\C^3$};
	\node at (3,3.50) {$\C^2$};
	\node at (3.50,3) {$\C^2$};
	\draw[thin,->] (3,3.75) -- (3,4.25);
	\draw[thin,->] (3.75,3) -- (4.25,3);
	
	\node at (1,5) {$\C$};

	\node at (7,1) {$\C$};
	
	\node at (5.5,2.5) {$\C^3$};
	\node at (4.5,3) {$\C^2$};
	\node at (5,3.50) {$\C^2$};
	\draw[thin,->] (5,3.75) -- (5,4.25);
	
	\node at (2.5,5.5) {$\C^3$};
	\node at (3,4.5) {$\C^2$};
	\node at (3.50,5) {$\C^2$};
	\draw[thin,->] (3.75,5) -- (4.25,5);
	
	\node at (1,7) {$\C$};
	
	\node at (7,3) {$\C$};
	
	\node at (5.5,5.5) {$\C^3$};
	\node at (5,4.50) {$\C^2$};
	\node at (4.5,5) {$\C^2$};
	
	\node at (3,7) {$\C$};
	
	\node at (7,5) {$\C$};
	
	\node at (5,7) {$\C$};
	
	\node at (7,7) {$\C$};

	\end{tikzpicture}}
	\caption{Graphic representation of the decomposition stated in Theorem \ref{decomposition} related to the solvmanifold $G_1/\Gamma$ with underlying Lie algebra $\g_1$ in case $ (ii) $.}
	\label{figure2}
\end{minipage}
\end{figure}

\begin{figure}[H]
	\begin{minipage}{.45\textwidth}
		\centering
		\captionsetup{width=\linewidth}
	\scalebox{.7}{\begin{tikzpicture}
	\draw [step=2cm] (0,0) grid (8,8);
	\draw[thick,->](0,0) -- (8.25,0);
	\draw[thick,->](0,0) -- (0,8.25);
	\node at (8.5,0) {$q$};
	\node at (0,8.5) {$p$};
	\node at (-0.3,1) {$0$};
	\node at (-0.3,3) {$1$};
	\node at (-0.3,5) {$2$};
	\node at (-0.3,7) {$3$};
	\node at (1,-0.3) {$0$};
	\node at (3,-0.3) {$1$};
	\node at (5,-0.3) {$2$};
	\node at (7,-0.3) {$3$};
	
	\node at (1,1) {$\C$};
	
	\node at (3,1) {$\C$};
	
	\node at (1,3) {$\C$};
	
	\node at (5,1) {$\C$};
	
	\node at (3,3) {$\C^3$};
	
	\node at (1,5) {$\C$};
	
	\node at (7,1) {$\C$};
	
	\node at (5,3) {$\C^3$};
	
	\node at (3,5) {$\C^3$};
	
	\node at (1,7) {$\C$};
	
	\node at (7,3) {$\C$};
	
	\node at (5,5) {$\C^3$};
	
	\node at (3,7) {$\C$};
	
	\node at (7,5) {$\C$};
	
	\node at (5,7) {$\C$};
	
	\node at (7,7) {$\C$};

	\end{tikzpicture}}
	\caption{Graphic representation of the decomposition stated in Theorem \ref{decomposition} related to the solvmanifold $G_1/\Gamma$ with underlying Lie algebra $\g_1$ in case $ (iii) $.}
	\label{figure3}
\end{minipage}
\end{figure}

\subsection{\for{toc}{Decompositions for $\g_2^\alpha$}\except{toc}{}}\textbf{Decompositions for $\g_2^\alpha$:} starting from $\alpha=0$, we have:

\begin{description}[itemsep=5pt]
	\item[$x_3=\frac{\pi}{2}$] \begin{description}[itemsep=5pt]
		\item[dots] $\bigoplus_{j=0}^3 D^{j,0} \oplus \bigoplus_{j=0}^3 D^{0,j}\oplus (D^{1,1})^{\oplus 5}\oplus (D^{2,1})^{\oplus 5}\oplus (D^{1,2})^{\oplus 5}\oplus (D^{2,2})^{\oplus 5}$;
		
		\item[horizontal lines] $\emptyset$;
		
		\item[vertical lines] $\emptyset$.
		
	\end{description}
	
	\item[$x_3\in\{\frac{\pi}{3}, \frac{\pi}{4}, \frac{\pi}{6}\}$] \begin{description}[itemsep=5pt]
		\item[dots] $\bigoplus_{j=0}^3 D^{j,0} \oplus \bigoplus_{j=0}^3 D^{0,j}\oplus (D^{1,1})^{\oplus 3}\oplus (D^{2,1})^{\oplus 3}\oplus (D^{1,2})^{\oplus 3}\oplus (D^{2,2})^{\oplus 3}$;
		
		\item[horizontal lines] $\emptyset$;
		
		\item[vertical lines] $\emptyset$.
		
	\end{description}

\end{description}

For $\alpha_n=\frac{1}{\pi}(|\log(\frac{n+\sqrt{n^2-4}}{2})|)$, we have:
\begin{description}[itemsep=5pt]
	\item[$b=\frac{2k+1}{2\mathfrak{Re}A_n}\pi$] \begin{description}[itemsep=5pt]
		\item[dots] $\bigoplus_{j=0}^3 D^{j,0} \oplus \bigoplus_{j=0}^3 D^{0,j}\oplus (D^{1,1})^{\oplus 3}\oplus (D^{2,1})^{\oplus 3}\oplus (D^{1,2})^{\oplus 3}\oplus (D^{2,2})^{\oplus 3}$;
		
		\item[horizontal lines] $(S_h^{1,0})^{\oplus 2}\oplus (S_h^{0,1})^{\oplus 2}\oplus (S_h^{1,1})^{\oplus 4}\oplus (S_h^{0,2})^{\oplus 2}\oplus (S_h^{2,1})^{\oplus 2}\oplus (S_h^{1,2})^{\oplus 4}\oplus (S_h^{2,2})^{\oplus 2}\oplus (S_h^{1,3})^{\oplus 2}$;
		
		\item[vertical lines] $(S_v^{1,0})^{\oplus 2}\oplus (S_v^{0,1})^{\oplus 2}\oplus (S_v^{2,0})^{\oplus 2}\oplus (S_v^{1,1})^{\oplus 4}\oplus (S_v^{2,1})^{\oplus 4}\oplus (S_v^{1,2})^{\oplus 2}\oplus (S_v^{3,1})^{\oplus 2}\oplus (S_v^{2,2})^{\oplus 2}$.
		
	\end{description}
	
	\item[$b=\frac{2k+1}{2\mathfrak{Re}A_n}\pi$] \begin{description}[itemsep=5pt]
		\item[dots] $\bigoplus_{j=0}^3 D^{j,0} \oplus \bigoplus_{j=0}^3 D^{0,j}\oplus (D^{1,1})^{\oplus 3}\oplus (D^{2,1})^{\oplus 3}\oplus (D^{1,2})^{\oplus 3}\oplus (D^{2,2})^{\oplus 3}$;
		
		\item[horizontal lines] $(S_h^{1,1})^{\oplus 2}\oplus (S_h^{1,2})^{\oplus 2}$;
		
		\item[vertical lines] $(S_v^{1,1})^{\oplus 2}\oplus (S_v^{2,1})^{\oplus 2}$.
		
	\end{description}
	
	\item[$b\neq\frac{k}{2\mathfrak{Re}A_n}\pi$] \begin{description}[itemsep=5pt]
		\item[dots] $\bigoplus_{j=0}^3 D^{j,0} \oplus \bigoplus_{j=0}^3 D^{0,j}\oplus (D^{1,1})^{\oplus 3}\oplus (D^{2,1})^{\oplus 3}\oplus (D^{1,2})^{\oplus 3}\oplus (D^{2,2})^{\oplus 3}$;
		
		\item[horizontal lines] $\emptyset$;
		
		\item[vertical lines] $\emptyset$.
		
	\end{description}
		
\end{description}

\begin{figure}[H]
	\begin{minipage}{.45\textwidth}
		\centering
		\captionsetup{width=\linewidth}
	\scalebox{.7}{\begin{tikzpicture}
	\draw [step=2cm] (0,0) grid (8,8);
	\draw[thick,->](0,0) -- (8.25,0);
	\draw[thick,->](0,0) -- (0,8.25);
	\node at (8.5,0) {$q$};
	\node at (0,8.5) {$p$};
	\node at (-0.3,1) {$0$};
	\node at (-0.3,3) {$1$};
	\node at (-0.3,5) {$2$};
	\node at (-0.3,7) {$3$};
	\node at (1,-0.3) {$0$};
	\node at (3,-0.3) {$1$};
	\node at (5,-0.3) {$2$};
	\node at (7,-0.3) {$3$};
	
	\node at (1,1) {$\C$};
	
	\node at (3,1) {$\C$};
	
	\node at (1,3) {$\C$};
	
	\node at (5,1) {$\C$};
	
	\node at (3,3) {$\C^5$};
	
	\node at (1,5) {$\C$};
	
	\node at (7,1) {$\C$};
	
	\node at (5,3) {$\C^5$};
	
	\node at (3,5) {$\C^5$};
	
	\node at (1,7) {$\C$};
	
	\node at (7,3) {$\C$};
	
	\node at (5,5) {$\C^5$};
	
	\node at (3,7) {$\C$};
	
	\node at (7,5) {$\C$};
	
	\node at (5,7) {$\C$};
	
	\node at (7,7) {$\C$};

	\end{tikzpicture}}
	\caption{Graphic representation of the decomposition stated in Theorem \ref{decomposition} related to the solvmanifold $G_2^0/\Gamma$ with underlying Lie algebra $\g_2^0$ and with $x_3=\frac{\pi}{2}$.}
	\label{figure4}
\end{minipage}
\hspace{1.2cm}
\begin{minipage}{.45\textwidth}
	\centering
	\captionsetup{width=\linewidth}
\scalebox{.7}{	\begin{tikzpicture}
	\draw [step=2cm] (0,0) grid (8,8);
	\draw[thick,->](0,0) -- (8.25,0);
	\draw[thick,->](0,0) -- (0,8.25);
	\node at (8.5,0) {$q$};
	\node at (0,8.5) {$p$};
	\node at (-0.3,1) {$0$};
	\node at (-0.3,3) {$1$};
	\node at (-0.3,5) {$2$};
	\node at (-0.3,7) {$3$};
	\node at (1,-0.3) {$0$};
	\node at (3,-0.3) {$1$};
	\node at (5,-0.3) {$2$};
	\node at (7,-0.3) {$3$};

	\node at (1,1) {$\C$};
	
	\node at (3,1) {$\C$};
	
	\node at (1,3) {$\C$};
	
	\node at (5,1) {$\C$};
	
	\node at (3,3) {$\C^3$};
	
	\node at (1,5) {$\C$};
	
	\node at (7,1) {$\C$};
	
	\node at (5,3) {$\C^3$};
	
	\node at (3,5) {$\C^3$};
	
	\node at (1,7) {$\C$};
	
	\node at (7,3) {$\C$};
	
	\node at (5,5) {$\C^3$};
	
	\node at (3,7) {$\C$};
	
	\node at (7,5) {$\C$};
	
	\node at (5,7) {$\C$};
	
	\node at (7,7) {$\C$};

	\end{tikzpicture}}
	\caption{Graphic representation of the decomposition stated in Theorem \ref{decomposition} related to the solvmanifold $G_2^0/\Gamma$ with underlying Lie algebra $\g_2^0$ and with $x_3\in\{\frac{\pi}{3},\frac{\pi}{4},\frac{\pi}{6}\}$.}
	\label{figure5}
\end{minipage}
\end{figure}

\begin{figure}[H]
	\begin{minipage}{.45\textwidth}
	\centering
	\captionsetup{width=\linewidth}
	\scalebox{.7}{\begin{tikzpicture}
	\draw [step=2cm] (0,0) grid (8,8);
	\draw[thick,->](0,0) -- (8.25,0);
	\draw[thick,->](0,0) -- (0,8.25);
	\node at (8.5,0) {$q$};
	\node at (0,8.5) {$p$};
	\node at (-0.3,1) {$0$};
	\node at (-0.3,3) {$1$};
	\node at (-0.3,5) {$2$};
	\node at (-0.3,7) {$3$};
	\node at (1,-0.3) {$0$};
	\node at (3,-0.3) {$1$};
	\node at (5,-0.3) {$2$};
	\node at (7,-0.3) {$3$};
	
	\node at (1,1) {$\C$};
	
	\node at (2.5,0.5) {$\C$};
	\node at (3.5,1) {$\C^2$};
	\node at (3,1.5) {$\C^2$};
	\draw[thin, ->] (3.75,1) -- (4.25,1);
	\draw[thin, ->] (3,1.75) -- (3,2.15);
	
	\node at (0.5,2.5) {$\C$};
	\node at (1,3.5) {$\C^2$};
	\node at (1.5,3) {$\C^2$};
	\draw[thin, ->] (1,3.75) -- (1,4.25);
	\draw[thin, ->] (1.75,3) -- (2.15,3);
	
	\node at (5.5,0.5) {$\C$};
	\node at (4.5,1) {$\C^2$};
	\node at (5,1.5) {$\C^2$};
	\draw[thin, ->] (5,1.75) -- (5,2.15);
	
	\node at (3,3) {$\C^3$};
	\node at (3,2.40) {$\C^2$};
	\node at (2.4,3) {$\C^2$};
	\node at (3,3.60) {$\C^4$};
	\node at (3.60,3) {$\C^4$};
	\draw[thin,->] (3,3.85) -- (3,4.15);
	\draw[thin,->] (3.85,3) -- (4.15,3);
	
	\node at (0.5,5.5) {$\C$};
	\node at (1,4.5) {$\C^2$};
	\node at (1.5,5) {$\C^2$};
	\draw[thin, ->] (1.75,5) -- (2.15,5);
	
	\node at (7,1) {$\C$};
	
	\node at (5,3) {$\C^3$};
	\node at (5,2.40) {$\C^2$};
	\node at (4.4,3) {$\C^4$};
	\node at (5,3.60) {$\C^4$};
	\node at (5.60,3) {$\C^2$};
	\draw[thin,->] (5,3.85) -- (5,4.15);
	\draw[thin,->] (5.85,3) -- (6.25,3);
	
	\node at (3,5) {$\C^3$};
	\node at (3,4.40) {$\C^4$};
	\node at (2.4,5) {$\C^2$};
	\node at (3,5.60) {$\C^2$};
	\node at (3.60,5) {$\C^4$};
	\draw[thin,->] (3,5.85) -- (3,6.25);
	\draw[thin,->] (3.85,5) -- (4.15,5);
	
	\node at (1,7) {$\C$};
	
	\node at (7.5,2.5) {$\C$};
	\node at (7,3.5) {$\C^2$};
	\node at (6.5,3) {$\C^2$};
	\draw[thin, ->] (7,3.75) -- (7,4.25);
	
	\node at (5,5) {$\C^3$};
	\node at (5,4.40) {$\C^4$};
	\node at (4.4,5) {$\C^4$};
	\node at (5,5.60) {$\C^2$};
	\node at (5.60,5) {$\C^2$};
	\draw[thin,->] (5,5.85) -- (5,6.25);
	\draw[thin,->] (5.85,5) -- (6.25,5);
	
	\node at (2.5,7.5) {$\C$};
	\node at (3.5,7) {$\C^2$};
	\node at (3,6.5) {$\C^2$};
	\draw[thin, ->] (3.75,7) -- (4.25,7);
	
	\node at (7.5,5.5) {$\C$};
	\node at (7,4.5) {$\C^2$};
	\node at (6.5,5) {$\C^2$};
	
	\node at (5.5,7.5) {$\C$};
	\node at (5,6.5) {$\C^2$};
	\node at (4.5,7) {$\C^2$};
	
	\node at (7,7) {$\C$};

	\end{tikzpicture}}
	\caption{Graphic representation of the decomposition stated in Theorem \ref{decomposition} related to the solvmanifold $G_2^{\alpha_n}/\Gamma$ with underlying Lie algebra $\g_2^{\alpha_n}, \alpha_n=\frac{1}{\pi}(|\log(\frac{n+\sqrt{n^2-4}}{2})|)$ and $b=\frac{2k+1}{2\mathfrak{Re}A_n}\pi$.}
	\label{figure6}
\end{minipage}
\hspace{1.2cm}
\begin{minipage}{.45\textwidth}
\centering	
\captionsetup{width=\linewidth}
\scalebox{.7}{\begin{tikzpicture}
	\draw [step=2cm] (0,0) grid (8,8);
	\draw[thick,->](0,0) -- (8.25,0);
	\draw[thick,->](0,0) -- (0,8.25);
	\node at (8.5,0) {$q$};
	\node at (0,8.5) {$p$};
	\node at (-0.3,1) {$0$};
	\node at (-0.3,3) {$1$};
	\node at (-0.3,5) {$2$};
	\node at (-0.3,7) {$3$};
	\node at (1,-0.3) {$0$};
	\node at (3,-0.3) {$1$};
	\node at (5,-0.3) {$2$};
	\node at (7,-0.3) {$3$};
	
	\node at (1,1) {$\C$};
	
	\node at (3,1) {$\C$};
	
	\node at (1,3) {$\C$};
	
	\node at (5,1) {$\C$};
	
	\node at (2.5,2.5) {$\C^3$};
	\node at (3,3.50) {$\C^2$};
	\node at (3.50,3) {$\C^2$};
	\draw[thin,->] (3,3.75) -- (3,4.25);
	\draw[thin,->] (3.75,3) -- (4.25,3);
	
	\node at (1,5) {$\C$};
	
	\node at (7,1) {$\C$};
	
	\node at (5.5,2.5) {$\C^3$};
	\node at (4.5,3) {$\C^2$};
	\node at (5,3.50) {$\C^2$};
	\draw[thin,->] (5,3.75) -- (5,4.25);
	
	\node at (2.5,5.5) {$\C^3$};
	\node at (3,4.5) {$\C^2$};
	\node at (3.50,5) {$\C^2$};
	\draw[thin,->] (3.75,5) -- (4.25,5);
	
	\node at (1,7) {$\C$};
	
	\node at (7,3) {$\C$};
	
	\node at (5.5,5.5) {$\C^3$};
	\node at (5,4.50) {$\C^2$};
	\node at (4.5,5) {$\C^2$};
	
	\node at (3,7) {$\C$};
	
	\node at (7,5) {$\C$};
	
	\node at (5,7) {$\C$};
	
	\node at (7,7) {$\C$};

	\end{tikzpicture}}
	\caption{Graphic representation of the decomposition stated in Theorem \ref{decomposition} related to the solvmanifold $G_2^{\alpha_n}/\Gamma$ with underlying Lie algebra $\g_2^{\alpha_n}, \alpha_n=\frac{1}{\pi}(|\log(\frac{n+\sqrt{n^2-4}}{2})|)$ and $b=\frac{2k}{2\mathfrak{Re}A_n}\pi$.}
	\label{figure7}
\end{minipage}
\vspace{2cm}
\end{figure}

\begin{figure}[H]
	\begin{minipage}{.45\textwidth}
		\captionsetup{width=\linewidth}
		\centering
	\scalebox{.7}{\begin{tikzpicture}
	\draw [step=2cm] (0,0) grid (8,8);
	\draw[thick,->](0,0) -- (8.25,0);
	\draw[thick,->](0,0) -- (0,8.25);
	\node at (8.5,0) {$q$};
	\node at (0,8.5) {$p$};
	\node at (-0.3,1) {$0$};
	\node at (-0.3,3) {$1$};
	\node at (-0.3,5) {$2$};
	\node at (-0.3,7) {$3$};
	\node at (1,-0.3) {$0$};
	\node at (3,-0.3) {$1$};
	\node at (5,-0.3) {$2$};
	\node at (7,-0.3) {$3$};
	
	\node at (1,1) {$\C$};
	
	\node at (3,1) {$\C$};
	
	\node at (1,3) {$\C$};
	
	\node at (5,1) {$\C$};
	
	\node at (3,3) {$\C^3$};
	
	\node at (1,5) {$\C$};
	
	\node at (7,1) {$\C$};
	
	\node at (5,3) {$\C^3$};
	
	\node at (3,5) {$\C^3$};
	
	\node at (1,7) {$\C$};
	
	\node at (7,3) {$\C$};
	
	\node at (5,5) {$\C^3$};
	
	\node at (3,7) {$\C$};
	
	\node at (7,5) {$\C$};
	
	\node at (5,7) {$\C$};
	
	\node at (7,7) {$\C$};

	\end{tikzpicture}}
	\caption{Graphic representation of the decomposition stated in Theorem \ref{decomposition} related to the solvmanifold $G_2^{\alpha_n}/\Gamma$ with underlying Lie algebra $\g_2^{\alpha_n}, \alpha_n=\frac{1}{\pi}(|\log(\frac{n+\sqrt{n^2-4}}{2})|)$ and $b\neq\frac{k}{2\mathfrak{Re}A_n}\pi$.}
	\label{figure8}
\end{minipage}
\end{figure}

\subsection{\for{toc}{Decompositions for $\g_8$}\except{toc}{}}\textbf{Decompositions for $\g_8$:} distinguishing the cases involved, we have:

\begin{description}[itemsep=5pt]
	\item[case \ref{caso1}] \begin{description}[itemsep=5pt]
		\item[dots] $\bigoplus_{i,j=0}^3 D^{i,j}$;
		
		\item[horizontal lines] $\emptyset$;
		
		\item[vertical lines] $\emptyset$.
		
	\end{description}
	
	\item[case \ref{caso2}] \begin{description}[itemsep=5pt]
		\item[dots] $\bigoplus_{i,j=0}^3 D^{i,j}$;
		
		\item[horizontal lines] $\bigoplus_{j=0}^3(S_h^{1,j})^{\oplus 2}$;
		
		\item[vertical lines] $\bigoplus_{j=0}^3(S_v^{j,1})^{\oplus 2}$.
		
	\end{description}
	
	\item[case \ref{caso3}] \begin{description}[itemsep=5pt]
		\item[dots] $\bigoplus_{i,j=0}^3 D^{i,j}$;
		
		\item[horizontal lines] $(S_h^{0,1})^{\oplus 2}\oplus (S_h^{0,2})^{\oplus 2}\oplus (S_h^{2,1})^{\oplus 2}\oplus (S_h^{2,2})^{\oplus 2}$;
		
		\item[vertical lines] $(S_v^{1,0})^{\oplus 2}\oplus (S_v^{2,0})^{\oplus 2}\oplus (S_v^{1,2})^{\oplus 2}\oplus (S_v^{2,2})^{\oplus 2}$.
		
	\end{description}
	
	\item[case \ref{caso4}] \begin{description}[itemsep=5pt]
		\item[dots] $\bigoplus_{j=0}^3 D^{j,0} \oplus \bigoplus_{j=0}^3 D^{0,j}\oplus (D^{1,1})^{\oplus 3}\oplus (D^{2,1})^{\oplus 3}\oplus (D^{1,2})^{\oplus 3}\oplus (D^{2,2})^{\oplus 3}$;
		
		\item[horizontal lines] $\emptyset$;
		
		\item[vertical lines] $\emptyset$.
		
	\end{description}
	
	\item[case \ref{caso5}] \begin{description}[itemsep=5pt]
		\item[dots] $\bigoplus_{i,j=0}^3 D^{i,j}$;
		
		\item[horizontal lines] $(S_h^{1,1})^{\oplus 2}\oplus (S_h^{1,2})^{\oplus 2}$;
		
		\item[vertical lines] $(S_v^{1,1})^{\oplus 2}\oplus (S_v^{2,1})^{\oplus 2}$.
		
	\end{description}
	
	\item[case \ref{caso6}] \begin{description}[itemsep=5pt]
		\item[dots] $\bigoplus_{j=0}^3 D^{j,0} \oplus \bigoplus_{j=0}^3 D^{0,j}\oplus (D^{1,1})^{\oplus 3}\oplus (D^{2,1})^{\oplus 3}\oplus (D^{1,2})^{\oplus 3}\oplus (D^{2,2})^{\oplus 3}$;
		
		\item[horizontal lines] $(S_h^{1,1})^{\oplus 2}\oplus (S_h^{1,2})^{\oplus 2}$;
		
		\item[vertical lines] $(S_v^{1,1})^{\oplus 2}\oplus (S_v^{2,1})^{\oplus 2}$.
		
	\end{description}
	
	\item[case \ref{caso7}] \begin{description}[itemsep=5pt]
		\item[dots] $\bigoplus_{j=0}^3 D^{j,0} \oplus \bigoplus_{j=0}^3 D^{0,j}\oplus (D^{1,1})^{\oplus 3}\oplus (D^{2,1})^{\oplus 3}\oplus (D^{1,2})^{\oplus 3}\oplus (D^{2,2})^{\oplus 3}$;
		
		\item[horizontal lines] $(S_h^{1,0})^{\oplus 2}\oplus (S_h^{0,1})^{\oplus 2}\oplus (S_h^{1,1})^{\oplus 4}\oplus (S_h^{0,2})^{\oplus 2}\oplus (S_h^{2,1})^{\oplus 2}\oplus (S_h^{1,2})^{\oplus 4}\oplus (S_h^{2,2})^{\oplus 2}\oplus (S_h^{1,3})^{\oplus 2}$;
		
		\item[vertical lines] $(S_v^{1,0})^{\oplus 2}\oplus (S_v^{0,1})^{\oplus 2}\oplus (S_v^{2,0})^{\oplus 2}\oplus (S_v^{1,1})^{\oplus 4}\oplus (S_v^{2,1})^{\oplus 4}\oplus (S_v^{1,2})^{\oplus 2}\oplus (S_v^{3,1})^{\oplus 2}\oplus (S_v^{2,2})^{\oplus 2}$.
		
	\end{description}
	
\end{description}

\begin{figure}[H]
	\begin{minipage}{.45\textwidth}
		\centering
		\captionsetup{width=\linewidth}
	\scalebox{.7}{\begin{tikzpicture}
	\draw [step=2cm] (0,0) grid (8,8);
	\draw[thick,->](0,0) -- (8.25,0);
	\draw[thick,->](0,0) -- (0,8.25);
	\node at (8.5,0) {$q$};
	\node at (0,8.5) {$p$};
	\node at (-0.3,1) {$0$};
	\node at (-0.3,3) {$1$};
	\node at (-0.3,5) {$2$};
	\node at (-0.3,7) {$3$};
	\node at (1,-0.3) {$0$};
	\node at (3,-0.3) {$1$};
	\node at (5,-0.3) {$2$};
	\node at (7,-0.3) {$3$};
	
	\node at (1,1) {$\C$};
	
	\node at (3,1) {$\C$};
	
	\node at (1,3) {$\C$};
	
	\node at (5,1) {$\C$};
	
	\node at (3,3) {$\C$};
	
	\node at (1,5) {$\C$};
	
	\node at (7,1) {$\C$};
	
	\node at (5,3) {$\C$};
	
	\node at (3,5) {$\C$};
	
	\node at (1,7) {$\C$};
	
	\node at (7,3) {$\C$};
	
	\node at (5,5) {$\C$};
	
	\node at (3,7) {$\C$};
	
	\node at (7,5) {$\C$};
	
	\node at (5,7) {$\C$};
	
	\node at (7,7) {$\C$};

	\end{tikzpicture}}
	\caption{Graphic representation of the decomposition stated in Theorem \ref{decomposition} related to the solvmanifold $G_8/\Gamma$ with underlying Lie algebra $\g_8$ in case \ref{caso1} of Theorem \ref{teorema}.}
	\label{figure9}
\end{minipage}
\hspace{1.2cm}
\begin{minipage}{.45\textwidth}
	\centering
	\captionsetup{width=\linewidth}	
\scalebox{.7}{\begin{tikzpicture}
	\draw [step=2cm] (0,0) grid (8,8);
	\draw[thick,->](0,0) -- (8.25,0);
	\draw[thick,->](0,0) -- (0,8.25);
	\node at (8.5,0) {$q$};
	\node at (0,8.5) {$p$};
	\node at (-0.3,1) {$0$};
	\node at (-0.3,3) {$1$};
	\node at (-0.3,5) {$2$};
	\node at (-0.3,7) {$3$};
	\node at (1,-0.3) {$0$};
	\node at (3,-0.3) {$1$};
	\node at (5,-0.3) {$2$};
	\node at (7,-0.3) {$3$};
	
	\node at (1,1) {$\C$};
	
	\node at (2.5,0.75) {$\C$};
	\node at (3.5,1.25) {$\C^2$};
	\draw[thin,->] (3.75,1.25) -- (4.25,1.25);
	
	\node at (0.75,2.5) {$\C$};
	\node at (1.25,3.5) {$\C^2$};
	\draw[thin,->] (1.25,3.75) -- (1.25,4.25);
	
	\node at (5.5,0.75) {$\C$};
	\node at (4.5,1.25) {$\C^2$};
	
	\node at (2.5,2.5) {$\C$};
	\node at (3,3.50) {$\C^2$};
	\node at (3.50,3) {$\C^2$};
	\draw[thin,->] (3,3.75) -- (3,4.25);
	\draw[thin,->] (3.75,3) -- (4.25,3);
	
	\node at (0.75,5.5) {$\C$};
	\node at (1.25,4.5) {$\C^2$};
	
	\node at (7,1) {$\C$};
	
	\node at (5.5,2.5) {$\C$};
	\node at (4.5,3) {$\C^2$};
	\node at (5,3.50) {$\C^2$};
	\draw[thin,->] (5,3.75) -- (5,4.25);
	
	\node at (2.5,5.5) {$\C$};
	\node at (3,4.5) {$\C^2$};
	\node at (3.50,5) {$\C^2$};
	\draw[thin,->] (3.75,5) -- (4.25,5);
	
	\node at (1,7) {$\C$};
	
	\node at (7.25,2.5) {$\C$};
	\node at (6.75,3.5) {$\C^2$};
	\draw[thin,->] (6.75,3.75) -- (6.75,4.25);
	
	\node at (5.5,5.5) {$\C$};
	\node at (5,4.50) {$\C^2$};
	\node at (4.5,5) {$\C^2$};
	
	\node at (2.5,7.25) {$\C$};
	\node at (3.5,6.75) {$\C^2$};
	\draw[thin,->] (3.75,6.75) -- (4.25,6.75);
	
	\node at (7.25,5.5) {$\C$};
	\node at (6.75,4.5) {$\C^2$};
	
	\node at (5.5,7.25) {$\C$};
	\node at (4.5,6.75) {$\C^2$};
	
	\node at (7,7) {$\C$};
	
	\end{tikzpicture}}
	\caption{Graphic representation of the decomposition stated in Theorem \ref{decomposition} related to the solvmanifold $G_8/\Gamma$ with underlying Lie algebra $\g_8$ in case \ref{caso2} of Theorem \ref{teorema}.}
	\label{figure10}
	\end{minipage}
\end{figure}

\begin{figure}[H]
	\begin{minipage}{.45\textwidth}
		\centering
		\captionsetup{width=\linewidth}
	\scalebox{.7}{\begin{tikzpicture}
	\draw [step=2cm] (0,0) grid (8,8);
	\draw[thick,->](0,0) -- (8.25,0);
	\draw[thick,->](0,0) -- (0,8.25);
	\node at (8.5,0) {$q$};
	\node at (0,8.5) {$p$};
	\node at (-0.3,1) {$0$};
	\node at (-0.3,3) {$1$};
	\node at (-0.3,5) {$2$};
	\node at (-0.3,7) {$3$};
	\node at (1,-0.3) {$0$};
	\node at (3,-0.3) {$1$};
	\node at (5,-0.3) {$2$};
	\node at (7,-0.3) {$3$};
	
	\node at (1,1) {$\C$};
	
	\node at (3.5,0.5) {$\C$};
	\node at (3,1.5) {$\C^2$};
	\draw[thin,->] (3,1.75) -- (3,2.25);
	
	\node at (0.5,3.5) {$\C$};
	\node at (1.5,3) {$\C^2$};
	\draw[thin,->] (1.75,3) -- (2.25,3);
	
	\node at (4.5,0.5) {$\C$};
	\node at (5,1.5) {$\C^2$};
	\draw[thin,->] (5,1.75) -- (5,2.25);
	
	\node at (3.5,3.5) {$\C$};
	\node at (3,2.50) {$\C^2$};
	\node at (2.50,3) {$\C^2$};

	\node at (0.5,4.5) {$\C$};
	\node at (1.5,5) {$\C^2$};
	\draw[thin,->] (1.75,5) -- (2.25,5);
	
	\node at (7,1) {$\C$};
	
	\node at (4.5,3.5) {$\C$};
	\node at (5,2.50) {$\C^2$};
	\node at (5.50,3) {$\C^2$};
	\draw[thin,->] (5.75,3) -- (6.25,3);
	
	\node at (3.5,4.5) {$\C$};
	\node at (2.50,5) {$\C^2$};
	\node at (3,5.50) {$\C^2$};
	\draw[thin,->] (3,5.75) -- (3,6.25);
	
	\node at (1,7) {$\C$};
	
	\node at (7.5,3.5) {$\C$};
	\node at (6.5,3) {$\C^2$};
	
	\node at (4.5,4.5) {$\C$};
	\node at (5.50,5) {$\C^2$};
	\node at (5,5.50) {$\C^2$};
	\draw[thin,->] (5,5.75) -- (5,6.25);
	\draw[thin,->] (5.75,5) -- (6.25,5);
		
	\node at (3.5,7.5) {$\C$};
	\node at (3,6.5) {$\C^2$};
	
	\node at (7.5,4.5) {$\C$};
	\node at (6.5,5) {$\C^2$};
	
	\node at (4.5,7.5) {$\C$};
	\node at (5,6.5) {$\C^2$};
	
	\node at (7,7) {$\C$};
	
	\end{tikzpicture}}
	\caption{Graphic representation of the decomposition stated in Theorem \ref{decomposition} related to the solvmanifold $G_8/\Gamma$ with underlying Lie algebra $\g_8$ in case \ref{caso3} of Theorem \ref{teorema}.}
	\label{figure11}
\end{minipage}
\hspace{1.2cm}
\begin{minipage}{.45\textwidth}
	\centering
	\captionsetup{width=\linewidth}
\scalebox{.7}{\begin{tikzpicture}
	\draw [step=2cm] (0,0) grid (8,8);
	\draw[thick,->](0,0) -- (8.25,0);
	\draw[thick,->](0,0) -- (0,8.25);
	\node at (8.5,0) {$q$};
	\node at (0,8.5) {$p$};
	\node at (-0.3,1) {$0$};
	\node at (-0.3,3) {$1$};
	\node at (-0.3,5) {$2$};
	\node at (-0.3,7) {$3$};
	\node at (1,-0.3) {$0$};
	\node at (3,-0.3) {$1$};
	\node at (5,-0.3) {$2$};
	\node at (7,-0.3) {$3$};
	
	\node at (1,1) {$\C$};
	
	\node at (3,1) {$\C$};
	
	\node at (1,3) {$\C$};
	
	\node at (5,1) {$\C$};
	
	\node at (3,3) {$\C^3$};
	
	\node at (1,5) {$\C$};
	
	\node at (7,1) {$\C$};
	
	\node at (5,3) {$\C^3$};
	
	\node at (3,5) {$\C^3$};
	
	\node at (1,7) {$\C$};
	
	\node at (7,3) {$\C$};
	
	\node at (5,5) {$\C^3$};
	
	\node at (3,7) {$\C$};
	
	\node at (7,5) {$\C$};
	
	\node at (5,7) {$\C$};
	
	\node at (7,7) {$\C$};

	\end{tikzpicture}}
	\caption{Graphic representation of the decomposition stated in Theorem \ref{decomposition} related to the solvmanifold $G_8/\Gamma$ with underlying Lie algebra $\g_8$ in case \ref{caso4} of Theorem \ref{teorema}.}
	\label{figure12}
\end{minipage}
\end{figure}

\begin{figure}[H]
	\begin{minipage}{.45\textwidth}
		\centering
		\captionsetup{width=\linewidth}
	\scalebox{.7}{\begin{tikzpicture}
	\draw [step=2cm] (0,0) grid (8,8);
	\draw[thick,->](0,0) -- (8.25,0);
	\draw[thick,->](0,0) -- (0,8.25);
	\node at (8.5,0) {$q$};
	\node at (0,8.5) {$p$};
	\node at (-0.3,1) {$0$};
	\node at (-0.3,3) {$1$};
	\node at (-0.3,5) {$2$};
	\node at (-0.3,7) {$3$};
	\node at (1,-0.3) {$0$};
	\node at (3,-0.3) {$1$};
	\node at (5,-0.3) {$2$};
	\node at (7,-0.3) {$3$};
	
	\node at (1,1) {$\C$};
	
	\node at (3,1) {$\C$};
	
	\node at (1,3) {$\C$};
	
	\node at (5,1) {$\C$};
		
	\node at (2.5,2.5) {$\C$};
	\node at (3,3.50) {$\C^2$};
	\node at (3.50,3) {$\C^2$};
	\draw[thin,->] (3,3.75) -- (3,4.25);
	\draw[thin,->] (3.75,3) -- (4.25,3);
	
	\node at (1,5) {$\C$};
	
	\node at (7,1) {$\C$};
	
	\node at (5.5,2.5) {$\C$};
	\node at (4.5,3) {$\C^2$};
	\node at (5,3.50) {$\C^2$};
	\draw[thin,->] (5,3.75) -- (5,4.25);
	
	\node at (2.5,5.5) {$\C$};
	\node at (3,4.5) {$\C^2$};
	\node at (3.50,5) {$\C^2$};
	\draw[thin,->] (3.75,5) -- (4.25,5);
	
	\node at (1,7) {$\C$};
	
	\node at (7,3) {$\C$};
	
	\node at (5.5,5.5) {$\C$};
	\node at (5,4.50) {$\C^2$};
	\node at (4.5,5) {$\C^2$};
	
	\node at (3,7) {$\C$};
	
	\node at (7,5) {$\C$};
	
	\node at (5,7) {$\C$};
	
	\node at (7,7) {$\C$};
	
	\end{tikzpicture}}
	\caption{Graphic representation of the decomposition stated in Theorem \ref{decomposition} related to the solvmanifold $G_8/\Gamma$ with underlying Lie algebra $\g_8$ in case \ref{caso5} of Theorem \ref{teorema}.}
	\label{figure13}
\end{minipage}
\hspace{1.2cm}
\begin{minipage}{.45\textwidth}
	\centering
	\captionsetup{width=\linewidth}
	\scalebox{.7}{\begin{tikzpicture}
	\draw [step=2cm] (0,0) grid (8,8);
	\draw[thick,->](0,0) -- (8.25,0);
	\draw[thick,->](0,0) -- (0,8.25);
	\node at (8.5,0) {$q$};
	\node at (0,8.5) {$p$};
	\node at (-0.3,1) {$0$};
	\node at (-0.3,3) {$1$};
	\node at (-0.3,5) {$2$};
	\node at (-0.3,7) {$3$};
	\node at (1,-0.3) {$0$};
	\node at (3,-0.3) {$1$};
	\node at (5,-0.3) {$2$};
	\node at (7,-0.3) {$3$};
	
	\node at (1,1) {$\C$};
	
	\node at (3,1) {$\C$};
	
	\node at (1,3) {$\C$};
	
	\node at (5,1) {$\C$};
	
	\node at (2.5,2.5) {$\C^3$};
	\node at (3,3.50) {$\C^2$};
	\node at (3.50,3) {$\C^2$};
	\draw[thin,->] (3,3.75) -- (3,4.25);
	\draw[thin,->] (3.75,3) -- (4.25,3);
	
	\node at (1,5) {$\C$};
	
	\node at (7,1) {$\C$};
	
	\node at (5.5,2.5) {$\C^3$};
	\node at (4.5,3) {$\C^2$};
	\node at (5,3.50) {$\C^2$};
	\draw[thin,->] (5,3.75) -- (5,4.25);
	
	\node at (2.5,5.5) {$\C^3$};
	\node at (3,4.5) {$\C^2$};
	\node at (3.50,5) {$\C^2$};
	\draw[thin,->] (3.75,5) -- (4.25,5);
	
	\node at (1,7) {$\C$};
	
	\node at (7,3) {$\C$};
	
	\node at (5.5,5.5) {$\C^3$};
	\node at (5,4.50) {$\C^2$};
	\node at (4.5,5) {$\C^2$};
	
	\node at (3,7) {$\C$};
	
	\node at (7,5) {$\C$};
	
	\node at (5,7) {$\C$};
	
	\node at (7,7) {$\C$};

	\end{tikzpicture}}
	\caption{Graphic representation of the decomposition stated in Theorem \ref{decomposition} related to the solvmanifold $G_8/\Gamma$ with underlying Lie algebra $\g_8$ in case \ref{caso6} of Theorem \ref{teorema}.}
	\label{figure14}
\end{minipage}
\end{figure}

\begin{figure}[H]
	\begin{minipage}{.45\textwidth}
		\centering
		\captionsetup{width=\linewidth}
	\scalebox{.7}{\begin{tikzpicture}
	\draw [step=2cm] (0,0) grid (8,8);
	\draw[thick,->](0,0) -- (8.25,0);
	\draw[thick,->](0,0) -- (0,8.25);
	\node at (8.5,0) {$q$};
	\node at (0,8.5) {$p$};
	\node at (-0.3,1) {$0$};
	\node at (-0.3,3) {$1$};
	\node at (-0.3,5) {$2$};
	\node at (-0.3,7) {$3$};
	\node at (1,-0.3) {$0$};
	\node at (3,-0.3) {$1$};
	\node at (5,-0.3) {$2$};
	\node at (7,-0.3) {$3$};
	
	\node at (1,1) {$\C$};
	
	\node at (2.5,0.5) {$\C$};
	\node at (3.5,1) {$\C^2$};
	\node at (3,1.5) {$\C^2$};
	\draw[thin, ->] (3.75,1) -- (4.25,1);
	\draw[thin, ->] (3,1.75) -- (3,2.15);
	
	\node at (0.5,2.5) {$\C$};
	\node at (1,3.5) {$\C^2$};
	\node at (1.5,3) {$\C^2$};
	\draw[thin, ->] (1,3.75) -- (1,4.25);
	\draw[thin, ->] (1.75,3) -- (2.15,3);
	
	\node at (5.5,0.5) {$\C$};
	\node at (4.5,1) {$\C^2$};
	\node at (5,1.5) {$\C^2$};
	\draw[thin, ->] (5,1.75) -- (5,2.15);
	
	\node at (3,3) {$\C^3$};
	\node at (3,2.40) {$\C^2$};
	\node at (2.4,3) {$\C^2$};
	\node at (3,3.60) {$\C^4$};
	\node at (3.60,3) {$\C^4$};
	\draw[thin,->] (3,3.85) -- (3,4.15);
	\draw[thin,->] (3.85,3) -- (4.15,3);
	
	\node at (0.5,5.5) {$\C$};
	\node at (1,4.5) {$\C^2$};
	\node at (1.5,5) {$\C^2$};
	\draw[thin, ->] (1.75,5) -- (2.15,5);
	
	\node at (7,1) {$\C$};
	
	\node at (5,3) {$\C^3$};
	\node at (5,2.40) {$\C^2$};
	\node at (4.4,3) {$\C^4$};
	\node at (5,3.60) {$\C^4$};
	\node at (5.60,3) {$\C^2$};
	\draw[thin,->] (5,3.85) -- (5,4.15);
	\draw[thin,->] (5.85,3) -- (6.25,3);
	
	\node at (3,5) {$\C^3$};
	\node at (3,4.40) {$\C^4$};
	\node at (2.4,5) {$\C^2$};
	\node at (3,5.60) {$\C^2$};
	\node at (3.60,5) {$\C^4$};
	\draw[thin,->] (3,5.85) -- (3,6.25);
	\draw[thin,->] (3.85,5) -- (4.15,5);
	
	\node at (1,7) {$\C$};
	
	\node at (7.5,2.5) {$\C$};
	\node at (7,3.5) {$\C^2$};
	\node at (6.5,3) {$\C^2$};
	\draw[thin, ->] (7,3.75) -- (7,4.25);
	
	\node at (5,5) {$\C^3$};
	\node at (5,4.40) {$\C^4$};
	\node at (4.4,5) {$\C^4$};
	\node at (5,5.60) {$\C^2$};
	\node at (5.60,5) {$\C^2$};
	\draw[thin,->] (5,5.85) -- (5,6.25);
	\draw[thin,->] (5.85,5) -- (6.25,5);
	
	\node at (2.5,7.5) {$\C$};
	\node at (3.5,7) {$\C^2$};
	\node at (3,6.5) {$\C^2$};
	\draw[thin, ->] (3.75,7) -- (4.25,7);
	
	\node at (7.5,5.5) {$\C$};
	\node at (7,4.5) {$\C^2$};
	\node at (6.5,5) {$\C^2$};
	
	\node at (5.5,7.5) {$\C$};
	\node at (5,6.5) {$\C^2$};
	\node at (4.5,7) {$\C^2$};
	
	\node at (7,7) {$\C$};

	\end{tikzpicture}}
	\caption{Graphic representation of the decomposition stated in Theorem \ref{decomposition} related to the solvmanifold $G_8/\Gamma$ with underlying Lie algebra $\g_8$ in case \ref{caso7} of Theorem \ref{teorema}.}
	\label{figure15}
\end{minipage}
\end{figure}

\bibliography{splitting_type_ARTICOLO}

\end{document}